\documentclass[12pt]{article}
\usepackage{amsmath, latexsym, amsthm, amssymb}
\usepackage{color}
\usepackage[normalem]{ulem}
\usepackage[toc,page]{appendix}
\usepackage{graphicx}

\usepackage[
colorlinks=true, urlcolor=blue, citecolor=blue, urlcolor=red, citecolor=red, pdfstartview=FitH]{hyperref}

\numberwithin{equation}{section}

\newtheorem{theorem}{Theorem}[section]
\newtheorem{lemma}[theorem]{Lemma}
\newtheorem{proposition}[theorem]{Proposition}
\newtheorem{corollary}{Corollary}[section]
\newtheorem{remark}{Remark}[section]
\newtheorem{assumption}{Assumption}[section]
\newtheorem{definition}{Definition}[section]
\newtheorem{example}{Example}[section]
\newcommand{\NN}{\mathbb{N}}
\newcommand{\m}{\mathbf{m}}
\newcommand{\rr}{\mathbb{R}}
\newcommand{\W}{W}
\newcommand{\U}{u}

\newcommand{\Halmos}{$\Box$}

\title{ Admission Control for Double-ended Queues}

\author{ Xin Liu\hspace{2mm} and Ananda Weerasinghe\footnote{Research partially supported by  Simons Foundation Collaboration Grant  317381.}\\
 Clemson University and Iowa State University}

\date{\today}
\begin{document}

\maketitle \pagestyle{plain}

\vspace{-0.4in}

\abstract {We consider a controlled double-ended queue consisting of two classes of customers, labeled sellers and buyers. The sellers and buyers arrive in a trading market according to two independent renewal processes. Whenever there is a  seller and buyer pair, they are matched and leave the system instantaneously. The matching follows first-come-first-match service discipline. Those customers who cannot be matched immediately need to wait in the designated queue, and they are assumed to be impatient with generally distributed patience times. The control problem is concerned with the trade-off between blocking and abandonment for each class and the interplay of statistical behaviors of the two classes, and its objective is to choose optimal queue-capacities (buffer lengths) for sellers and buyers to minimize an infinite horizon discounted linear cost functional which consists of holding costs and penalty costs for blocking and abandonment. 

When the arrival intensities of both customer classes tend to infinity in concert, we use  a  heavy traffic approximation to formulate an approximate diffusion control problem (DCP), and develop an optimal threshold policy for the DCP. Finally, we employ the DCP solution to establish an easy-to-implement asymptotically optimal threshold policy  for the original queueing control problem.}

\medskip

\medskip

{\bf Abbreviated Title:}  Controlled double-ended queues.



\section{Introduction}

 We consider a mathematical model of a matching platform that matches two classes of customers, labeled as buyers and sellers. The customers of each class arrive sequentially and wait in their respective queue to be matched with a customer of the other class. They are matched according to the order of  arrival which is known as the first-come-first-match service discipline. Once matched, a trade occurs and the pair leaves the system immediately. Both classes of customers are assumed to be  impatient and they leave the system without being matched, if their patience runs out. Due to instantaneous matching, there cannot be positive numbers of buyers and sellers simultaneously in the system. Such queueing models are known as double-ended queues (or matching queues) with impatient customers.

The arrivals of buyers and sellers are sequential and form two independent renewal processes and  their patience-times are represented by  two  independent IID (independent and identically distributed) sequences.  To avoid long queues, we introduce an admission control mechanism of controlling the \textit{queue-capacity} by blocking the incoming customers.  On one hand, when there is no blocking, long queues can occur, which leads to heavy customer abandonment. On the other hand, when too many arrivals are blocked,  profits of the operation decrease. Furthermore, blocking one side of the system affects the other side. In this work the queue-capacities for the two customer classes are the system manager's choice and at a given time,  they are represented by a vector of two positive integer-valued, time-dependent random variables, which may depend on the past history as well as the current state of the system. The two components of the queue-capacity vector behave as barriers on the seller queue and the buyer queue. Incoming customers of each class are blocked (i.e., rejected), when the queue is full at capacity at the time of their arrival. Once in the queue, they abandon if their patience expires, which happens after random times, IID across each customer class. We introduce  
an infinite horizon discounted linear cost functional which consists of holding costs, abandonment costs, and  blocking costs. 
We are interested in analyzing such a controlled system and to develop near optimal control policies which minimize the above described cost structure of the queueing control problem (QCP). However, this problem is too complex for direct analysis and therefore, we resort to a heavy traffic approximation. Approximating queueing systems in heavy traffic by Brownian models is an effective strategy in obtaining  both qualitative and quantitative insights of the original queueing systems (cf. \cite{amr, kocaga, reed, weera3}). Our techniques are closely related to the analysis of queueing systems in Halfin-Whitt heavy traffic regime where the number of servers tends to infinity (cf. \cite{dai, halfin, reed2, weera4}).
 
To establish a heavy traffic approximation, we develop an asymptotic framework, under which the class-dependent arrival intensities are increasing to infinity in concert so that the system is  critically loaded. To capture the behavior of both customer classes, we define the state of the system by an \textit{imbalance process}  whose value at time $t$ is given by the  number of sellers at time $t$ minus the number of buyers at  time $t$. We first establish the tightness of a sequence of diffusion-scaled state processes and identify the limit points of such a sequence as solutions to a stochastic differential equation (SDE). An associated diffusion control problem (DCP) is then formulated using such limiting SDEs. The DCP turns out to be a singular control problem, and its solution is obtained by finding a smooth solution to the corresponding Hamilton-Jacobi-Bellman (HJB) equation. Using the solution to the HJB equation, we obtain an optimal threshold policy for the DCP. Finally, we employ it to establish a threshold regime which describes four different types of asymptotically optimal strategies for the original QCP. Therefore, our solution yields a simple, asymptotically optimal control strategy which is easy to implement in a double-ended queue and the involved threshold values are easy to compute from the given system parameters. 

The threshold parameters are characterized as the ratios $T_b\equiv (c_b + r_b\delta_b)/(\alpha + \delta_b)$ and $T_s\equiv (c_s + r_s\delta_s)/(\alpha + \delta_s)$, where $c_b, r_b,\delta_b$ (resp. $c_s, r_s,\delta_s$) represent the holding cost rate per buyer (resp. seller), the abandonment cost per buyer (resp. seller), and the abandonment rate for buyers (resp. sellers), respectively, and the parameter $\alpha$ is the discount factor in the cost functional \eqref{cost-1}. Now let $p_b$ and $p_s$ denote the penalty costs per each blocked buyer and seller, respectively. The strategy we develop says when $p_b \ge T_b$ (blocking is expensive), there should be no blocking on buyers, and when $p_b < T_b$ (blocking is cheap), an asymptotically optimal buffer size for buyers can be designed using the solution of the DCP. The policy for the seller side is the same with $p_b, T_b$ replaced by $p_s, T_s.$ To gain an insight into our solution, let us consider a special Markovian setting where the double-ended queue has Poisson arrivals for both types of customers with exponential patience times. Considering the buyer side, the overall cost rate per buyer becomes $c_b+r_b \delta_b$. Our cost structure in \eqref{cost-1} can be thought of as the expected total cost of the three types of costs over an ``observation period $[0, \tau]$", where $\tau$ is an exponential random variable which is independent of all the other processes in the system, and has a rate parameter given as the discount factor $\alpha$ in \eqref{cost-1}. Hence the waiting time for a buyer who is facing an ``extremely long queue'' upon arrival is the minimum of two independent exponential random variables, with parameters  $\delta_b$ (for abandonment) and  $\alpha$ (for the observation period), respectively. Thus $T_b = {(c_b + r_b\delta_b)}/{(\alpha+\delta_b)}$ represents the expected cost for this buyer. Now if the blocking cost $p_b$ is greater than $T_b,$  it is reasonable to admit  this buyer to  the queue. On the other hand, if $p_b < T_b,$ then it is better to block the buyer. Similar explanation can be given to the seller side. {However, this intuition does not yield the values of optimal queue-capacities. In this work, relying on the DCP obtained by the heavy traffic approximation, we show that this intuition remains valid under more general assumptions for the arrival processes and patience time distributions, and develop an asymptotically optimal admission control policy.  It should be noted that the proposed asymptotically optimal queue-capacity for buyers depend on the statistical behavior of sellers and vice versa.}

This work is related to the authors' previous works \cite{weera4} and \cite{L}. The work \cite{L} develops the diffusion approximation of the uncontrolled double-ended queue with renewal arrivals and generally distributed patience times and establish a linear relationship between the diffusion-scaled queue length and the offered waiting time (the asymptotic Little's law). In \cite{weera4}, a controlled $G/M/n/B + GI$ queue is studied in the Halfin-Whitt heavy traffic (also known as Quality and Efficiency-Driven) regime, where $B$ is the control process representing the queue-capacity.  Both \cite{weera4} and our work study the trade-off between abandonment and blocking  and use the heavy traffic approximation to find near-optimal cost minimization strategies. 
However, a significant challenge faced in this work will be the dependence of the random quantities of the buyer queue (resp. seller queue) such as the queue-length, virtual waiting times etc, on the statistical behaviors of the quantities of the sellers (resp. buyers). This leads to develop original proofs such as Proposition  \ref{prop-3.3}.   We summarize the novelty of the current work as follows. (i) We establish the tightness and the moment bounds for the diffusion-scaled state processes under any admissible control satisfying Assumption \ref{assump:admiss} (see Theorem \ref{thm3.2}), while \cite{L} only studies the uncontrolled double-ended queue and \cite{weera4} establishes the tightness result under constant control policies; (ii) To establish Theorem \ref{thm3.2}, we introduce the virtual waiting time processes for the unblocked buyers and sellers and establish the asymptotic Little's law in Proposition \ref{prop-3.9}; (iii) The proofs of both Theorem \ref{thm3.2} and the asymptotic optimality result in Theorem \ref{thm:AO} depend on the properties of the two-sided Skorokhod problem (SP) over a time varying interval. This is in contrast with \cite{weera4} where only the one-sided SP is used. We establish new weak convergence and oscillation results for the two-sided SP with time-dependent barriers, which will be of independent interest (see Appendix \ref{sec:sm}); (iv) The study of the HJB equation of the DCP is much more involved and it leads to a free boundary problem. In the most interesting case, this free boundary includes two points $a^*<0< b^*,$ and the values  $|a^*|$ and $b^*$ describe the optimal boundaries for the DCP. Obtaining these free boundary points is quite complex and one has to carefully analyze the solution profiles of the corresponding differential equations on $(-\infty, 0) $ and $(0, \infty)$. Then use ``the principle of smooth fit'' at the origin  to find the smooth solution for the HJB equation.  In comparison, the free boundary associated with the HJB equation in \cite{weera4}  consists of a single point and proving its existence is relatively easy.

Many articles in the literature of double-ended queues are driven by applications. In \cite{kashyap}, a taxi queueing system with a limited waiting space is modeled as a double-ended queue with Markovian assumptions on the arrival processes of taxis and passengers, and the steady state behavior of the system is studied. In  \cite{conolly}, the effect of impatience behavior of customers in double-ended queues with Markovian arrivals and exponential patience times was studied. Double-ended queues have also been used to study perishable production-inventory systems (cf. \cite{perry0, perry}), where one side of the queue represents the inventory of products and the other side accepts the arrivals of orders. In a recent work \cite{lllz}, a production rate control problem is studied to minimize a finite horizon cost functional consisting of linear costs for inventory and waiting and a cost that penalizes rapid fluctuations of production rates. An asymptotic optimal production rate is developed under the fluid scaling given that the demand arrival rate is time and state dependent. 
Double-ended queues are the simplest matching systems and can be naturally generalized to have more than two classes of customers. The generalized multi-class matching queues have occured in assembled products in manufacturing setting where each product is completed by combining multiple components upon their arrivals (cf. \cite{harrison1, plumbeck, GW}). In his early work, Harrison  \cite{harrison1} studied the behavior of vector waiting times in  a model for an assembly line product made with several components in heavy traffic. Each input component  arrives according to an independent renewal process. Once the server has one component from each category, it takes a random processing time to finish the product. In \cite{GW}, such a matching system is studied with instantaneous processing, for which the authors consider the problem of minimizing finite horizon cumulative holding costs. A myopic discrete-review matching control is developed and shown to be asymptotically optimal in heavy traffic. Plumbeck and Ward \cite{plumbeck} study a control problem of an assemble-to-order system to maximize an expected infinite-horizon discounted profit by choosing product prices, component production capacities, and a dynamic policy for sequencing customer orders for assembly.    
   If both sides of the double-ended queue are generalized to have multiple classes, we end up with a bipartite matching system. Such systems are widely used to model the organ transplant systems, where one side represents patients with multiple classes and the other side represents organs of multiple types. We refer to \cite{aaa, kl} for optimal allocation problems of the bipartite matching systems in fluid scaling. 
   An application of the simple double-ended queue to organ transplant systems is studied in \cite{boxma}.


The rest of this article is organized as follows: In Section \ref{sec:model}, we introduce the model, the blocking control structure,  and the control problem for the double-ended queue. Section \ref{sec:weak} is devoted to the moment bounds and C-tightness of the diffusion-scaled processes. In particular, in Theorem \ref{thm3.2}, we identify the limit points of any given sequence of state processes as a solution to a SDE. In Section \ref{sec:dcp}, the DCP is formulated using such limiting SDE, and its explicit solution is established in Theorem \ref{thm:solution-DCP}. Section \ref{sec:as} establishes the asymptotic optimality. Theorem \ref{thm5.4} shows that the DCP provides a lower bound for the diffusion-scaled QCP. Next in Theorem \ref{thm5.2}, we employ the solution to the DCP to obtain asymptotically optimal policy under four different parameter regimes.  We establish new convergence and oscillation results for the two-sided SP in Appendix \ref{sec:sm}, which are of independent interest. Finally, Appendix \ref{sec:proofprop-R} collects a series of lemmas to find a solution to the HJB equation.

\paragraph{Notation.} Let $\mathbb{N}$ denote the set of positive integers and $\mathbb{R}$ denote the one dimensional Euclidean space.  For a function $f:G\rightarrow \mathbb{R}$ where $G$ is an open set,  we write $f\in C^k(G)$ if the $k^{th}$ derivative  of  $f$ is continuous on $G.$ For $0< T\le \infty$,   we denote the function space of $\mathbb{R}$-valued right-continuous functions with left limits (RCLL), defined  on $[0, T]$, by ${D}[0, T].$ This function space is endowed with  the standard Skorokhod $J_1$ topology. For  $k \in \mathbb{N},$   let $D^k[0, T]$ be the product of $k$ of $D[0, T]$ spaces.  The uniform norm on $[0, T]$  for a stochastic process  $X$ in $D^k[0, T]$ is defined by $\|X\|_T =\sup_{ 0 \leq t \leq T} |X(t)|.$  To describe the processes associated with the   $n^{\rm th}$  system,  we typically use the superscript, such as in the case of $X^n(\cdot),$ etc. Throughout, we use $\Rightarrow$ to denote weak convergence of  processes in ${D}^k[0, T]$. For each $f\in {D}^k[0, T]$, we let the oscillation of $f$ in a  sub-interval $[t_1, t_2]$ be defined by  
\begin{equation}
Osc(f, [t_1, t_2])\equiv \sup\{ |f(t)-f(s)|;  s, t \in [t_1, t_2] \}.
\label{Osc1}
\end{equation}
We simply denote $Osc(f, [0, T])$ by $Osc(f, T).$
For a given $\delta>0,$ its modulus of continuity $\omega(f, \delta, T)$ is defined by
\begin{equation}
\omega(f, \delta, T)= \sup\{ |f(t)-f(s)|:  |t-s| <\delta \text{  and }  s, t \in[0,  T] \}.
\label{omega}
\end{equation}
Using the modulus of continuity  $\omega( \cdot, \delta, T)$ is advantageous in our arguments, because of its sub-additive property: For $f, g \in {D}^k[0, T]$ and $\delta > 0, $ 
$$ \omega(f+g, \delta, T) \leq  \omega(f, \delta, T) + \omega(g, \delta, T). $$  
This also helps us  to establish the C-tightness of  several processes considered here. We also follow the convention that the infimum of an empty set is infinity. For any real number $x$,  $x^+=\max \{ 0,  x \} $ and  $x^-=\max \{ 0,  -x \}.$  For any two real numbers $a$ and $b$, $a \wedge b = \min \{ a, b \}$ and $a\vee b = \max\{a, b\}.$

\section{Double-ended Queues}\label{sec:model}

All our stochastic processes and random variables are defined on  a complete probability
space $ ( \Omega ,\mathfrak F, P) $. We have a sequence of
queueing systems indexed by $n\in\mathbb{N}$,  where the scaling parameter $n$ is used to model the scale and traffic intensity of the system (see the heavy traffic condition in Assumption \ref{assump:htc}).  The
 $n^{\rm th}$  system  represents  a simple trading market  where buyers and sellers arrive   according to two  independent renewal  processes  $A^n_b(\cdot)$   and   $A^n_s(\cdot).$  
 A trade occurs when a buyer meets a seller and thereafter the pair leaves the system instantaneously. The buyers and sellers are matched according to  first-come-first-match policy and  they wait in their respective queues if not getting matched immediately. Since the matching is instantaneous,   it is not possible to have  positive numbers of  buyers and sellers waiting in their queues simultaneously. It is assumed that  both buyers and sellers are impatient and  if they have to wait in the queue, they abandon the system when their patience expires.  This abandonment mechanism works as follows: with each customer, there is an associated
clock. This clock rings after a random time, and
if the clock rings while the customer is waiting in the queue, then the
customer abandons the system.  These clocks are all IID and independent of the arrival
 processes, as well as the history of the system up to that time. The cumulative distribution function of the patience-time for buyers is represented by $F_b$ and for  sellers, it is given by $F_s.$
At  any time instant, either  the buyer queue or the seller queue is empty and therefore, the state description of the system can be  given by the \emph{imbalance process} $ {X}^n.$ At any time instant $t \geq 0,$  if $X^n(t) \geq 0,$  then there are no buyers at time $t$ and $X^n(t)$ represents the queue length of sellers. Similarly, if $X^n(t) <0$ then there are no sellers at time $t$ and   $-X^n(t)$ represents the queue length of buyers.  Without loss of generality, we assume that the initial number of  sellers is given by  $X^n(0-) \geq 0.$ We refer to  \cite{LGK, L} and  \cite{GW} for similar representations of the state process.  Let the quantities $G^n_s(t)$ and $G^n_b(t)$  represent the numbers of sellers and buyers abandoning the system during $[0,  t]$, respectively.   

 To minimize the costs associated with abandonment and waiting, the management is permitted to control the system  by blocking new arrivals whenever the queue-length is sufficiently large. Therefore, we allow the system manager to choose  \emph{queue-capacities or buffer lengths} for buyers and sellers. Throughout, we use these two words intermittently.
 For each customer class, when the buffer  is full, incoming customers will be rejected. Each rejected customer incurs a
 loss in {\color{red}{profit.}} For the $n^{\rm th}$ system, the vector valued stochastic
 process   ${\bf m}_n(\cdot)=(m^n_b(\cdot),   m^n_s(\cdot))$ represents the processes of controlled queue-capacities (waiting-room sizes or buffer lengths), where $m^n_b(t) <0 <m^n_s(t)$ for each $t\ge 0$.  More precisely,  $-m^n_b(t)$  represents the buffer length for buyers at time $t$, and   $m^n_s(t)$  represents  the buffer length for sellers at time $t.$
  This queue-capacity process
${\bf m}_n(\cdot)$ is the only control at the disposal of the manager of the
$n^{\rm th}$ system. The choice of  large queue-capacities reduces
the blocking of customers and increases the profit margins. However, at the same
time, such large capacities are likely to give rise to long queues,
which will in turn, increase the number of abandonments from the system,
leading to a loss of income. This trade-off, between blocking and
abandonment, naturally leads  to a cost minimization  problem, which is the underlying theme of our paper.

Associated with a queue-capacity process ${\bf m}_n(\cdot)=(m^n_b(\cdot),   m^n_s(\cdot)),$  we introduce a pair of non-decreasing processes  $ U^n_s$ and $U^n_b.$   For any time $t \geq 0,$  $U^n_s(t)$ and $U^n_b(t)$ represent the numbers of  sellers and buyers rejected during $[0, t]$, respectively.  To describe the dynamics of the controlled system, we assume the stochastic primitives  $A^n_b,  A^n_s,  G^n_b$ and $G^n_s$ to be RCLL processes, and make the following assumptions about the model.  

\begin{assumption}[Initial condition]\label{assump:initial}


The  number of  initial customers $X^n(0) $ is assumed  to be deterministic, non-abandoning, and for some real $x$,  
\begin{equation}
\lim_{n \rightarrow \infty} \frac{X^n(0)}{\sqrt{n}} =x.
\label{initial}
\end{equation}

\end{assumption}

Throughout, we simply assume  $X^n(0)\geq 0$  so that there are no buyers initially in the system.  The non-abandonment of initial customers is not a restrictive assumption  and it can be easily relaxed following  the proof of   Lemma 4.1 of \cite{weera3}.  

\begin{assumption}[Arrival processes]\label{assump:arrival} 
We assume that the arrival processes  $A^n_b$ and  $A^n_s$ are independent renewal processes.  More precisely, there exist two positive sequences of real numbers  $\{\lambda^n_b\}_{n\in\NN}$ and  $\{\lambda^n_s\}_{n\in\NN}$ and positive constants  $\varsigma_b$ and $\varsigma_s$ so that for $n\in\NN$,  the inter-arrival times of buyers and sellers in the $n^{\text{th}}$ are independent IID sequences with mean-variance pairs $({1}/{\lambda^n_b}, ({\varsigma_b}/{\lambda^n_b})^2)$ and $({1}/{\lambda^n_s}, ({\varsigma_s}/{\lambda^n_s})^2)$, respectively. 
  
  \end{assumption}
  
  \begin{assumption}[Heavy traffic conditions]\label{assump:htc}

There exists a constant $\lambda_0 >0$ and $\beta_b, \beta_s\in \mathbb{R}$ so that 
\begin{align}
\lim_{n \rightarrow \infty} \frac{\lambda^n_b -\lambda_0 n}{\sqrt{n}} =\beta_b,
\label{HT-b}\\
\lim_{n \rightarrow \infty} \frac{\lambda^n_s -\lambda_0 n}{\sqrt{n}} =\beta_s.
\label{HT-s}
\end{align}

\end{assumption}
  

From the well known functional central limit theorem for renewal processes, we have  
\begin{equation}
\left(\frac{A^n_b(t)-\lambda^n_b t}{\sqrt{n}},   \frac{A^n_s(t)-\lambda^n_s t}{\sqrt{n}}\right)  \text{ converges weakly to } (\sigma_b B_1(t),  \sigma_s B_2(t))
\label{arrivals}
\end{equation}
in the  space ${D}^2[0, T]$, where $\sigma_b = \sqrt{\varsigma_b^2\lambda_0}$  and $\sigma_s = \sqrt{\varsigma_s^2\lambda_0}$ and $(B_1,  B_2)$ is a standard two dimensional Brownian motion. We also have the following  moment condition: For $T>0$,
\begin{equation}
E\left[\sup\limits_{t\in [0, T]} \left( \big({A^n_b(t)-\lambda^n_b t} \big)^2 +  \big( {A^n_s(t)-\lambda^n_s t}\big)^2\right)\right]  \leq C n(1+T^m),
\label{arrivals-M}
\end{equation}
where $C>0$ and the integer $m>1$ are constants independent of $T$ (for details, we refer to Lemma 2 of \cite{amr} and \cite{taksar}).

\begin{assumption}[Patience time distributions]


The patience-times of buyers and sellers are independent of each other. They are IID with cumulative distribution functions  $F_b$ and $F_s$, respectively.  We further assume  $F_b(0)=F_s(0)=0$ and they are right-differentiable at the origin with positive derivatives. Thus for some positive constants $\delta_b$ and $\delta_s$,
\begin{align}
\lim\limits_{h \rightarrow 0^+} \frac{F_b(h)}{h} =\delta_b,
\label{Patience-b}\\
\lim\limits_{h \rightarrow 0^+} \frac{F_s(h)}{h} =\delta_s.
\label{Patience-s}
\end{align}

\end{assumption}
Similar assumptions on patience-times of customers in many-server queues were imposed in the articles \cite{dai,  amr, weera4}.


\begin{assumption}[Admission control]\label{assump:admiss}

 At a given time $t \geq 0,$ the controlled queue-capacity of sellers is represented by $m^n_s(t)\geq 1$   and  the queue-capacity of buyers is represented by $-m^n_b(t)\geq 1.$ {Thus,  if a customer queue is empty,  then an  incoming arrival from the same  customer  class is always admitted.}
We assume that both $m^n_s(t)$ and $m^n_b(t) $ are integer valued, and the controlled queue-capacity process  ${\bf m}_n(\cdot)=(m^n_b(\cdot),   m^n_s(\cdot))$ has paths which are RCLL with a  finite number of jumps in each finite interval.
  Furthermore, the  control ${ \bf m}_n(t)$   is allowed to depend on the current state,
 as well as the whole history of the system up to time $t$. Therefore, 
 the process ${ \bf m}_n(t)$ is assumed to be  adapted to  $\mathcal{F}^n_{t}$ for each $t \geq 0$, where 
\begin{equation}
\mathcal{F}^n_t=\sigma (X^n(u), A^n_b(u),  A^n_s(u),  G^n_b(u),  G^n_s(u), U^n_s(u-), U^n_b(u-); 0 \leq u \le t ), 
\label{filtration}
\end{equation}
completed by all the null sets.  This $\sigma$-algebra
represents all the information available to the system manager at
time $t$. 

The controller is allowed to remove the initial customers if necessary. Therefore, the initial customer population is represented by $X^n(0-)$ and after the initial removal, the  customer population  at time $t=0$  is represented by $X^n(0).$ We assume  the process ${ \bf m}_n $ also adheres to the following conditions:
\begin{itemize}

\item[\rm (i)]  There is a  constant  $M >0$ independent of $n,$ but which depends on the initial data $x$ in \eqref{initial}  so that for each $ n,$    $-M \sqrt{n} < m^n_b(0)<0 < m^n_s(0) < M \sqrt{n}$  and it is assumed that   $ X^n(0) \in [m^n_b(0), m^n_s(0)].$
Thus the number of initially removed customers is given by $\max\{X^n(0-)-m^n_s(0), m^n_b(0)-X^n(0-)\}.$


\item[\rm (ii)]  At any time $t>0,$  $X^n(t) \in [m^n_b(t), m^n_s(t)].$
Thus, once allowed to enter the system, no customer will be removed from the queue in the future.

\item[\rm (iii)] The process ${ \bf m}_n$ 
satisfies 
\begin{equation}
\lim\limits_{\delta \rightarrow 0^+} \limsup\limits_{n \rightarrow \infty}\frac{1}{\sqrt{n}}\omega({ \bf m}_n,  \delta, T) = 0      \text{ in probability, }
\label{cap-3}
\end{equation}
\end{itemize}
where  the modulus of continuity $\omega$  is defined in \eqref{omega}.
 Loosely speaking, \eqref{cap-3} imposes that a change
of queue-capacities of sellers and buyers   at any time $t$ will be at most of order
$\sqrt{n}$.  If the controlled queue-capacity
is  a deterministic time-dependent function  ${ \bf m}_n(\cdot)$, then
 (\ref{cap-3}) can be replaced by the following
simple sufficient condition: For all $s, t$ in $[0, T]$,  
$$  |{ \bf m}_n(t)- { \bf m}_n(s) | \leq \sqrt{n}\hspace{1mm}p(T)[\rho(|t-s|) +f(n)].$$  
 Here $p: [0,\infty)\to \mathbb{R}_+$ is a positive
continuous  function; $\rho: [0,\infty)\to \mathbb{R}_+$ is a positive bounded  continuous
function, which satisfies $
\lim_{r\rightarrow 0}\rho(r)=0$; the function $f: \mathbb{N}\to  \mathbb{R}_+$ satisfies $\lim_{n\rightarrow\infty}f(n)=0$.

\end{assumption}

 Given such an   ${ \bf m}_n ,$   the solution to the SP with time dependent barriers described in \cite{burdzy}  guarantees the existence of the state process $X^n$ in $D[0, \infty).$
Note that the above assumptions allow constant ${ \bf m}_n$ policies. 
They  also accommodate the situation where  no buyers or sellers are ever rejected. This is typically associated with  the infinite buffer capacity.  However,  this can also be achieved  by simply choosing the finite buffer capacities of $m^n_b(t)=\min \{ \inf_{u\in [0, t]} X^n(u) -2, -1\}$ and    $m^n_s(t)= \max \{\sup_{u\in [0, t]} X^n(u) +2, 1\}$, where $\{X^n(t); t\ge 0\}$ is the queue length process for the uncontrolled $n^{th}$ system.

For a given buffer length policy  ${ \bf m}_n(\cdot)$, the controlled state process $X^n$ has integer-valued RCLL paths and it satisfies the following equation: For all $t \geq 0$,
\begin{equation*}
X^n(t)=X^n(0) + A^n_s(t) - A^n_b(t) -G^n_s(t) + G^n_b(t) -U^n_s(t) +U^n_b(t),
\end{equation*}
where the processes $U^n_b$ and $U^n_s$  are given by 
\begin{align}
U^n_b(t)&=\int_{0}^{t}1_{\{X^n(u)=m^n_b(u)\}}dA^n_b(u)
 + [m^n_b(0)-X^n(0-)]^+,  
 \label{eq2.8b}\\
U^n_s(t)&=\int_{0}^{t}1_{\{X^n(u)=m^n_s(u)\}}dA^n_s(u)
+ [X^n(0-)-m^n_s(0)]^+.  
 \label{eq2.8s}
\end{align}
The associated infinite-horizon discounted cost functional is defined as 
\begin{align*}
 \tilde{J}^n( X^n(0),U^n_s,U^n_b) & =E\biggm(\int_{0}^{\infty} e^{- \alpha t}[\tilde C(X^n(t))dt +r_s dG^n(t) +r_b dG^n_b(t) \\
 & \qquad +p_s dU^n_s(t) +p_b dU^n_b(t)]\biggm),
\end{align*}
where  $\tilde C(x)=c_s x^+ +c_b x^- $ for all $x\in\mathbb{R},$  and $\alpha$,  $ c_s, c_b, r_s, r_b, p_s, p_b$ are all positive constants. {In particular, $c_s$ and $c_b$ are the linear holding cost rates per each waiting seller and buyer, $r_s, r_b$ are linear penalty costs for each abandoning seller and buyer, $ p_s, p_b$ are linear penalty costs for each blocked seller and buyer, and finally, $\alpha$ is the discount factor.} 
 Our objective is to find optimal strategies which minimize the above cost functional and are  easy to implement from the design point of view.  
 
  We introduce the following fluid and diffusion scaled quantities, which are similar to those in Halfin-Whitt  heavy traffic regime. \\
  {\bf Fluid-scaled processes:} For a process $Y^n\in\{A^n_s, A^n_b, G^n_b, G^n_s, U^n_b, U^n_s, X^n\}$, its fluid-scaled version is defined to be $\bar Y^n(t) = Y^n(t)/n$ for each $t\ge 0.$ \\
  {\bf Diffusion-scaled processes:} The diffusion-scaled renewal processes are given as $\hat{A}^n_s(t)= (A^n_s(t)-\lambda_0 nt)/{\sqrt{n}}$, and $\hat{A}^n_b(t)= {(A^n_b(t)-\lambda_0 n t)}/{\sqrt{n}},$ and for any other process $Z^n \in \{G^n_b, G^n_s, U^n_b, U^n_s, X^n, {\bf m}_n\}$, its diffusion-scaled version is given as $\hat Z^n(t) = Z^n(t)/\sqrt{n}$ for $t\ge 0.$

The diffusion-scaled state process can now be formulated as 
\begin{equation}
\hat{X}^n(t)=\hat{X}^n(0) +\hat{ A}^n_s(t) - \hat{A}^n_b(t) -\hat{G}^n_s(t) + \hat{G}^n_b(t) -\hat{U}^n_s(t) +\hat{U}^n_b(t), \ \ t\ge 0.
\label{state-2}
\end{equation}
The corresponding diffusion-scaled cost function is given by $\hat{J}^n = \tilde{J}^n/\sqrt{n}$. Since the holding cost function is  piecewise linear, the diffusion-scaled cost functional can be formulated as a functional of the diffusion-scaled processes as follows:
\begin{equation}\label{cost-1}
\begin{aligned}
\hat{J}^n( \hat{X}^n(0), \hat{U}^n_s, \hat{U}^n_b) & = E \biggm(\int_{0}^{\infty} e^{- \alpha t} [\tilde C(\hat{X}^n(t))dt +r_s d\hat{G}^n_s(t) +r_b d\hat{G}^n_b(t)  \\
& \qquad +p_s d\hat{U}^n_s(t) +p_b d\hat{U}^n_b(t)]\biggm). 
\end{aligned}
\end{equation}
 The corresponding value function is given by 
 \begin{equation}
 \hat{V}^n(x)=\inf\limits_{\mathcal{A}^n_x} \hat{J}^n(x, \hat{U}^n_s, \hat{U}^n_b),
 \label{2.20}
 \end{equation}
where $\mathcal{A}^n_x $ is the collection of all admissible processes $(\hat{X}^n, \hat{U}^n_s, \hat{U}^n_b)$ with $\hat{X}^n(0)=x$, and the process $(\hat X^n, \hat{U}^n_s, \hat{U}^n_b)$ is said to be \textit{admissible} if the corresponding admission control $\m_n$ satisfies Assumption \ref{assump:admiss}.   

In the next Section \ref{sec:weak}, we establish the tightness of the diffusion-scaled processes. In particular, the computations illustrate that $\{\hat{X}_n\}_{n\ge 1}$ is stochastically bounded and hence $\tilde{J}^n(X^n(0), U^n_s, U^n_b)$ is of order $\sqrt{n}$ and $ \hat{J}^n( \hat{X}^n(0), \hat{U}^n_s, \hat{U}^n_b)$ is of order $1$. Next in Section \ref{sec:dcp}, we develop the DCP and derive an explicit optimal solution. Finally, in Section \ref{sec:as}, the DCP is shown to be a good approximation for the $n^{\rm th}$ QCP. In particular, we show that  $\hat{V}^n$, the value function of the $n^{\rm th}$ QCP, converges to the value function of the DCP as $n$ tends to infinity. We then propose a threshold type admission control policy described by the optimal 
solution of the DCP for the $n^{\rm th}$ QCP, and prove that it is asymptotically optimal for the QCP as $n$ tends to infinity. Here a sequence of admissible control policies $\{(\hat U^{n,*}_s, \hat U^{n,*}_b)\}_{n\in\mathbb{N}}$ is said to be \textit{asymptotically optimal} if for any sequence of admissible control policies $\{(\hat U^{n}_s, \hat U^{n}_b)\}_{n\in\mathbb{N}}$, 
\begin{align*}
\lim_{n\to\infty} \hat{J}^n( \hat{X}^{n}(0), \hat{U}^{n,*}_s, \hat{U}^{n,*}_b) \le \liminf_{n\to\infty}\hat{J}^n( \hat{X}^n(0), \hat{U}^n_s, \hat{U}^n_b).
\end{align*}

\section{Weak convergence}\label{sec:weak}
This section is devoted to establishing the tightness of the diffusion-scaled processes and to characterizing their weak limits. For a given queue-capacity process  ${\bf m}_n(\cdot)$,  the controlled diffusion-scaled state process  $\hat{X}^n$ described in \eqref{state-2}  can be written as follows: For all $t \geq 0,$ 
\begin{equation}
\hat{X}^n(t)=\hat{\zeta}^n(t)  -\hat{G}_s^n(t) + \hat{G}_b^n(t)  -\hat{U}_s^n(t) + \hat{U}_b^n(t),
\label{3.1}
\end{equation}
 where 
\begin{align}
\hat{\zeta}^n(t)& =\hat{X}^n(0) +\hat{ A}^n_s(t) - \hat{A}^n_b(t). \label{3.2}
\end{align}

We present the main result of weak convergence in the following theorem, and the rest of this section will focus on its proof. 
\begin{theorem}\label{thm3.2} Any sequence of the controlled diffusion-scaled processes {$\{(\hat{X}^n,  \hat{\zeta}^n, \hat{G}^n_b, \hat{G}^n_s, \hat{U}^n_b,  \hat{U}^n_s)\}_{n\in\NN}$} is C-tight in $D^6[0, T]$ for each $T\ge 0$. In particular, there exist a constant $C>0$ and an integer $m>1$ independent of the queue-capacity ${\bf m}_n$ such that for each $n\in \NN$ and $T>0$, 
\[
E[\|\hat{X}^n\|^2_T +  \|\hat{\zeta}^n\|^2_T+ \|\hat{G}^n_b\|^2_T + \|\hat{G}^n_s)\|^2_T] \le C(1+ T^m).
\]
Furthermore, let $({X},  {\zeta}, {G}_b, {G}_s, {U}_b,  {U}_s)$ denote a limit point. Then the following hold.
\begin{itemize}
\item[\rm (i)] There exists a standard one-dimensional Brownian motion $B$ such that for $t\ge 0$, 
\begin{align}\label{3.4}
\zeta(t) = x +  \sigma B(t) +\beta t, 
\end{align}
where $x$ is the limit initial value given in \eqref{initial}, $\sigma^2=  \sigma_s^2 + \sigma_b^2$, $\beta= \beta_s - \beta_b$ and the constants $\sigma_s,$  $\sigma_b,$  $\beta_s,$ and $\beta_b$ are given in \eqref{arrivals}, \eqref{HT-b} and \eqref{HT-s}.  
\item[\rm (ii)] For $t\ge 0$, 
\begin{align*}
G_s(t) = \delta_s \int_0^t X^+(u) du, \ \ G_b(t) = \delta_b \int_0^t X^-(u) du,
\end{align*}
where $\delta_s$ and $\delta_b$ are as in \eqref{Patience-s} and \eqref{Patience-b}.
\item[\rm (iii)] The process $X$ satisfies  the It$\hat{o}$ equation
\begin{equation} \label{3.90}
X(t)= x +\sigma B(t) +  \int_{0}^{t}[\beta - h(X(s))]ds -U(t),
\end{equation}
where $h$ is a piecewise linear function given by $h(x)= \delta_s x^+ - \delta_b x^-$ for $x\in\mathbb{R}$ and $U = U_s - U_b$ is a process of bounded variation, which is adapted to the filtration generated by $(X, B)$ and satisfies 
\begin{align*}
U_s(t)= \int_{0}^{t} 1_{\{X(u)>0\}} dU_s(u), \ \  U_b(t)= \int_{0}^{t} 1_{\{X(u)<0\}} dU_b(u). 
\end{align*}
\end{itemize}
\end{theorem}

\begin{remark}\label{rem:wc}
Theorem \ref{thm3.2} (i) follows from Assumptions \ref{assump:arrival} and \ref{assump:htc}. More precisely, from \eqref{arrivals},  \eqref{HT-b} and \eqref{HT-s}, 
 \begin{equation} 
 \big(\hat{A}^n_b(t),  \hat{A}^n_s(t)\big)   \text{ converges weakly to } (\sigma_b B_1(t)+ \beta_b t, \   \sigma_s B_2(t)+\beta_s t),
 \label{wk-cgce}
 \end{equation}
in the space  ${D}^2[0, T]$  for any $T>0.$  
From the initial condition \eqref{initial} and  \eqref{wk-cgce},  it follows that the process  $\hat{\zeta}^n$ is convergent weakly in $D[0, T]$ and its limiting process is given by $x+ \sigma B(t) +\beta t$ for $t\in[0,T].$
\end{remark}

The remainder  of the  proof of Theorem \ref{thm3.2} is divided into five subsections. Section \ref{subsec:sb} establishes the stochastic boundedness of $\hat X^n$. In Section \ref{subsec:vwt}, we introduce the virtual waiting times for unblocked customers and show that the diffusion-scaled virtual waiting time processes are stochastically bounded. Section \ref{subsec:tightG} is devoted to the C-tightness of $(\hat G^n_s, \hat G^n_b).$ The C-tightness of $(\hat X^n, \hat U^n_s, \hat U^n_b)$ and the asymptotic relationships between $\hat X^n$ and $\hat G^n_s$ and $\hat G^n_b$ are obtained in Section \ref{subsec:tightX}. Finally, in Section \ref{subsec:lastproof}, we  complete the proof of Theorem \ref{thm3.2} based on the results derived in Sections \ref{subsec:sb} -- \ref{subsec:tightX}. For notation convenience, we introduce 
\begin{align}
    \hat{G}^n(t)& =\hat{G}^n_s(t) - \hat{G}^n_b(t), \label{3.3} \\
\hat{U}^n(t)& = \hat{U}^n_s(t) -\hat{U}^n_b(t). \label{3.4a}
\end{align}

\subsection{Stochastic boundedness of \texorpdfstring{$\hat X^n$}{lg}}\label{subsec:sb}
For a given queue-capacity process  ${\bf m}_n(\cdot)$, consider the controlled diffusion-scaled  process $(\hat X^n, \zeta^n, \hat G^n_s, \hat G^n_b, \hat U^n_s, \hat U^n_b)$  satisfyng \eqref{3.1} and \eqref{3.2}.
 When the state  process   $ \hat{X}^n$  deviates  far away from the origin, the processes $ \hat{G}^n_s, \hat{G}^n_b, \hat{U}^n_s,$ and $\hat{U}^n_b$ act as frictional forces.  The proof of the following result is based on this fact.

\begin{proposition} \label{prop-3.1}
For any state process $\hat{X}^n$  in $ D[0, T]$, we have
\begin{equation}
E[\|\hat{X}^n\|_T^2] \leq C(1+ T^m),
\label{3.6}
\end{equation}
 where the constant $C>0$ and the integer $m > 1$ are independent of  $n,  T$ and of the queue-capacity  ${\bf m}_n(\cdot).$
Consequently, the sequence  $\{\hat{X}^n\}_{n\ge 1} $  is stochastically bounded.
\end{proposition}

\proof
From \eqref{arrivals-M},   \eqref{HT-b}  and  \eqref{HT-s}, it follows that
\begin{equation}
E \big[\|\hat{A}^n_b\|_T^2 +  \|\hat{A}^n_s\|^2_T\big]  \leq \tilde C (1+ T^{\tilde m}),
\label{arrivals-M2}
\end{equation}
where $\tilde C$ is a positive constant and $m>1$ is an integer constant, and both constants are independent of $n$,  $T$, and the queue-capacity ${\bf m}_n$. Thus there exists a constant $C>0$ so that  $E[\|\hat{\zeta}^n\|^2] \leq C(1+T^m).$  The constant $C$ is also independent of  $n,  T$ and the queue-capacity  ${\bf m}_n(\cdot).$  We let $Y_n(t)= \hat{X}^n(t)- \hat{\zeta}^n(t) $ for $t \geq 0$ and  then use a path-wise argument to obtain  \eqref{3.6}.
We  choose $M \equiv M(\omega)=1+ \|\hat{\zeta}^n\|_T <\infty $ a.s. by \eqref{arrivals-M2}. We claim that $ \|Y_n\|_T \leq 2 M.$  Suppose it doesn't hold. Then there exists a $t_0$ in  $[0, T]$ so that  $|Y_n(t_0)|> 2(1+ \epsilon)M$ for some $\epsilon >0.$

First we consider the case  $Y_n(t_0) > 2(1+ \epsilon)M$ for some $t_0.$ Then, for any such $t_0,$ $ \hat{X}^n(t_0)= Y_n(t_0) + \hat{\zeta}^n(t_0) > M. $
Since $Y_n(0)=0,$  we let   $\hat{t}=\inf \{t \in(0, T] : Y_n(t) >2(1+ \epsilon)M \}.$  Then $0< \hat{t} \leq T,$  $Y_n(\hat{t}) > 2M$  and $ Y_n(\hat{t})-Y_n(\hat{t}-) >0 $ since $Y_n$ has piecewise constant RCLL paths. Since $ \hat{X}^n(\hat{t}) >M>0$, we have $\hat{X}^n(\hat{t}-) \geq \hat X^n(\hat t) - \frac{1}{\sqrt{n}} \ge 0$ and $\hat G^n_b(\hat t)+ \hat U^n_b(\hat t) = \hat G^n_b(\hat t-)+ \hat U^n_b(\hat t-)$.   Then $Y_n(\hat{t}) - Y_n(\hat{t}-) =-[(\hat{G}^n_s(\hat{t}) +\hat{U}^n_s(\hat{t})) - (\hat{G}^n_s(\hat{t}-) +\hat{U}^n_s(\hat{t}-))] \leq 0.$ This is a contradiction. 

A similar argument shows that  $Y_n(t) < -2(1+\epsilon)M$ also not possible when $0 \leq t \leq T.$ Therefore $\|Y_n\| \leq 2M $ holds a.s.  Since  $Y_n(t)= \hat{X}^n(t)- \hat{\zeta}^n(t), $ this yields that  $\|\hat{X}^n\|_T \leq 3M=3( 1+ \|\hat{\zeta}^n\|_T). $  Now the conclusion \eqref{3.6} can be obtained using the moment bound in \eqref{arrivals-M2}.  This completes the proof.
\Halmos

The following corollary is  an immediate consequence of the above result.

\begin{corollary} \label{cor-3.1}
For any state process $\hat{X}^n$  in $ D[0, T]$, 
\begin{equation*}
E[\|\bar{X}^n\|_T^2] \leq \frac{C(1+ T^m)}{n},
\end{equation*}
where $C$ and $m$ are as in Proposition \ref{prop-3.1}.
Consequently, 
\begin{equation}
\lim\limits_{n \rightarrow \infty} E[\|\bar{X}^n\|_T^2] =0.   
 \label{3.7}
\end{equation}
\end{corollary}

\qquad

\subsection{Virtual waiting times}\label{subsec:vwt}

We need to introduce the virtual waiting time processes and obtain their stochastic boundedness to guarantee the tightness of the state process. However, the virtual waiting times  can be undefined on the time intervals during which the corresponding buffer is full. To circumvent this difficulty, we  introduce the auxiliary arrival processes of  unblocked buyers and sellers  by 
\begin{equation*}
E^n_b(t)=A^n_b(t)-U^n_b(t) \qquad\text{and}\qquad 
E^n_s(t)=A^n_s(t)-U^n_s(t), \qquad \text{for $t\ge 0$.}
\end{equation*} 
 It is evident that they are adapted to the filtration  $ \{\mathcal{F}^{n}_t\}$ defined in \eqref{filtration}. Their corresponding fluid-scaled and diffusion-scaled  processes are respectively given by 
\begin{align*}
& \bar{E}^n_b(t)=\bar{A}^n_b(t)-\bar{U}^n_b(t), \qquad\text{and}\qquad 
\bar{E}^n_s(t)=\bar{A}^n_s(t)-\bar{U}^n_s(t), \\
& \hat{E}^n_b(t)=\hat{A}^n_b(t)-\hat{U}^n_b(t), \qquad\text{and}\qquad 
\hat{E}^n_s(t)=\hat{A}^n_s(t)-\hat{U}^n_s(t).
\end{align*} 
We use these  arrival processes of unblocked customers to introduce the virtual waiting times.  Let $V^n_b(t)$ be the waiting time of an infinitely-patient buyer who arrives 
from  $E^n_b(\cdot)$  at time $t \geq 0.$ We define the waiting time $V^n_s(t)$ of an infinitely-patient seller who arrives
from  $E^n_s(\cdot)$  at time $t \geq 0$ accordingly. The diffusion-scaled virtual waiting times are defined by 
\begin{equation}
\hat{V}^n_b(t)=\sqrt{n} V^n_b(t), \qquad\text{and}\qquad 
\hat{V}^n_s(t)=\sqrt{n} V^n_s(t) \qquad \text{for $t\ge 0$.}
\label{3.14}
\end{equation}

In the following result, we obtain the moment bounds for $\hat{V}^n_b(t)$    and   $\hat{V}^n_s(t)$  in $D[0, T].$
\begin{proposition}
 \label{prop-3.3}
For each  $T>0,$  let $\hat{V}^n_b$    and   $\hat{V}^n_s$   be the virtual waiting time processes in  $D[0, T].$  Then  
\begin{equation}
    E[\| \hat{V}^n_b\|_T^2 + \|\hat{V}^n_s\|_T^2  ] \leq C(1+T^m),
    \label{3.99}
    \end{equation}
     where the constant $C>0$ and the integer  $m>1$ are independent of  $n$, $T$, and the queue-capacity ${\bf m}_n$. Consequently, the  processes  $\{\hat{V}^n_b\}$    and   $\{\hat{V}^n_s\}$     are   stochastically bounded     in $ D[0, T].$ 

%
\end{proposition}

\proof
 We establish a moment bound for   $ E[\| \hat{V}^n_b\|_T^2 ].$    The moment bound for   $ E[\| \hat{V}^n_s\|_T^2 $ is similar.  Let $t^n_j$ and $d^n_j$  represent the arrival time and the patience time of the  $j$th  buyer (according to the order of arrival), respectively.  Similarly, we let $\tilde{t}^n_j$ and $\tilde{d}^n_j$  represent the arrival time and the patience time of the  $j$th  seller (according to the order of arrival), respectively. 
 
 If there are $K_n$ sellers waiting in the queue initially, the first $K_n$ buyers will be matched upon their arrival.
In this case, using Assumption \ref{assump:initial}, $0 \leq K_n \leq M \sqrt{n}$ and they are all non-abandoning. Let $T_n$ be the time it takes to match all these buyers. Then $T_n \leq  \sum_{j=1}^{K_n} \tilde{\tau}^n_j $ and hence using Cauchy-Shwartz inequality, $E [T_n^2] \leq K_n^2 E [\tilde{\tau}_1^2]. $  By Assumption \ref{assump:arrival}, $E [\tilde{\tau}_1^2] \leq {(1+\varsigma_b^2)}/{(\lambda^n_b)^2}$  and
consequently, $E [(\sqrt{n}T_n)^2] \leq  M^2   {(1+\varsigma_b^2)n^2}/{(\lambda^n_b)^2}.$  Now using Assumption \ref{assump:htc}, $\lim_{n \rightarrow \infty}{\lambda^n_b}/{n}=\lambda_0>0$ and hence, $\sup_{ n \geq 1} E [(\sqrt{n}T_n)^2] \leq C <\infty,$ where $C>0$ is  a generic constant. Therefore, we can simply assume that there are no buyers initially.

    Let $v^n_j$ be the amount of time the $j$th buyer spent as the head of the queue (i.e.  the  time spent  in the first place of the queue).  We simply take $v^n_j=0$ if the $j$th buyer  did not reach the head position due to blocking or abandonment. Then we can write
 \begin{equation*}
V^n_b(t)= \sum\limits_{j=1}^{A^n_b(t)} v^n_j  -  \int_{0}^{t}  1_{\{V^n_b(s)>0\}} ds.
\end{equation*}
We introduce $S^n_b(t)$ to be the virtual waiting time to reach the head of the queue for an infinitely patient unblocked buyer who arrived at time $t.$ To simplify the notation, we write $s^n_j \equiv S^n_b(t^n_j-).$   
 With a simple algebraic manipulation, we observe that  the diffusion scaled virtual waiting time $\hat{V}^n_b(t)= \sqrt{n}V^n_b(t)$ satisfies
 \begin{equation*}
\hat{V}^n_b(t)= Y_n(t) +  \sqrt{n} \int_{0}^{t}  1_{\{ \hat V^n_b(s)=0\}} ds,
\end{equation*}
where  $Y_n$ is given by 
\begin{equation}
Y_n(t)=\sqrt{n} \left[ \sum\limits_{j=1}^{A^n_b(t)} v^n_j  -t \right]
\label{3.101}
\end{equation}
for all $t\geq 0.$ Therefore, the pair $\{(\hat{V}^n_b(t),    \sqrt{n} \int_{0}^{t}  I_{\{\hat V^n_b(s)=0\}} ds); {t \geq 0}\}$ is  the unique solution to the SP with the input process $Y_n$ and the reflection barrier  at the origin (see \cite{kruk}). {Let $\Gamma$ denote the one-sided Skorokhod map with reflection barrier at the origin. Then $\hat{V}^n_b(t)= \Gamma(Y_n)(t)$ for $t\in [0, T]$. From the Lipschitz continuity of $\Gamma$ (see \cite{kruk} and \cite{weera4}), we have $\|\hat V^n_b\|_T \leq 2 \|Y_n\|_T$. In the following, we estimate the second moment of $\|Y_n\|_T$. }


For  each $j \geq 1,$ introduce the set $\mathcal{V}^n_j=[s^n_j, s^n_j+v^n_j).$  Notice that the length (Lebesgue measure) of the set $|\mathcal{V}^n_j|=v^n_j$ and this holds even when $v^n_j=0.$   The sets $\{\mathcal{V}^n_j\}$ are disjoint and for each $j,$ $\mathcal{V}^n_j \subseteq [\tilde{t}^n_k,  \tilde{t}^n_{k+1})$ for some $k.$ It is important to observe that the same set $ [\tilde{t}^n_k,  \tilde{t}^n_{k+1}) $ may contain several $ \mathcal{V}^n_j$ intervals. Next we consider the collection of the intervals $$ \mathcal{A}(t)=\{ [\tilde{t}^n_k,  \tilde{t}^n_{k+1}) : {k\in\NN \ \text{and}} \ [\tilde{t}^n_k,  \tilde{t}^n_{k+1}) \cap \mathcal{V}^n_j \  \text{is nonempty for some } j = 1,\ldots, A^n_b(t)\}.$$
This collection  of intervals is finite and we let $N^n_b(t)\equiv |\mathcal{A}(t)|$ which represents  the number of elements in this set.  For a given $t\geq 0,$  the intervals in $ \mathcal{A}(t)$ are called ``good'' intervals. For $t \geq 0$, let $h(t, k)$ denote the number of $\mathcal{V}^n_j$'s in $[\tilde{t}^n_k,  \tilde{t}^n_{k+1}) \in \mathcal{A}(t)$ for each $k$. Then $h(t, k)$ is non-decreasing in $t$ and $A^n_b(t)-N^n_b(t)= \sum_{k=1}^{\infty}(h(t, k)-1)^+.$  Thus,  $A^n_b(t)-N^n_b(t)$ is a non-negative, non-decreasing process and we will use this fact to obtain our final estimate (\ref{10}). We further introduce the scaled quantities $\bar{N}^n_b(t)= N^n_b(t)/n$ and $\hat{N}^n_b(t) = (N^n_b(t)-\lambda^n_b t)/\sqrt{n}.$
 Next introduce the non-negative integer valued random variables  $$\tau^n_1= \min\{k\geq 0 : [\tilde{t}^n_k, \tilde{t}^n_{k+1}) \cap \mathcal{V}^n_j \text{ is non-empty for some  } j= 1, \ldots, A^n_b(t)\},$$ and $\tau^n_1$ is infinite if the above set is empty. Let 
$$\tau^n_j= \min\{k> \tau^n_{j-1} : [\tilde{t}^n_k, \tilde{t}^n_{k+1}) \cap \mathcal{V}^n_l \text{ is non-empty for some  } l = 1, \ldots, A^n_b(t)\},$$  
and $\tau^n_j$ is infinite if the above set is empty. Introduce the filtration  $\{\mathcal{G}^n_k\}_{k\in\NN}$  by $ \mathcal{G}^n_k=\sigma\{ (\tilde{t}^n_1, \tilde{d}^n_1), ....., (\tilde{t}^n_{k}, \tilde{d}^n_{k}), (t^n_l, d^n_l) \text{ for all  } l\geq 1\}$. Then $\tau^n_j=k$ if and only if there are exactly $j-1$ ``good'' intervals among $[0, \tilde{t}^n_1),....[\tilde{t}^n_{k-1}, \tilde{t}^n_{k})  $ and the interval $[\tilde{t}^n_k, \tilde{t}^n_{k+1})$ is also ``good''. Consequently,  $\tau^n_j+1$ is a $\{\mathcal{G}^n_k\}$ stopping time.
Since the sets $\{\mathcal{V}^n_j\}$ are disjoint and for each $j,$ $\mathcal{V}^n_j \subseteq [\tilde{t}^n_k,  \tilde{t}^n_{k+1})$ for some $k,$  it follows that   
 $$ Y_n(t) \leq W_n(t)\equiv\sqrt{n} \left[ \sum_{j=1}^{N^n_b(t)} \tilde{u}^n_{\tau^n_j+1}  -t \right],$$
 where $\{\tilde{u}^n_k\}_{k\in\NN}$ are the inter-arrival times of the arrival process $A^n_s(t)$ of the sellers.    
Focusing on $W_n(t)$ now, we establish  that $E[\tilde{u}^n_{\tau^n_j+1} ] = {1}/{\lambda^n_s}$ and  the Var$(  \tilde{u}^n_{\tau^n_j+1} )= ({\varsigma_s}/{\lambda^n_s})^2$  for each $j.$ Indeed, keep $k$ fixed and consider the filtration $\{\mathcal{G}^n_k\}$ described above. The random variables $\{\tilde{u}^n_{1},...., \tilde{u}^n_{k-1}\}$ are $\mathcal{G}^n_k$ adapted and the  sequence $\{ \tilde{u}^n_{k}, \tilde{u}^n_{k+2},...\} $ are independent of $\mathcal{G}^n_k.$ Hence, we can employ the Wald's equation for random sums to conclude that $E[  \tilde{u}^n_{\tau^n_j+1} ] = {1}/{\lambda^n_s}.$  A similar proof yields  $E[(\tilde{u}^n_{\tau^n_j+1})^2] = ({1}/{\lambda^n_s})^2 +({\varsigma_s}/{\lambda^n_s})^2$  and consequently, Var$(\tilde{u}^n_{\tau^n_j+1} ) = ({\varsigma_s}/{\lambda^n_s})^2.$ 
 Next we define the stopped $\sigma$-algebra $\mathcal{H}^n_j= \mathcal{G}^n_{\tau^n_j+1}.$ Then  $\{\mathcal{H}^n_j\}$ is a filtration.  We use this filtration to introduce a discrete-time  martingale  $\{({M}^n(k),   \mathcal{H}^n_k)\}_{k \in\NN}$  given by  ${M}^n(0)=0$  and 
$ {M}^n(k)= {\sqrt{n}} \sum_{j=1}^{k}[\tilde{u}^n_{\tau^n_j+1}- {1}/{\lambda^n_s}].$
  It is straight forward to check that   $\{({M}^n(k),   \mathcal{H}^n_k)\}_{k \in\NN}$  is a square integrable martingale. 
   Next, we introduce $ \mathcal{H}^n_t= \mathcal{H}^n_{[nt]}$  and  $ M^n(t)= {M}^n([nt])$ for all $t \geq 0$,  where $[x]$ represents the integer part of $x.$
 Then it can be easily checked that $\{(M^n(t),  \mathcal{H}^n_t))\}_{t \geq 0}$ is a square integrable, pure-jump martingale  and its quadratic variation process is given by
$[M^n, M^n](t)= n  \sum_{j=1}^{[nt]}[\tilde{u}^n_{\tau^n_j+1}- {1}/{\lambda^n_s}]^2$
for all $t \geq 0.$  
With a simple algebraic manipulation, we can represent $W_n$  by
\begin{equation}
W_n(t)= M^n(\bar{N}^n_b(t)) +  \frac{n}{\lambda^n_s}\left[ \hat{N}^n_b(t)+\left(\frac{\lambda^n_b-\lambda^n_s}{\sqrt{n}}\right)t\right], \text{ for} \quad 0 \leq t \leq T.
\label{3.104}
\end{equation}
{Based on the formulation of $W_n$,} we introduce the process 
\begin{align}\label{eq:z-process}
    Z_n(t)= M^n(\bar{N}^n_b(t)) +  \frac{n}{\lambda^n_s}\left[ \hat{A}^n_b(t)+\left(\frac{\lambda^n_b-\lambda^n_s}{\sqrt{n}}\right)t\right], \text{ for} \quad 0 \leq t \leq T.
\end{align}    
Since $ A^n_b(t)-N^n_b(t)$ is a non-negative, non-decreasing process, it follows that 
\begin{equation}
Y_n(t) \leq W_n(t) \leq Z_n(t), \text{ for} \quad 0 \leq t \leq T.
\label{9}
\end{equation}
Moreover,  $W_n(t)-Y_n(t)$ and $Z_n(t)-Y_n(t)$ are non-negative, non-decreasing processes. Consequently, {using the monotonicity property of $\Gamma$ (that is, for $f, g\in D[0, T]$ and if $g$ is a non-negative, nondecreasing function,  then $ \Gamma(f)(t) \leq \Gamma(f+g)(t)$ (see again \cite{kruk} and section 2.3 of \cite{weera4}))}, we obtain 
\begin{equation}
 \hat{V}^n_b(t)= \Gamma(Y_n)(t) \leq \Gamma(Z_n)(t), \text{  for all   } 0 \leq t \leq T.
 \label{10}
 \end{equation}
{Thus it suffices to estimate the second moment of $\|Z_n\|_T$.} Applying the heavy traffic condition in Assumption \ref{assump:htc}, and the fact that $ \bar{N}^n_b(t)\leq \bar{A}^n_b(t) \leq \bar{A}^n_b(T)$ for all $0 \leq t \leq T$ to \eqref{eq:z-process}, we have
\begin{equation*}  
0 \leq \|Z_n\|_T \leq   \sup\limits_{0\leq s \leq \bar{A}^n_b(T)} |M^n(s)| +    \frac{n}{\lambda^n_s}  \|\hat{A}^n_b\|_T + K T, 
\end{equation*}
where $K>0$ is a constant which depends on the heavy traffic limit $\beta_b-\beta_s$.
Notice that  $\{\bar{A}^n_b(T)\leq r\}= \{{A}^n_b(T) \leq [nr]\} \in\mathcal{H}^n_r  $ for each $r>0,$  and hence,   $\bar{A}^n_b(T)$ is a stopping time with respect to the filtration $ \{\mathcal{H}^n_t\}_{t \geq 0}.$  In the rest of the proof, the generic constants $C_i>0$ where $i=1,...,5$  are independent of $T$ and $n.$ We then obtain   
\begin{align*}
&E \left[ \sup_{0\leq t \leq \bar{A}^n_b(T)} | M^n(t)|^2\right]  \leq C_1 E( [M^n, M^n](\bar{A}^n_b(T))) = C_1 E[E( [M^n, M^n](\bar{A}^n_b(T))|\bar{A}^n_b(T)]\\
& = C_1 n({\varsigma_s }/{\lambda^n_s})^2 E[\bar{A}^n_b(T)] \le  \frac{C_1}{n}\left(\frac{n\varsigma_s }{\lambda^n_s}\right)^2\left(\frac{\lambda^n_bT}{n}  +C_2 \sqrt{\frac{(1+T^m)}{n}}\right),
\end{align*} 
where the first inequality follows from the Burkholder's inequality (see \cite{protter}) and the second inequality follows from \eqref{arrivals-M}. 
Next, from \eqref{arrivals-M}, 
$E[\|\hat{A}^n_b\|^2_T] \leq C_3 (1+T^m)$ and hence, it follows that  $E[\|Z_n\|^2_T] \leq C_4(1+T^m)$   where the constant  $ C_4 >0$ and the integer $m>1$ are independent of $n$ and $T>0.$
 Combining this estimate   with  \eqref{10} and using the properties of Skorokhod map $\Gamma$, we obtain 
\begin{equation*}
 E[ \|\hat V^n_b\|^2_T] \leq E[\|\Gamma(Z_n)\|^2_T] \leq 4 E[\|Z_n\|^2_T] \leq C_5(1+T^m). 
\end{equation*}
This yields the moment bound for   $ E[ \|\hat V^n_b\|^2_T].$
 The  estimate for   $E[  \|\hat{V}^n_s\|^2_T  ]$ is similar. Hence \eqref{3.99} follows and consequently,  $\{\hat{V}^n_b\}_{n\in\NN}$ and $\{\hat{V}^n_s\}_{n\in\NN}$ are stochastically bounded  in  $D[0, T].$ This completes the proof. \Halmos

\subsection{ Tightness of  \texorpdfstring{$\{(\hat{G}^n_b, \hat{G}^n_s)\}_{n\ge 1}$}{lg} }\label{subsec:tightG}

We first introduce a pair of ``eventual" abandonment processes  $\hat{L}^n_b$ and  $\hat{L}^n_s$. For $t \geq 0,$ let $L^n_b(t)$ (resp. $L^n_s(t)$) denote the number of buyers (resp. sellers) who enter the queue by time $t$ and eventually abandon the system. 
Then $L^n_b(0)=L^n_s(0)=0$  and  it is evident that  $  G^n_b(t)\leq  L^n_b(t) \leq  E^n_b(t)$  and  $  G^n_s(t)\leq  L^n_s(t) \leq  E^n_s(t)$  for all $t \geq 0.$ Notice that the processes  ${L}^n_b$ and  ${L}^n_s$ are not adapted to the natural filtration   $ \{\mathcal{F}^{n}_t\}$ defined in \eqref{filtration}. Their diffusion-scaled  versions are given by $\hat{L}^n_b(t)={L^n_b(t)}/{\sqrt{n}}$ and $\hat{L}^n_s(t)={L^n_s(t)}/{\sqrt{n}}$ for $t\ge 0.$

 Our first step is to show that the processes $\{\hat{L}^n_b\}$ and   $\{\hat{L}^n_s\}$  are tight. Here we appropriately modify a proof originally developed in \cite{reed} for a single-server queue in conventional heavy traffic.

\begin{proposition}
 \label{prop-3.4}
For the  non-decreasing processes  $\hat{L}^n_b,$   $\hat{G}^n_b,$     $\hat{L}^n_s$  and   $\hat{G}^n_s$  in $ D[0, T]$,
 \begin{align}
    E[ \hat{L}^n_b(T)^2 + \hat{L}^n_s(T)^2] &\leq C(1+T^m),
    \label{3.21} \\
    E[ \hat{G}^n_b(T)^2 + \hat{G}^n_s(T)^2] &\leq C(1+T^m),
    \label{3.22}
    \end{align} 
    where the constant $C>0$ and the integer  $m> 1$ are independent of $n,T$, {and the queue-capacity ${\bf m}_n$}. Consequently, the processes  $\{(\hat{L}^n_b, \hat{G}^n_b, \hat{L}^n_s , \hat{G}^n_s)\}_{n\in\NN}$ are stochastically bounded in $ D^4[0, T].$
\end{proposition}

\proof
We will  prove  the results for buyers and the case of sellers is very similar. Let $t^n_j$ and $d^n_j$  represent the arrival time and the patience time of the  $j$th unblocked  buyer respectively, according to the order of arrival.  Since initial customers are of infinite patience, we   can write   $\hat{L}^n_b(t)$ in the form 
\begin{equation*}
 \hat{L}^n_b(t)=
 \begin{cases} 0  &  \text{ if }  \quad 0 \leq t <t^n_1, \\
  \frac{1} {\sqrt{n}} \sum_{j=1}^{E^n_b(t)}  1_{\{V^n_b(t^n_j-)>d^n_j\}}, & \text{  if  }  \quad t \geq t^n_1. 
  \end{cases}
\end{equation*}
Since  $\hat{L}^n_b$ is not adapted to the filtration $\{\mathcal{F}^{n}_t\}_{t\ge 0}$ defined in \eqref{filtration}, we now  construct a new filtration $\{\mathcal{G}^n_t\}_{t\ge 0}$ where it will be adapted. Let $\tilde{t}^n_j$ and $\tilde{d}^n_j$  represent the arrival time and the patience time of the  $j$th unblocked  seller, respectively, according to the order of their arrival.  Consider $ \mathcal{H}$ to be the $\sigma-$algebra generated by the sequence $\{(\tilde{t}^n_j,\tilde{d}^n_j)\}_{j\geq 1}$ and introduce 
$ \mathcal{G}^n_0=\sigma(t^n_1) \vee \mathcal{H}$
to be the $\sigma$-algebra generated by the random variable $t^n_1$ and the collection  $\mathcal{H}.$  Similarly, for each $k \geq 1,$ we introduce the $\sigma$-algebra 
 \begin{equation*}
 \mathcal{G}^n_k=\sigma((t^n_1, d^n_1),.., (t^n_k, d^n_k), t^n_{k+1}) \vee \mathcal{H}.
\end{equation*}
Notice $V^n_b(t^n_j-)$  is adapted to  $ \mathcal{G}^n_{j-1}$ and  $d^n_j$ is independent of    $ \mathcal{G}^n_{j-1}.$ Let  $\hat{M}^n(0)=0$  and 
 \begin{equation*}
 \hat{M}^n(k)= \frac{1}{\sqrt{n}} \sum_{j=1}^{k}\left[ 1_{\{V^n_b(t^n_j-) \geq d^n_j\}}-F_b(V^n_b(t^n_j-))\right],
\end{equation*}
where $F_b$ is in \eqref{Patience-b}. Then $\{(\hat{M}^n(k),   \mathcal{G}^n_k)\}_{k \geq 1}$  is a square integrable mean zero martingale.  Next, let  $ \mathcal{G}^n_t= \mathcal{G}^n_{[nt]}$  and  $ M^n(t)= \hat{M}^n([nt])$  for all $t \geq 0.$ Consequently, $\{(M^n(t),  \mathcal{G}^n_t))\}_{t \geq 0}$ is a square integrable, mean zero, pure-jump martingale  and  
 \begin{equation}
[M^n, M^n](t) \leq 2t 
 \label{3.34}
\end{equation}
for all $t \geq 0$ (for details, we refer to \cite{reed}).  We can represent  $\hat{L}^n_b$ by 
 \begin{equation}
\hat{L}^n_b(t)=M^n(\bar{E}^n_b(t)) +   \frac{1} {\sqrt{n}} \sum_{j=1}^{E^n_b(t)} F_b(V^n_b(t^n_j-))
 \label{3.36}
\end{equation}
for all $t \geq 0$, where  $ \sum_{j=1}^{0}$  is interpreted as zero. Next, we use this representation to  prove \eqref{3.21}. 
Let $T >0$ that is fixed. First we show that  $\bar{E}^n_b(T)$ is a discrete-valued stopping time with respect to the filtration $\{\mathcal{G}^n_t\}_{t \geq 0}.$  When $nr$ is a positive integer,  
$\{\bar{E}^n_b(T)= r\}=\{ {E}^n_b(T)= nr\}=\{ t^n_{nr}\leq T < t^n_{nr+1}\} \in  \mathcal{G}^n_r.$
 Similar argument holds when $r=0.$ Hence   $\bar{E}^n_b(T)$ is a stopping time. In what follows, the positive constants $C_i, i=1, \ldots, 5$ and integer $m>1$ are independent of $T$ and $n.$ By \eqref{arrivals-M2}, 
\begin{equation}
 E\left(\sup\limits_{t\in [0, T]} \big[ |\bar{A}^n_b(t)- \lambda_0 t|^2  +  |\bar{A}^n_s(t)- \lambda_0 t|^2 ] \right)\leq \frac{C_1(1+T^m)}{n}. 
\label{arrivals-M22}
\end{equation}
Therefore,
$E[(\bar{E}^n_b(T)])^2 \leq E[(\bar{A}^n_b(T)] )^2\leq (\lambda_0 T)^2 + C_1 (1 +T^m)/n.$
  Using  \eqref{3.34} together with  the Burkholder's inequality (see \cite{protter}),  we obtain  
$E[\sup_{0\leq t \leq \bar{E}^n_b(T)} | M^n(t)|^2] \leq C_2 E( [M^n, M^n](\bar{E}^n_b(T))) \leq 2C_2 E [ \bar{E}^n_b(T) ] \leq C_3 (1+T^m).$  
Consequently,   
\begin{align}\label{potenaband-est1}
 E \left[\sup\limits_{0\leq t \leq T} |M^n(\bar{E}^n_b(t))|^2\right] \leq E\left[\sup_{0\leq t \leq \bar{E}^n_b(T)} | M^n(t)|^2\right] \le C_3 (1+T^m).
 \end{align}
We now consider the second term in the RHS of \eqref{3.36}.
By \eqref{Patience-b}, we have $F_b(x) \leq Kx$ for all  $x \geq 0,$ where $K>0$ is a constant independent  of  $n.$
Since $ 0 \leq \bar{E}^n_b(T) \leq   \bar{A}^n_b(T)$, we have the following estimates.
  \begin{align*}
  \sqrt{n} F_b(\|V^n_b\|_T)  \bar{E}^n_b(T) & \leq \sqrt{n} [F_b(\|V^n_b\|_T) |\bar{A}^n_b(T)- \lambda_0 T| +   \lambda_0 T  \sqrt{n} F_b(\|V^n_b\|_T) \\
  &  \leq \sqrt{n}  |\bar{A}^n_b(T)- \lambda_0 T|  + K \lambda_0 T \|\hat{V}^n_b\|_T, 
  \end{align*}
and thus 
\begin{align*}E\left[ \left(\sup_{0\leq t \leq T} \frac{1} {\sqrt{n}} \sum_{j=1}^{E^n_b(t)} F_b(V^n_b(t^n_j-))\right)^2\right] & \le E\left[(\sqrt{n} F_b(\|V^n_b\|_T)  \bar{E}^n_b(T))^2\right] \\
& \leq 2 n E( |\bar{A}^n_b(T)- \lambda_0 T|^2) +2(  \lambda_0 KT)^2 \|E(\|\hat{V}^n_b\|_T^2).
\end{align*}
  Now using  \eqref{arrivals-M22} and \eqref{3.99}, we obtain   
  \begin{align}\label{potenaband-est2}
   E\left[ \left(\sup\limits_{0\leq t \leq T}\frac{1} {\sqrt{n}} \sum_{j=1}^{E^n_b(t)} F_b(V^n_b(t^n_j-))\right)^2\right] \leq C_4(1+T^m).
   \end{align}  
      Finally, using the estimates in \eqref{potenaband-est1} and \eqref{potenaband-est2} in \eqref{3.36}, we obtain $ E[\hat{L}^n_b(T)^2] \leq    C_5 (1+T^m). $ Similar estimate can be obtained for  $ E[\hat{L}^n_s(T)^2] $ and hence \eqref{3.21} follows.
   
Since $0 \leq \hat{G}^n_b(t) \leq \hat{L}^n_b(t)$ holds for all $t \geq 0$,  \eqref{3.22}  follows from \eqref{3.21}, and the stochastic boundedness is an  immediate consequence of  \eqref{3.21} and \eqref{3.22}.  \Halmos

The next proposition is concerned with the fluid-scaled control processes and the fluid-scaled arrival processes of unblocked customers.
\begin{proposition}\label{prop-3.44}
 The process sequrence $\{(\bar{U}^n_b, \bar{U}^n_s)\}_{n\in\NN}$ is stochastically bounded and C-tight in $ D^2[0, T]$ and  
\begin{equation}
\lim\limits_{n \rightarrow \infty} E[\bar{U}^n_b(T)^2  +  \bar{U}^n_s(T)^2] =0.
 \label{3.364}
\end{equation}
Consequently, 
\begin{equation}
\lim\limits_{n \rightarrow \infty} E \left(\sup\limits_{t\in[0, T]} \left[ |\bar{E}^n_b(t)- \lambda_0 t|^2  +  |\bar{E}^n_s(t)- \lambda_0 t|^2 \right] \right)=0.
\label{arrivals-M4}
\end{equation}

%
%
%
%
%
%

\end{proposition}

\proof
We focus on $\bar{U}^n_b$ and the proof for $\bar{U}^n_s$ is similar. For a given controlled  buffer length policy   ${ \bf m}_n(\cdot)$, the fluid-scaled state equation is given by 
\begin{equation}
\bar{X}^n(t)=\bar{X}^n(0) + \bar{A}^n_s(t) - \bar{A}^n_b(t) -\bar{G}^n_s(t) + \bar{G}^n_b(t)  +\bar{U}^n_b(t) -\bar{U}^n_s(t), \ \ t\ge 0.
\label{state-1}
\end{equation}
From \eqref{eq2.8b} and \eqref{eq2.8s}, the processes $\bar{U}^n_b$ and $\bar{U}^n_s$  satisfy 
\begin{equation}
\bar{U}^n_b(t)=\int_{0}^{t}1_{\{\bar{X}^n(s)=\bar{m}^n_b(s)\}}d\bar{A}^n_b(s)
 + [\bar{m}^n_b(0)-\bar{X}^n(0-)]^+,  
 \label{eq3.381}
\end{equation}
and 
\begin{equation}
\bar{U}^n_s(t)=\int_{0}^{t}1_{\{\bar{X}^n(s)=\bar{m}^n_s(s)\}}d\bar{A}^n_s(s)
 + [\bar{X}^n(0-)- \bar{m}^n_s(0)]^+.  
 \label{eq3.382}
\end{equation}
We can use  $ \bar{U}^n_b(T)\leq |\bar{X}^n(0)| + \bar{A}^n_b(T) $ together with  \eqref{initial} and \eqref{arrivals-M2} to obtain
\begin{equation}
E[\bar{U}^n_b(T)^2 ] \leq  C[1 +T^m],
 \label{eq3.384}
\end{equation}
where the constants $C>0,$ and the integer $ m > 1$ are independent of $n, T$ {and the queue-capicity ${\bf m}_n.$}  
Since  the process  $\bar{U}^n_b$  is non-decreasing, \eqref{eq3.384} yields that  $\{\bar{U}^n_b\}_{n\in\NN}$  is stochastically bounded in $ D[0, T].$ 

To prove the C-tightness of    $\{\bar{U}^n_b\}_{n\in\NN}$ in $D[0, T],$  let  $0 \leq t_1 \leq t_2 \leq T,$ then from \eqref{eq3.381}, we have 
 $0 \leq \bar{U}^n_b(t_2) - \bar{U}^n_b(t_1) \leq   \bar{A}^n_b(t_2) - \bar{A}^n_b(t_1).$
 Consequently, for any $\delta>0,$ 
 \begin{align*}
  \sup\limits_{|t_1-t_2|<\delta, t_1, t_2   \in[0, T] } | \bar{U}^n_b(t_2) - \bar{U}^n_b(t_1)| \leq  \lambda_0 \delta +   2 \sup\limits_{t\in [0, T]} | \bar{A}^n_b(t) -\lambda_0 t|.    
 \end{align*}
Using \eqref{arrivals-M2}, it follows that
 \begin{equation}
\limsup\limits_{n \rightarrow \infty} E\left[ \sup\limits_{|t_1-t_2|<\delta } | \bar{U}^n_b(t_2) - \bar{U}^n_b(t_1)|\right]^2  \leq  \lambda_0^2 \delta^2 
\label{3.388}
\end{equation}
for any $\delta >0.$ The conditions  \eqref{eq3.384}  and \eqref{3.388} (see \cite{billingsley})   imply that the process   $\{\bar{U}^n_b\}$  is C-tight in $ D[0, T].$  Similarly, $\{\bar{U}^n_s\}$ can also been shown to be C-tight in $ D[0, T].$ Let  $(\bar{u}_b, \bar{u}_s)$ be a weak limit  of the joint process through a subsequence.  We relabel this subsequence as the original sequence $ (\bar{U}^n_b,  \bar{U}^n_s)$  for convenience. We next show that  $\bar{u}_b$ and  $\bar{u}_s$  are identically zero in $[0, T].$

Now we let 
\begin{equation}
\bar{\xi}^n(t)= \bar{X}^n(0) + \bar{A}^n_s(t) - \bar{A}^n_b(t) -\bar{G}^n_s(t) + \bar{G}^n_b(t), \ t\ge 0.
\label{eq3.380}
\end{equation}
 Following Theorem 2.6 and Corollary 2.4 in \cite{burdzy},
 given the input process $\bar{\xi}^n$ in \eqref{eq3.380}, the equations \eqref{state-1}, \eqref{eq3.382} and \eqref{eq3.384} yield a unique solution  $( \bar{X}^n,  \bar{U}^n_b-\bar{U}^n_s)$  for the extended SP (ESP) associated with the reflection barriers  $ \{(\bar{m}^n_b(t),  \bar{m}^n_s(t)) : 0 \leq t \leq T \}.$  By  \eqref{3.4} and \eqref{3.7}, $\lim_{n \rightarrow\infty}( \bar{\xi}^n,  \bar{X}^n)=(0, 0)$  uniformly on $[0, T].$  Moreover, by \eqref{cap-3}, $\lim_{n \rightarrow\infty}( \bar{m}^n_b,  \bar{m}^n_s)=(0, 0)$  uniformly on $[0, T].$  Now by the closure property  in Proposition 2.5 of \cite{burdzy},  the pair of zero functions $(0, 0)$ yields a unique solution to the ESP in the degenerate  interval $[0, 0]. $  Hence, by the uniqueness of the solution in Theorem 2.6 of \cite{burdzy}, the limiting processes $\bar{u}_b$ and  $\bar{u}_s$  are identically zero in $[0, T].$  Consequently, 
 $ (\bar{U}^n_b,  \bar{U}^n_s)$  converges to  $(0, 0)$    in $ D^2[0, T].$  
 
 To obtain  \eqref{3.364}, we prove $ \lim_{n \rightarrow \infty} E[\bar{U}^n_b(T)^2 ] =0.$ The proof of $ \lim_{n \rightarrow \infty} E[\bar{U}^n_s(T)^2 ] =0$ is similar.   By  \eqref{eq3.384},  
$\lim_{a \rightarrow \infty} \sup_{n\geq 1} P[ \bar{U}^n_s(T)^2 >a]=0.$   
 Let $\epsilon >0$ be arbitrary. We can pick $a>0$ so that    $P[ \bar{U}^n_s(T)^2 >a]< \epsilon$ for all $n.$   
 Using $\bar{U}^n_b(T)\leq |\bar{X}^n(0)| + \bar{A}^n_b(T), $ we have  
 \[
 E\left[ \bar{U}^n_b(T)^2 1_{\{\bar{U}^n_b(T)^2 >a\}} \right] \leq2 E[|\bar{X}^n(0)|^2] + 2E\left[ \bar{A}^n_b(T)^2  1_{\{\bar{U}^n_b(T)^2 >a\}}\right].
 \]  
From Assumption \ref{assump:initial},  $\lim_{n \rightarrow \infty}  E[|\bar{X}^n(0)|^2] =0.$
Next, we observe that    
  \begin{align*}
   E\left[ \bar{A}^n_b(T)^2  1_{\{\bar{U}^n_s(T)^2 >a\}}\right] & \leq  2 E\left[ |\bar{A}^n_b(T)-\lambda_0T|^2  1_{\{\bar{U}^n_b(T)^2 >a\}} + (\lambda_0 T)^2 P[\bar{U}^n_b(T)^2 >a]\right]\\          
   & \leq  2 E[ |\bar{A}^n_b(T)-\lambda_0T|^2   + (\lambda_0 T)^2 \epsilon.
   \end{align*}
Therefore,  we have 
 $ \limsup_{n \rightarrow \infty} E[\bar{U}^n_s(T)^2 1_{\{\bar{U}^n_b(T)^2 > a\}}] \leq  (\lambda_0 T)^2 \epsilon.$
  Next noting that $\bar U^n_b$ converges to zero in $D[0,T]$, by the bounded convergence theorem,   
\[\limsup_{n \rightarrow \infty} E\left[\bar{U}^n_b(T)^2 1_{\{\bar{U}^n_b(T)^2 \leq a\}}\right]=0.\]
Hence, $\lim_{n\to\infty}E[\bar U^n_b(T)^2] = 0.$ Similarly, $\lim_{n\to\infty}E[\bar U^n_s(T)^2] = 0.$ 

Since $ \bar{U}^n_b(t) $ and $ \bar{U}^n_s(t)$  are  non-decreasing processes,  \eqref{arrivals-M2} together with  \eqref{3.364} yields the conclusion  \eqref{arrivals-M4}. This completes the proof.
\Halmos

\begin{proposition} \label{prop-3.48}
For  the  process sequences  $\{\hat{L}^n_b\}_{n\in\NN}$ and   $\{\hat{L}^n_s\}_{n\in\NN}$ in $D[0, T]$, we have
 \begin{align}
&\lim\limits_{n \rightarrow \infty} \sup_{t\in[0,T]} \left| \hat{L}^n_b(t)- \sqrt{n} \int_{0}^{t} F_b(V^n_b(u-)) d \bar{A}^n_b(u) \right|=0, \ \text{in  probability,} \label{3.202}\\
&\lim\limits_{n \rightarrow \infty} \sup_{t\in[0,T]} \left| \hat{L}^n_s(t)- \sqrt{n} \int_{0}^{t} F_s(V^n_s(u-)) d \bar{A}^n_s(u) \right|=0, \ \text{in  probability.}
\label{3.204}
\end{align}
\end{proposition}

\proof
To prove \eqref{3.202}, we use the representation  \eqref{3.36} for $ \hat{L}^n_b(t)$ and first show that   
\begin{equation}
\lim\limits_{n \rightarrow \infty} \sup_{t\in[0,T]}\left| \hat{L}^n_b(t)- \sqrt{n} \int_{0}^{t} F_b(V^n_b(u-)) d \bar{E}^n_b(u) \right|=0,    \text{  in  probability. }
\label{3.40}
\end{equation}
Noting that  
$n^{-1/2} \sum_{j=1}^{E^n_b(t)} F_b(V^n_b(t^n_j-)) =  \sqrt{n} \int_{0}^{t} F_b(V^n_b(u-)) d \bar{E}^n_b(u),$ 
 \eqref{3.40}  will follow  from  \eqref{3.36} if we can show 
$ \lim_{n \rightarrow \infty} \sup_{t\in[0,T]}| M^n(\bar{E}^n_b(t))| =0$ in probability,  where the martingale $\{M^n(t)\}_{t\ge 0}$ is  described in  \eqref{3.36}. We intend to show that $ \lim_{n \rightarrow \infty} E[\sup_{t\in[0,T]}|M^n(\bar{E}^n_b(t))|^2]=0.$  Following the proof  of  Proposition \ref{prop-3.4}, we notice that $\{(M^n(t), \mathcal{G}^n_t )\}_{t\ge 0}$ is a square integrable, pure jump martingale and its quadratic variation process  is given by 
$$ [M^n, M^n](T)=  \frac{1}{n}  \sum_{j=1}^{[nT]}\left[ 1_{\{V^n_b(t^n_j-) \geq d^n_j\}}-F_b(V^n_b(t^n_j-))\right]^2.$$  
Moreover,  $\bar{E}^n_b(T)$ is a stopping time with respect to the filtration  $\{\mathcal{G}^n_t\}_{t \geq 0}.$  Using the Burkholder's inequality (see \cite{protter}),  we have   
\begin{align}\label{new3.56}
 E\left[\sup_{t\in[0,T]}|M^n(\bar{E}^n_b(t))|^2\right] \le E \left[ \sup\limits_{0\leq t \leq \bar{E}^n_b(T)} | M^n(t)|^2\right] \leq C_1 E( [M^n, M^n](\bar{E}^n_b(T))), 
 \end{align}
 where $C_1>0$ is a constant independent of $n, T$ {and the queue-capacity ${\bf m}_n$.}   Since $ [M^n, M^n](t) \leq 2t$  and   $\bar{E}^n_b(t) \leq  \bar{A}^n_b(t)$ for all $0 \leq t \leq T, $ we first pick a fixed  constant $ K > \lambda_0 +2$ such that  
\begin{equation}
  [M^n, M^n](\bar{E}^n_b(T)) \leq [M^n, M^n](\bar{E}^n_b(T))1_{\{ \bar{E}^n_b(T) <KT\}} +2  \bar{A}^n_b(T) 1_{\{ \bar{E}^n_b(T) \geq KT\}}. 
  \label{3.41}
  \end{equation}
 For the first term in the right hand of \eqref{3.41}, we can write 
\[ 
[M^n, M^n](\bar{E}^n_b(T))1_{\{\bar{E}^n_b(T) <KT\}}  =   \frac{1}{n}  \sum_{j=1}^{KnT}H^n_j 1_{\{ \bar{E}^n_b(T) <KT\}},
\]
where $H^n_j=  [1_{\{V^n_b(t^n_j-) \geq d^n_j\}}- F_b(V^n_b(t^n_j-))]^2 1_{\{t^n_j \leq T\}}.$  Since $d^n_j$ is independent of  $ V^n_b(t^n_j-)$ and $ t^n_j,$    
$E[ H^n_j  |  V^n_b(t^n_j-),  t^n_j ]=  F_b(V^n_b(t^n_j-)) (1- F_b(V^n_b(t^n_j-)))1_{\{t^n_j \leq T\}}.$  
Therefore, 
\begin{align*}
E \left([M^n, M^n](\bar{E}^n_b(T))1_{\{ \bar{E}^n_b(T) <KT\}}\right) &  \leq  \frac{1}{n}  \sum_{j=1}^{KnT} E(H^n_j)  \le \frac{1}{n}  \sum_{j=1}^{KnT} F_b(\|V^n_b\|_T) \\
&  = K T E[ F_b(\|V^n_b\|_T) ]. 
\end{align*}
By  the assumption \eqref{Patience-b}, there exists a constant $C_0>0$ so that  $F_b(x) \leq C_0 x$ for all $x\geq 0,$  and thus  
$E[ F_b(\|V^n_b\|_T) ]\leq C_0  E[ \|V^n_b\|_T] \leq \frac{C_0}{\sqrt{n}} E[ \|\hat{V}^n_b\|_T] \to 0, \ \mbox{as $n\to\infty$},$  
using the moment bound in \eqref{3.99}. Consequently,  
$\lim_{n \rightarrow \infty} E\left( [M^n, M^n](\bar{E}^n_b(T))1_{\{ \bar{E}^n_b(T) <KT\}} \right) =0.$
For the second term in the RHS of \eqref{3.41}, recalling that $ K > \lambda_0 +2,$ we observe that 
\begin{align*}
& E\left( \bar{A}^n_b(T) 1_{\{ \bar{E}^n_b(T) \geq KT\}} \right) \leq E\left( \bar{A}^n_b(T) 1_{\{ \bar{A}^n_b(T) \geq KT\}} \right)  \le E\left( \bar{A}^n_b(T) 1_{\{ \bar{A}^n_b(T) \geq \lambda_0 T + 2 T\}} \right) \\
& = E\left( |\bar{A}^n_b(T) -\lambda_0 T| 1_{\{ \bar{A}^n_b(T) \geq \lambda_0 T + 2 T\}} \right) + \lambda_0 T P\left( \bar{A}^n_b(T) \geq \lambda_0 T + 2 T\right)\\
& \le E[ | \bar{A}^n_b(T)- \lambda_0 T|] +  \lambda_0 T P[ | \bar{A}^n_b(T)- \lambda_0 T| >2T ] \to 0,
 \end{align*} 
 where the last step follows by \eqref{arrivals-M22}. Hence, in view of \eqref{3.41} and \eqref{new3.56}, \eqref{3.40} is established. 

Next we establish  
\begin{equation}
\lim\limits_{n \rightarrow \infty}  \sqrt{n} \int_{0}^{T} F_b(V^n_b(u-)) d \bar{U}^n_b(u) =0,    \text{  in  probability. }
\label{3.42}
\end{equation}
Then \eqref{3.40} and \eqref{3.42} implies \eqref{3.202}. Let $\epsilon >0$ be arbitrary. Since $\{\hat{V}^n_b\}_{n\in\NN}$ is stochastically bounded, we pick  a large constant $M >1$ so that $ P [ \|\hat{V}^n_b\|_T \geq M] < \epsilon. $ To derive \eqref{3.42}, it suffices to consider the set $\{ \|\hat{V}^n_b\|_T \leq M\}.$  We observe that  
$\sqrt{n} F_b(V^n_b(t-) 1_{\{\|\hat{V}^n_b\|_T \leq M\}} \leq \sqrt{n} F_b\left({M}/{\sqrt{n}}\right).$  
By \eqref{Patience-b}, $\sqrt{n} F_b({M}/{\sqrt{n}})$ is a bounded sequence, say, bounded above by  $K_M>0.$  Let $\delta >0$ be arbitrary. Then 
\begin{align*}
 & P\left[ \sqrt{n} \int_{0}^{T} F_b(V^n_b(u-)) d \bar{U}^n_b(u) > \delta,  \|\hat{V}^n_b\|_T \leq M\right] \\
& \leq  P[ \sqrt{n} F_b({M}/{\sqrt{n}}) \bar{U}^n_b(T) > \delta ] \leq  P[ \bar{U}^n_b(T)  >  {\delta}/{K_M}] \to 0, 
 \end{align*}
where the convergence in the last step follows from \eqref{3.364}. 
Then \eqref{3.42} holds.  This completes the proof of \eqref{3.202}. The proof of \eqref{3.204} will be similar  and  is omitted. 
\Halmos
 
\begin{proposition} \label{prop-3.5}

The  process sequence  $\{(\hat{L}^n_b, \hat{L}^n_s)\}_{n\in\NN}$ is C-tight in $ D^2[0, T].$ 

\end{proposition}

\proof
Here we establish the tightness of  $\{\hat{L}^n_b\}$  in $D[0, T].$  The proof for $\{\hat{L}^n_s\}$ is similar.  Let $ \omega(f, \delta, T)$ be as in  \eqref{omega}. By  Proposition \ref{prop-3.4}, $\{\hat{L}^n_b\}$ is stochastically bounded, hence  it suffices to  prove  $\lim_{n \rightarrow \infty}  P[ \omega(\hat{L}^n_b, \delta, T) > \epsilon ] =0$ for any $\epsilon >0$ to obtain the desired C-tightness.

Let  $\epsilon >0 $ be arbitrary.  Since $\{\hat{V}^n_b\}$ is stochastically bounded,  for  a  given constant  $M >0,$ there is  an integer $n_M$ so that  $P[ \| \hat{V}^n_b\|_T \geq M] < \epsilon$ whenever $n \geq n_M.$   By \eqref{Patience-b},    $\{\sqrt{n} F_b({M}/{\sqrt{n}})\}_{n\in\NN}$ is a bounded sequence, say, bounded above by  $K_M>0.$  Since  $\bar{E}^n_b(\cdot)$ is a non-negative  increasing process,  then on the set $[ \| \hat{V}^n_b\|_T \leq M],$ we have       $\sqrt{n} \int_{s}^{t} F_b(V^n_b(u-)) d \bar{E}^n_b(u) \leq K_M( \bar{E}^n_b(t)- \bar{E}^n_b(s))$   for every $0 \leq s \leq t \leq T.$   
Using this estimate together with  \eqref{3.36} on the set $[ \| \hat{V}^n_b\|_T \leq M]$ and when $ s, t $ are in $[0, T]$,  we obtain 
$|\hat{L}^n_b(t) - \hat{L}^n_b(s) | \leq  \sup_{t\in [0, T]}| M^n(\bar{E}^n_b(t))| + K_M  | \bar{E}^n_b(t)- \bar{E}^n_b(s))|.$  
Therefore, when $n\ge n_M$, 
 \begin{equation}
 \begin{aligned}
&P[ \omega(\hat{L}^n_b, \delta, T) > 2 \epsilon ] \leq  P[ \omega(\hat{L}^n_b, \delta, T) > 2 \epsilon,  \| \hat{V}^n_b\|_T \leq M ] + P[  \| \hat{V}^n_b\|_T > M ] \\
& \leq  P\left[ \sup_{t\in[0,T]}| M^n(\bar{E}^n_b(t))| > \epsilon \right] +   P[ \omega(\bar{E}^n_b, \delta, T) > { \epsilon}/{K_M} ]  +\epsilon.
 \end{aligned}
\label{3.44}
\end{equation}
We know  that  $ \lim\limits_{n \rightarrow \infty} E[\sup_{t\in[0,T]}| M^n(\bar{E}^n_b(t))|^2 ]=0$, which was established above \eqref{3.42},  and   $\lim_{n \rightarrow \infty} \sup_{t\in[0,T]}| \bar{E}^n_b(t) - \lambda_0 t|= 0$ in probability from  \eqref{arrivals-M4}.  
Hence $P[ \omega(\hat{L}^n_b, \delta, T) > 2 \epsilon ]$ in \eqref{3.44} converges to $0$ as $n\to\infty$. 
This completes the proof.
\Halmos

In the next proposition, we prove the C-tightness of   $\{(\hat{G}^n_b, \hat{G}^n_s)\}_{n\in\NN}$ in $D^2[0, T].$

\begin{proposition} \label{prop-3.6}
The  process sequence $\{(\hat{G}^n_b, \hat{G}^n_s)\}_{n\in\NN}$ is $C$-tight in $D^2[0, T]$, and furthermore, $\lim_{n \rightarrow \infty} \| \hat{L}^n_b - \hat{G}^n_b \|_T = 0 $ and $\lim\limits_{n \rightarrow \infty} \| \hat{L}^n_s - \hat{G}^n_s \|_T = 0 $ in probability. Consequently, 
 \begin{align}
&\lim\limits_{n \rightarrow \infty} \sup_{t\in[0,T]}\left| \hat{G}^n_s(t)- \sqrt{n} \int_{0}^{t} F_s(V^n_s(u-)) d \bar{A}^n_s(u) \right|=0, \ \text{in probability,}\label{3.48}\\
&\lim\limits_{n \rightarrow \infty} \sup_{t\in[0,T]}\left| \hat{G}^n_b(t)- \sqrt{n} \int_{0}^{t} F_b(V^n_b(u-)) d \bar{A}^n_b(u) \right|=0, \ \text{in probability.} \label{3.50}
\end{align}


\end{proposition}
\proof
The proof of this proposition is similar to the discussion in Section 4.4 of \cite{dai} and with appropriate changes, one can also easily follow the proof of Proposition 4.3 in \cite{weera3}. Moreover, \eqref{3.48} and \eqref{3.50} follows directly by combining the first part of the proposition with  \eqref{3.202} and \eqref{3.204} of  Proposition \ref{prop-3.48}. Hence, it will be omitted. 
\Halmos

\subsection{Asymptotics for \texorpdfstring{$\hat X^n$}{lg} and \texorpdfstring{$(\hat{G}^n_b, \hat{G}^n_s)$}{lg} }\label{subsec:tightX} 

To obtain the tightness of the sequence of state processes  $\{\hat{X}^n\}_{n\in\NN},$ we are heavily dependent on the properties of the two-sided Skorokhod map  in $D[0, T]$  as described in  \cite{burdzy}. In particular, we make use of an oscillation inequality for the Skorokhod map in a time varying interval. The inequality is given in Proposition \ref{prop-3.7} in Appendix \ref{sec:sm} together with other results on the Skorokhod map. 

\begin{theorem}
The processes  $\{(\hat{X}^n, \hat{U}^n_b, \hat{U}^n_s)\}_{n\in\NN}$ are  C-tight in $D^3[0, T].$

\label{thm3.1}
\end{theorem}

\proof
To prove the C-tightness of  $\{\hat{X}^n\}$ in $D[0, T],$  we will verify the two conditions of Theorem 13.2  in \cite{billingsley}. Proposition \ref{prop-3.1} implies the first condition on stochastic boundedness.  To verify the second condition, by  \eqref{3.1}--\eqref{3.4a},  notice that  $(\hat{X}^n, \hat{U}^n)$ is the unique solution to the SP in $D[0, T]$ for the input process $ \hat{\zeta}^n -\hat{G}^n$ on the time-dependent  interval  $[\hat{l}^n,   \hat{r}^n],$   where $\hat{l}^n(t)=- {m^n_b(t)}/{\sqrt{n}}$ and  $\hat{r}^n(t)={m^n_s(t)}/{\sqrt{n}}.$ Since $\hat{l}^n(t)\leq -\frac{1}{\sqrt{n}} \leq \frac{1}{\sqrt{n}} \leq \hat{r}^n(t) $ for all  $t \geq 0,$ we can combine Corollary 2.4 and Theorem 2.6  of \cite{burdzy}  to obtain the above unique solution. Hence by \eqref{3.60},
\begin{equation}
 \omega(\hat{X}^n,  \delta, T) \leq 4[ \omega(\hat{\zeta}^n-\hat{G}^n, \delta, T) + \omega(\hat{l}^n, \delta, T) + \omega(\hat{r}^n, \delta, T)].
\label{3.74}
\end{equation}

Since $\hat{\zeta}^n$  converges weakly to a process with continuous paths as in \eqref{3.4}, and $\{\hat{G}^n\}$ is C-tight from Proposition \ref{prop-3.6}, using Corollary 3.33 of Chapter VI in \cite{jacod}, $\{\hat{\xi}^n -\hat{G}^n\}$ is C-tight, and it follows that  $ \lim_{n \rightarrow \infty} P[  \omega(\hat{\xi}_n-\hat{G}_n, \delta, T) > \epsilon ]=0. $
 Next, using \eqref{cap-3}, we have
$ \lim_{\delta \rightarrow 0} \limsup_{n \rightarrow \infty}  \omega(\hat{l}^n, \delta, T)=0, \ \mbox{and} \ \ \lim_{\delta \rightarrow 0} \limsup_{n \rightarrow \infty}  \omega(\hat{r}^n, \delta, T)=0 \quad \mbox{in probability.}$ 
Using these facts together with \eqref{3.74}, we conclude   $\{\hat{X}^n\}$ is C-tight in $D[0, T].$

Now by \eqref{3.1},  again from Corollary 3.33 of Chapter VI in \cite{jacod}, we conclude that $\{\hat{U}^n\}$ is C-tight in $D[0, T].$
 Recall that  $ \hat{U}^n(t)= \hat{U}^n_s(t) -\hat{U}^n_b(t)$ for all $t \geq 0,$    and    $\hat{U}^n_s$ and  $\hat{U}^n_b$  are  non-decreasing processes.    We would like to show that  the processes $\hat{U}^n_s$ and  $\hat{U}^n_b$ are also tight in $D[0, T].$  We prove it for $\{\hat{U}^n_b\}$ and the proof for  $\{\hat{U}^n_s\}$ is similar.   
We know $\hat{U}^n_b(0) = [\hat m_b^n(0) - \hat X^n(0-)]^+$ is bounded (see Assumptions \ref{assump:initial} and \ref{assump:admiss}). By  \eqref{7.291} of Corollary \ref{cor:3}, we have 
 $\omega( \hat{U}^n_b, \delta, T) \leq C[ \omega( \hat{\xi}_n, \delta, T) +\omega( \hat{l}^n, \delta, T)]$ 
 and thus    
$ \lim_{\delta \rightarrow 0} \limsup_{n \rightarrow \infty}  \omega( \hat{U}^n_b, \delta, T)=0$ {in probability}.
 Now we  can follow  Theorem 13.2 and its corollary  in Chapter 3,  page 140 of \cite{billingsley}  (or see the discussion underneath Theorem 3.2 of \cite{whitt2}) to conclude $\{\hat{U}^n_b\}$ is C-tight in $D[0, T]$.
 \Halmos

In the following result, we build {an asymptotic little's law for the queue length process and the waiting time process, and the asymptotic linear} relationship between the number of abandonments  in $[0, t]$ and the integral of the queue length on $[0, t].$ The proof is related to similar results in \cite{dai} and \cite{weera3}.

\begin{proposition} 
\label{prop-3.9}
Let $T>0.$ We have the following convergence results. 
\begin{enumerate}
    \item[\rm (i)] {Asymptotic little's law.} \begin{align}
    & \lim\limits_{n \rightarrow \infty}E [\|\hat{X}^{n,+} - \lambda_0 \hat{V}^n_s \|_T] =0,\  \lim\limits_{n \rightarrow \infty}E[ \|\hat{X}^{n,-} - \lambda_0 \hat{V}^n_b\|_T] =0,
    \end{align}
    \item[\rm (ii)] {Asymptotic linear relationship between the abandonment and the integral of queue length.}
     \begin{align}
    & \lim\limits_{n \rightarrow \infty}E \left[ \sup_{t\in[0,T]}\left| \hat{G}^n_s(t)- \delta_s \int_{0}^{t} (\hat{X}^n(s))^+ ds \right|\right]=0,  
\label{3.80}\\
& \lim\limits_{n \rightarrow \infty} E \left[\sup_{t\in[0,T]}\left| \hat{G}^n_b(t)- \delta_b \int_{0}^{t}( \hat{X}^n(s))^- ds \right|\right]=0, 
\label{3.82}
\end{align}
and consequently, 
\begin{equation}
\lim\limits_{n \rightarrow \infty}E\left[\sup_{t\in[0,T]}\left| \hat{G}^n(t)-  \int_{0}^{t}[ \delta_s (\hat{X}^n(s))^+ - \delta_b ( \hat{X}^n(s))^- ]ds\right|\right]=0.  
\label{3.82-1}
\end{equation}
\end{enumerate}

\end{proposition}
\proof
The idea of the proof of part (i) is similar to that of Theorem 4.5 (ii)  in \cite{weera3}. Any seller arrived by time $t$ will be served by the time $t+V^n_s(t)$ unless the seller has abandoned the system. Let $K^n_s(t, t+r)$ be the number of sellers arrived after time $t$, but abandoned the system  by time $t+r$ and let $\hat{K}^n_s(t, t+r)= {K^n_s(t, t+r)}/{\sqrt{n}}.$ Therefore,
 \begin{align*} 
 (X^n(t+V^n_s(t)))^+ & =  [A^n_s(t+V^n_s(t)) -A^n_s(t) ] -[ U^n_s(t+V^n_s(t)) -U^n_s(t) ]\\
 & \quad - K^n(t, t+V^n_s(t)),
 \end{align*}
 and
 \begin{align*}
 (\hat{X}^n(t+V^n_s(t)))^+& =  \lambda_0 \hat{V}^n_s(t) + [\hat{A}^n_s(t+V^n_s(t)) -\hat{A}^n_s(t) ] -[ \hat{U}^n_s(t+V^n_s(t)) -\hat{U}^n_s(t) ]\\
 &\quad  - \hat{K}^n(t, t+V^n_s(t)).
 \end{align*}  
 Noting that    $0 \leq  \hat{K}^n(t, t+V^n_s(t)) \leq [\hat{L}^n_s(t+V^n_s(t))- \hat{L}^n_s(t)]$,  we have 
 \begin{equation}
 \begin{aligned}
  & \|(\hat X^n(t+V^n_s(t)))^+ - \lambda_0 \hat{V}^n_s(t) \|_T\\
  &  \leq \| \hat{A}^n_s(t+V^n_s(t))-\hat{A}^n_s(t)\|_T +\| \hat{U}^n_s(t+V^n_s(t))-\hat{U}^n_s(t)\|_T\\
  &\quad + \| \hat{L}^n_s(t+V^n_s(t))- \hat{L}^n_s(t)\|_T. 
 \end{aligned}
 \label{3.78}
 \end{equation} 
 By \eqref{3.99}, $\lim_{n \rightarrow \infty} \|V^n_s\|_T =0$ in probability. Since $\hat A^n_s, \hat U^n_s, \hat L^n_s$ and $\hat{X}^n$ are all C-tight  in $D[0, T]$, the RHS of  \eqref{3.78} tends to zero in probability. Hence $ \lim_{n \rightarrow \infty} \|\hat{X}^{n,+} - \lambda_0 \hat{V}^n_s\|_T =0 $ in probability. Consequently, using \eqref{3.6} and \eqref{3.99},   $ \lim_{n \rightarrow \infty}E [\|\hat{X}^{n,+} - \lambda_0 \hat{V}^n_s\|_T] =0.$ The proof of the second result in part (i) is similar and is omitted.

We next establish \eqref{3.80} and the proof of \eqref{3.82} is similar.  The proof consists of several steps. First we show that 
 \begin{equation}
\lim\limits_{n \rightarrow \infty} \sup_{t\in[0,T]}\left| \hat{G}^n_s(t)- \sqrt{n} \delta_s\int_{0}^{t} V^n_s(u-) d \bar{A}^n_s(u) \right|=0  \quad \text{  in  probability. }
\label{3.84}
\end{equation}
With  \eqref{3.48} in hand, it suffices to show  
\begin{align}\label{new3.72}
\lim\limits_{n \rightarrow \infty} \sup_{t\in[0,T]}\left| \sqrt{n} \int_{0}^{t}[ F_s(V^n_s(u-)) - \delta_s V^n_s(u-)] d \bar{A}^n_s(u) \right|=0 \quad \mbox{in  probability}.
\end{align}
Let $\epsilon >0 $ be arbitrary and $I_{n, T}=\sup_{t\in[0,T]}| \sqrt{n} \int_{0}^{t}[ F_s(V^n_s(u-)) - \delta_s V^n_s(u-)] d \bar{A}^n_s(u)|.$  By  \eqref{Patience-s}, there is a $\delta>0$  so that  $|F_s(x)-\delta_s x | < \epsilon x $ whenever $0 < x <\delta.$ By \eqref{3.99},   $\{\hat{V}^n_s\}$ is stochastically bounded     in $ D[0, T]$  and hence we can find $M >0$ so that   $ P [ \|\hat{V}^n_s\|_T \geq M] < \epsilon $  for all $n.$ We pick $n_0 > 1$ so that   ${M}/{\sqrt{n_0}} < \delta.$  On the set  $\{ \|\hat{V}^n_s\|_T \leq M\},$ $ | F_s(V^n_s(u-)) - \delta_s V^n_s(u-)| \leq \epsilon V^n_s(u-) $ for all $ u \leq T$  when $n>n_0.$ Hence, when $n\ge n_0$, 
\[
E\left[I_{n,T} 1_{\{ \|\hat{V}^n_s\|_T \leq M\}}\right]\leq \epsilon E[ \|\hat{V}^n_s\|_T\|\bar{A}^n_s\|_T].
\] 
Since $\epsilon$ is arbitrary, \eqref{arrivals-M22} and \eqref{3.99} together with the Holder's inequality yield $\limsup_{n \rightarrow \infty} E[I_{n,T} 1_{\{ \|\hat{V}^n_s\|_T \leq M\}}]=0.$ Next, by  \eqref{Patience-s} there is a $K>\delta_s$ so that  $F_s(x)\leq K x$ for all $x.$  Hence, 
\[
I_{n,T} 1_{\{\|\hat{V}^n_s\|_T> M\}} \leq K \|\hat{V}^n_s\|_T 1_{\{ \|\hat{V}^n_s\|_T > M\}} \|\bar{A}^n_s\|_T.
\]
By the Holder's inequality, \eqref{arrivals-M22} and the uniform integrability of $\|\hat{V}^n_s\|_T$ from \eqref{3.99}, we obtain $\limsup_{n \rightarrow \infty} E[I_n 1_{\{ \|\hat{V}^n_s\|_T>M\}}]=0.$ Hence, \eqref{new3.72} follows.
Next, from part (i), we observe that 
\begin{equation}
\lim\limits_{n \rightarrow \infty} \sup_{t\in[0,T]} \left| \int_{0}^{t} (\hat{X}^n(u))^+ -\lambda_0 \hat{V}^n_s(u-)] d \bar{A}^n_s(u) \right|=0  \quad \text{  in  probability.  }
\label{3.86}
\end{equation}
Finally, we need to establish
\begin{equation}
\lim\limits_{n \rightarrow \infty} \sup_{t\in[0,T]} \left| \int_{0}^{t} (\hat{X}^n(u))^+  (d \bar{A}^n_s(u)- \lambda_0 du) \right|=0  \quad \text{  in  probability.  }
\label{3.88}
\end{equation}
We can use Lemma 4.6 of \cite{weera3} with the use of Lipschitz function $g(x)=x^+$ and $T_n(t)=\bar{A}^n_s(t)- \lambda_0t $ to obtain \eqref{3.88}. Now combining \eqref{3.84}, \eqref{3.86}  and  \eqref{3.88},  
\begin{align}\label{in-prob-conv}
\sup_{t\in[0,T]}\left| \hat{G}^n_s(t)- \delta_s \int_{0}^{t} (\hat{X}^n(s))^+ ds \right|\to 0, \ \ \mbox{in probability}.
\end{align}
Finally, to obtain \eqref{3.80}, we establish the uniform integrability of the LHS of \eqref{in-prob-conv}. 
Using  \eqref{3.6} and  \eqref{3.22}, 
{\[
E\left[\sup_{t\in[0,T]}| \hat{G}^n_s(t)- \delta_s \int_{0}^{t} (\hat{X}^n(s))^+ ds |^2 \right]\leq  2E[(\hat{G}^n_s(T))^2 + (\delta_s T  \|\hat{X}^n\|_T)^2] \le  C(1+T^m),
\]}
where $C, m$ are positive constants independent of $n$ and $T.$ Now \eqref{3.80}  follows. The proof of \eqref{3.82} is similar, and \eqref{3.82-1} follows from  \eqref{3.80} and \eqref{3.82}.  This completes the proof.
\Halmos

\subsection{Proof of Theorem \ref{thm3.2}}\label{subsec:lastproof} 


Part (i) follows from \eqref{3.4} in Remark \ref{rem:wc}. Introduce $\hat{Z}^n= (\hat{X}^n,  \hat{\zeta}^n, \hat{G}^n_b, \hat{G}^n_s, \hat{U}^n_b,  \hat{U}^n_s)$ in  $D^6[0, T].$   Then by  \eqref{3.4},  Proposition \ref{prop-3.6}  and Theorem \ref{thm3.1}, each component of $\hat Z^n$ is C-tight. Using Corollary 3.33 of Chapter VI in \cite{jacod}, it follows that the sequence  $\{\hat{Z}^n\}$ is C-tight in  $D^6[0, T].$ 

Let $Z = ({X},  {\zeta}, {G}_b, {G}_s, {U}_b,  {U}_s)$ be any limit point of  $\{\hat{Z}^n\}$ along a subsequence. Without loss of generality, we relabel the subsequence such that $\hat Z^n$ converges weakly to $Z$ as $n\to\infty$. The C-tightness guarantees the continuity of paths for $Z.$  Then, using the Skorokhod representation theorem (see \cite{kurtz} Chapter 3, Theorem 1.8), we simply assume that  $ \lim_{n \rightarrow \infty} \hat{Z}^n(t)= Z(t) \text{ a.s.} $ uniformly  on $[0, T].$  Moreover, $U_b, U_s$  are non-decreasing. By \eqref{3.82}, and the continuous mapping theorem, part (ii) follows. For part (iii), using \eqref{3.82}, we can write 
\begin{equation}
\hat{X}^n(t)=  \hat{\zeta}^n(t) -  \int_{0}^{t}h(\hat{X}^n(s))ds -  (\hat{U}^n_s(t) - \hat{U}^n_b(t))  +\epsilon_n(t),
\label{3.92}
\end{equation}
where  $h(x)= \delta_s x^+ - \delta_b x^- $  and  $ \lim_{n \rightarrow \infty} \| \epsilon_n\|_T = 0$ in probability. The processes $\hat{U}^n_s$ and  $\hat{U}^n_b$ are monotone increasing,  and with disjoint increment supports on the sets $\{X^n(s)> 0\}$  and  $\{X^n(s)< 0\}$, respectively.  Hence, $X$ satisfies
$X(t)=x +\sigma B(t) +  \int_{0}^{t}[\beta - h(X_x(s))]ds -U(t)$
 for each $0 \leq t \leq T.$  The constants $\sigma$ and $\beta$ are as described in  \eqref{3.4}.  The process $U$ is  of bounded variation in $D[0, T]$ adapted to the filtration generated by $(X, B)$  and using Proposition 2.3 of  \cite{burdzy}, it can be expressed as 
$U(t)=U_s(t)- U_b(t),$
where 
$U_s(t)= \int_{0}^{t} 1_{\{X(r)>0\}} dU_s(r) $ and   $U_b(t)= \int_{0}^{t} 1_{\{X(s)<0\}} dU_b(s). $  This completes the proof of part (iii).

\section{Diffusion Control Problem (DCP)}\label{sec:dcp}

\subsection{Problem Formulation}

We formulate  a one-dimensional stochastic control problem for diffusion processes (i.e., the DCP)  which can be considered as the
limiting form of the cost minimization problem for the queueing
systems (i.e., the QCP). An explicit solution for this DCP will be obtained here.  In Section \ref{sec:as}, we shall ``translate" the optimal strategy
of the DCP to obtain an asymptotically optimal strategy for the QCP.

We consider the limiting process of $(\hat X^n, \hat U^n)$ derived in Theorem \ref{thm3.2}. Rigorously, we define a controlled state-process
\begin{equation}
X_x(t)=x +\sigma B(t) +  \int_{0}^{t}[\beta - h(X_x(s))]ds -U(t),
\label{4.2}
\end{equation}
where $X_x(0)=x$ is a real number as in \eqref{initial}, the parameters $\sigma$ and $\beta$ are constants as in Theorem \ref{thm3.2} (i), $B$ is a standard one-dimensional
Brownian motion, adapted to a right-continuous filtration
$\mathcal{F}=\{\mathcal{F}_t : t \geq 0 \}$ on a probability space
$(\Omega, \mathfrak{F}, \mathbf{P})$, and the function 
\begin{equation}
h(x) = \delta_s x^+ - \delta_b x^-,
 \label{4.1}
\end{equation}
with $ \delta_b$ and $ \delta_s$ being  positive constants as in \eqref{Patience-s} and \eqref{Patience-b}. Furthermore, the $\sigma$-algebra
$\mathcal{F}_0$ contains all the null sets in $\mathfrak{F}$, the
Brownian increments $B(t+s)-B(t)$ are independent of $\mathcal{F}_t$
for all $t \geq 0$ and $s \geq 0$, and the control $U$ is a
 right-continuous process with paths of bounded variation  adapted to the
filtration  $\mathcal{F}.$  It is assumed that  $U$ can be expressed as 
\begin{equation}
U(t)=U_s(t)- U_b(t),
\label{control-1}
\end{equation}
where 
$U_s(t)= \int_{0}^{t} 1_{\{X(r)>0\}} dU_s(r) $ and   $    U_b(t)= \int_{0}^{t} 1_{\{X(s)<0\}} dU_b(s)$, and the processes $U_b$ and $U_s$ are adapted, non-decreasing  with RCLL paths. Thus no control can be enforced when the state process is at the origin where there is no holding cost in the cost structure. 
Throughout this section, the pair $(U_s, U_b)$ describes the control policy.

 We introduce the following cost functional for the state process in (\ref{4.2}):
\begin{equation}
J(X_x,  U_s, U_b)= E\left(\int_{0}^{\infty}e^{-\alpha t}[ (\theta_s X^+_x(t)+ \theta_b
X^-_x(t))dt + p_s\hspace{.5mm} dU_s(t) + p_b\hspace{.5mm} dU_b(t)]\right) , 
\label{cost-fn}
\end{equation}
{where $\theta_s = c_s + r_s\delta_s, \theta_b = c_b + r_b \delta_b$ and $\alpha$, $p_s$, and $p_b$ as well as $c_s, r_s, \delta_s, c_b, r_b, \delta_b$ are positive constants as given in the QCP \eqref{cost-1}.} Note that the holding cost 
\begin{align}\label{dcp-holding}
C(x) \equiv \theta_s x^+ + \theta_b x^-, x\in\rr,
\end{align} 
is  a piecewise linear  convex function. 

For $x\in\mathbb{R}$, we call the  sextuple $\left
((\Omega, \mathfrak{F}, \mathbf{P}),\mathcal{F}, B, X_x, U_s, U_b \right )$ \textit{an
admissible control system} if
\begin{itemize}
\item[(i)] $(X_x, U_s, U_b)$ is a weak solution to  \eqref{4.2}, and
\item[(ii)] the cost functional $J(X_x,  U_s, U_b)$ is finite.
\end{itemize}
When there is no ambiguity, 
we simply use  $(X_x,  U_s, U_b)$ to represent an admissible control
system.   To define the value
function, we introduce the set
\begin{equation}
\mathcal{A}(x)=\left \{ (X_x,  U_s, U_b) : (X_x,  U_s, U_b)\hspace{0.2cm} \textrm{is admissible}
\right \}. \label{ctrl-set}
\end{equation}
This set is nonempty for each $x$ in $\rr$ since the zero control policy $U_s=U_b=0$ leads to a diffusion process $X_x$ for which  $J(X_x, 0, 0)$ is
finite. The value function of the DCP is thus well defined and given by
\begin{equation}
V(x)=\inf\limits_{\mathcal{A}(x)}J(X_x,  U_s, U_b), \ x\in \rr. 
\label{value-fn}
\end{equation}

At last we describe the formal HJB equation
associated with our DCP. We introduce the differential operator
$\mathcal{G}$ by
\begin{equation}
\mathcal{G}= \frac{\sigma^2}{2}\frac{d^2}{dx^2} +(\beta - h(x))
\frac{d}{dx}-\alpha, 
\label{4.6}
\end{equation}
where the constants $\alpha$, $\beta$, and $ \sigma$ and the function $h(x)$ are as described earlier in \eqref{4.2}.
The formal HJB-equation for the above control problem can now be written as
\begin{equation}
\min \left \{\mathcal{G}F(x) +C(x), \  F'(x)+p_b   ,\  p_s-F'(x) \right \}=0, \label{4.8}
\end{equation}
where $C(x)$ is the holding cost function \eqref{dcp-holding}.
Our proofs show that the value function $V$ is the unique smooth solution to the HJB equation \eqref{4.8}.

\paragraph{A verification lemma.}
The following verification lemma  guarantees that any smooth function which 
satisfies  the HJB equation \eqref{4.8} is a lower bound for the
value function and helps us to identify an
optimal strategy. 
\begin{lemma} \label{lm4.1}
Assume that $F\in C^2(\rr)$  and satisfies  \eqref{4.8} on 
$ \rr.$    Then $F(x) \leq V(x)$ for all $x\in \rr$, where $V$ is the value
function defined in (\ref{value-fn}).
\end{lemma}
\proof 
Let $(X_x,U_s, U_b)$ be an admissible control system, and $F\in C^2(\rr)$ satisfy  \eqref{4.8} on 
$ \rr.$   Using the   generalized It$\hat{o}$'s lemma   (see \cite{meyer}, p.\ 285), 
\begin{equation}
\begin{aligned}
F(X_x(T))e^{-\alpha T}& =F(x)+\sigma \int_{0}^{T}e^{-\alpha s
}F'(X_x(s-))dB(s) \\
& \quad + \int_{0}^{T}e^{-\alpha s
}\mathcal{G}F(X_x(s-))ds- \int_{0}^{T}e^{-\alpha s
}F'(X_x(s-))dU(s) \\
&\quad + \sum\limits_{0 \leq s \leq T} e^{- \alpha
s}[\Delta F(X_x(s))+ F'(X_x(s-))\Delta U(s)],
\end{aligned} 
\label{4.10}
\end{equation}
where $\Delta F(X_x(s))=F(X_x(s))-F(X_x(s-))$, and $\Delta
U(s)=U(s)-U(s-)$.
We let  $U^c_b$, $U^c_s$, and $U^c = U^c_s - U^c_b$ be the continuous parts of the processes $U_b$, $U_s$ and $U$, respectively.  By \eqref{4.8}, $F'(x) \in [-p_b, p_s]$ and hence  $E[\int_{0}^{T}e^{-\alpha s
}F'(X_x(s-))dB(s)] =0$ for each $T.$ Then by \eqref{4.10}, we obtain 
\begin{equation}
\begin{aligned}
& E\left[F(X_x(T))e^{-\alpha T}\right]  = F(x)+E \left[ \int_{0}^{T}e^{-\alpha s
}\mathcal{G}F(X_x(s-))ds\right] \\
& - E\left[ \int_{0}^{T}e^{-\alpha s
}F'(X_x(s-))dU^c(s)\right] 
+ E\left[\sum\limits_{0 \leq s \leq T} e^{- \alpha
s} \Delta F(X_x(s))\right].
\end{aligned} 
\label{4.12}
\end{equation}

When  $X_x(t)>0$  for some $t>0,$    $\Delta X_x(t)=-\Delta U_s(t)$ and   $\Delta F(X_x(t))= F'(\xi) \Delta X_x(t),$ where   $\xi$ is between  $X_x(t)$ and $  X_x(t-).$   Since $F'(x) \leq p_s$, we obtain  
$\Delta F(X_x(t))\geq -p_s \Delta U_s(t)= -[ p_s \Delta U_s(t) + p_b \Delta U_b(t)],$  
where $\Delta U_b(t)=0$ when  $X_x(t)>0.$
By  a similar argument  when $X_x(t) <0$, 
 $\Delta F(X_x(t))\geq -[ p_s \Delta U_s(t) + p_b \Delta U_b(t)].$ 
Moreover, 
$F'(X_x(t-)) dU^c(t)=F'(X_x(t-))(dU^c_s(t)-dU^c_b(t)) \leq p_sdU^c_s(t) + p_b dU^c_b(t).$ By \eqref{4.8},  $\mathcal{G}F(x) \ge -C(x)$ holds for each $x\in\rr$, and hence from \eqref{4.12},  
\begin{align*} 
E\left[F(X_x(T))e^{-\alpha T}\right]& \ge  F(x)-E \left[ \int_{0}^{T}e^{-\alpha s}  C(x(s))ds\right] \\
& \quad - E\left[ \int_{0}^{T}e^{-\alpha s} (p_s dU_s(s) +p_b dU_b(s))\right]. 
\end{align*}
Consequently,  $ E[F(X_x(T))e^{-\alpha T}] + J(X_x, U_s, U_b) \geq F(x)$ for each $x$ and $T>0.$ Since $F'(\cdot)$ is bounded, there exist some $c_1, c_2 >0$, 
$E[|F(X_x(T))|e^{-\alpha T}] \leq (c_1 + c_2 E[|X_x(T)|) e^{-\alpha T}.$ 
 Noting that $xh(x) \geq 0$ and  $X(t) dU(t) \geq 0,$  we can apply It$\hat{o}$'s lemma for $X_x(t)^2$ to obtain  $\lim_{T \rightarrow \infty} e^{- \alpha T} E[|X_x(T)|^2] =0.$   Hence $\lim_{T \rightarrow \infty} e^{- \alpha T}E |F(X_x(T))| =0.$  Therefore,  $ J(X_x,  U_s, U_b) \geq F(x)$  for each $x.$  Consequently,  $V(x) \geq F(x)$  holds for each $x.$ This completes the proof.
\Halmos

\subsection{Optimal solutions of the DCP}
When the blocking costs $p_s$ and  $p_b$ are high relative to the holding and abandonment costs, blocking customers is not cost effective, whereas if the costs $p_s$ and  $p_b$ are sufficiently low, then it is optimal to use finite buffer sizes. Here we clarify the threshold values for $p_s$ and  $p_b$ and establish a threshold optimal solution to the DCP under four different regimes. 

\begin{theorem}\label{thm:solution-DCP}
Let $T_s = \theta_s/(\alpha + \delta_s)$ and $T_b = \theta_b/(\alpha + \delta_b)$. Then the DCP admits the following optimal solution. 
\begin{itemize}
    \item[\rm (i)] When $p_s \geq T_s$ and  $p_b \geq T_b$, an optimal strategy is given by the zero control policy, which is described by  $U^*_s(t)=U^*_b(t)=0$ for all $t,$ with the state equation 
\begin{equation}
X^*_x(t)=x +\sigma B(t) +  \int_{0}^{t}[\beta - h(X^*_x(s))]ds, \quad t\ge 0.
 \label{4.26}
\end{equation}
\item[\rm (ii)] When $0<p_s < T_s$ and  $0<p_b < T_b$, there exist two points $a^* <0 <b^*$ such that the reflected diffusion process $X^*_x$  on  $[a^*,  b^*] $  described by 
\begin{equation}
X^*_x(t)=x+ \sigma B(t) + \int_{0}^{t} [\beta-h(X^*_x(s))]ds -L_{b^*}(t) +L_{a^*}(t)
\label{reflected}
\end{equation}  
is an optimal state process, and the optimal control pair $(U^*_s,  U^*_b)$ is given by the local-time processes  $(L_{b^*},  L_{a^*}).$ 

\item[\rm (iii)] 
 When $p_s \ge T_s$ and  $0<p_b < T_b$, there exists a point $\tilde a^* < 0$ such that the reflected diffusion process $X^*_x$  on  $[\tilde a^*,   \infty)$  described by 
\begin{equation}
X^*_x(t)=x+ \sigma B(t) + \int_{0}^{t} [\beta-h(X^*_x(s))]ds +L_{\tilde a^*}(t),
\label{left-reflected}
\end{equation}  
is an optimal state process, and the optimal control pair $(U^*_s,  U^*_b)$ is given by the local-time processes  $(0,  L_{\tilde a^*}).$ 
\item[\rm (iv)] 
When $0<p_s < T_s$ and  $p_b \ge T_b$, there exists a point $\tilde b^* > 0$ such that the reflected diffusion process $X^*_x$  on  $(-\infty, \tilde b^*]$  described by 
\begin{equation}
X^*_x(t)=x+ \sigma B(t) + \int_{0}^{t} [\beta-h(X^*_x(s))]ds -L_{\tilde b^*}(t)
\label{right-reflected}
\end{equation}  
is an optimal state process, and the optimal control pair $(U^*_s,  U^*_b)$ is given by the local-time processes  $(L_{\tilde b^*}, 0).$ 
\end{itemize}
\end{theorem}
\begin{remark}
In Theorem \ref{thm:solution-DCP} (ii), (iii) and (iv), if the initial value $x$ is not in the desired interval, there will be an initial jump to the nearest point in the interval. For example, in (ii), if $x$ is outside the interval $[a^*, b^*]$, then there is an initial jump to $a^*$ when $x<a^*$ or $b^*$ when $x>b^*$.
\end{remark}

The proof of Theorem \ref{thm:solution-DCP} relies on the construction of solutions to the HJB equation \eqref{4.8} with different boundary conditions. In the following Sections \ref{sec:zero} and \ref{sec:reflected}, we construct the solution of the HJB \eqref{4.8} and provide the proof of Theorem \ref{thm:solution-DCP}. At last, in Section \ref{sec:num}, we summarize how the optimal buffer sizes $a^*, b^*, \tilde a^*$ and $\tilde b^*$ are computed and present a numerical example. 

\subsubsection{Optimality of Zero Control}\label{sec:zero}

In this section we verify that  the zero control (no blocking of customers) strategy is optimal when $p_b$ and $p_s$ are above the threshold values $T_b$ and $T_s$. To this end, 
we need to construct  a solution  $Q$ to \eqref{4.8}  of the the form:
\begin{equation}
\mathcal{G}Q(x) +C(x)=0 \text{ and }  -p_b< Q'(x) <p_s, \quad x\in\rr.
\label{4.14}
\end{equation}
Next let  $\gamma(x)=h'(x)$ and it can be written in the form
\begin{equation}
\gamma(x)=\delta_b 1_{\{x<0\}}+ \delta_s 1_{\{x>0\}}, \quad x\in\rr\backslash\{0\}.
\label{4.16}
\end{equation}
We let $\W(x)=Q'(x)$ and differentiating $Q$, we obtain 
\begin{equation}
\mathcal{G}\W(x) - \gamma(x) \W(x)+C'(x)=0, \quad x\in\rr\backslash\{0\}.
\label{4.18}
\end{equation}
 First we observe that $W(x)=-T_b$ is a constant solution on $(-\infty, 0)$, and  $W(x)= T_s$ is a constant solution on $(0, \infty).$ They play an important role in our analysis.  Since
$\alpha Q(0)= \frac{\sigma^2}{2} \W'(0) +\beta \W(0),$   $Q$  can be obtained by knowing $\W.$ 
Therefore, we  construct a  sufficiently smooth bounded solution $\W$ to \eqref{4.18}, 
with the boundary data
\begin{equation}\label{bound-zero}
 \W(-\infty)=  -T_b \text{ and }  \W(\infty)=  T_s. 
\end{equation}

The following proposition presents the existence and uniqueness of such a function $W$. We defer its proof in Appendix \ref{sec:proofprop-R}.
\begin{proposition} 
\label{prop-W}
There exists a unique   function  $\W \in C^1(\rr)$  satisfying  \eqref{4.18} {and \eqref{bound-zero}} with the following properties.
\begin{enumerate}
\item[\rm (i)]  The function $\W$ is bounded and strictly increasing.
\item[\rm (ii)] $\W $ is twice continuously differentiable everywhere except at the origin. 
\item[\rm (iii)] The one sided limits $\W''(0-)$ and $\W''(0+)$ exist and  they satisfy   
\begin{align} 
\frac{\sigma^2}{2} \W''(0-) +\beta \W'(0)&=(\alpha + \delta_b) \W(0) + \theta_b,
\label{W-right}\\
 \frac{\sigma^2}{2} \W''(0+) +\beta \W'(0)&=(\alpha + \delta_s) \W(0) - \theta_s.
 \label{W-left}
\end{align}
\end{enumerate}
\end{proposition}

\proof[Proof of Theorem \ref{thm:solution-DCP} (i)]
Using the function $\W$ derived in Proposition \ref{prop-W},  we introduce the function $Q$ by 
$$\alpha Q(x)=\frac{\sigma^2}{2}\W'(0) + \beta \W(0) + \alpha \int_{0}^{x} \W(u) du, \quad x \in \rr.$$
Moreover,   $Q'(x) \equiv \W(x)$ is strictly increasing. Then it is straightforward to check that  $Q$ satisfies 
$\mathcal{G}Q(x) +C(x)=0,$ and $-T_b<Q'(x) < T_s. $

By the assumptions in Theorem \ref{thm:solution-DCP} (i), $-p_b < Q'(x) < p_s$ for all $x$ and thus $Q\in C^2(\rr)$  and is a  convex  function which  satisfies  \eqref{4.8}.  Therefore, by using the verification Lemma \ref{lm4.1},  we can conclude that $Q(x) \leq V(x)$ for all $x.$

Next we consider the state process  $X_x$ described in \eqref{4.26} which corresponds to zero control strategy. We apply   the It$\hat{o}$'s lemma to   $Q(X_x(t))e^{-\alpha t}$  and use the facts that  $Q'$ is bounded and  $Q$ satisfies the HJB equation  \eqref{4.8}  to obtain  
\begin{align}\label{soln-zerocontrol}
Q(X_x(T))e^{-\alpha T}  =Q(x) - E \left(\int_{0}^{T} e^{-\alpha s} C(X_x(s)) ds\right).
\end{align}
From the proof of Lemma  \ref{lm4.1}, we have $\lim_{T \rightarrow \infty} e^{- \alpha T} E[|X_x(T)|] =0$. Since $Q'(\cdot)$ is bounded, it follows that  
$\lim_{T \rightarrow \infty} e^{- \alpha T} E[Q(X_x(T))] =0.$  
Therefore, by letting $T$ tends to infinity in \eqref{soln-zerocontrol}, we have  
$Q(x)= E \left[\int_{0}^{\infty} e^{-\alpha s} C(X_x(s)) ds\right],$ which says $Q(x)$ represents the pay-off function from the zero control strategy. Hence $Q(x) \geq V(x).$ This completes the proof.
\Halmos

\subsubsection{Optimality of Reflected Diffusion Processes}\label{sec:reflected}




We first focus on the case when both costs $p_b$ and $p_s$ are below the threshold values $T_b$ and $T_s$. The following proposition lays the groundwork to derive the optimal policy in Theorem \ref{thm:solution-DCP} (ii), and its proof is presented in Appendix \ref{sec:proofprop-R}.

\begin{proposition} 
\label{prop-R}

Let the cost parameters $p_b$ and  $p_s$ satisfy $0< p_b < T_b $ and $0< p_s < T_s$ 
Then  there exist two points $a^* <0 <b^*$ and a function  $\W^*:[a^*,   b^*] \rightarrow \mathbf{R} $  satisfying the following free boundary problem: 
\begin{align}
& \mathcal{G}\W^*(x) - \gamma(x) \W^*(x)+C'(x)=0, \   a^*\le x \le b^*,   
\label{4.32} \\
& \W^*(a^*)=  -p_b,   \W^*(b^*)= p_s, \ \text{  and  } \ {\W^*}'(a^*)={\W^*}'(b^*)=0.
\label{4.34}
\end{align}
Moreover,  $-p_b <\W^*(x)<p_s$ when $a^*<x<b^*.$
\end{proposition}



\proof[Proof of Theorem \ref{thm:solution-DCP} (ii)]
We first extend the function  $\W^*$ in Proposition \ref{prop-R}  to $\rr $ by defining   $\W^*(x)=-p_b$ if $x<a^*$ and  $\W^*(x)=p_s$ if $ x>b^*.$  Hence, $ \W^* \in C^1(\rr)\cap C^2(\rr-\{a^*, b^*\}).$  Then introduce $Q\in C^2(\rr)$  by 
$$ 
\alpha Q(x)=\left[ \frac{\sigma^2}{2}{\W^*}'(0) + \beta \W^*(0)\right] + \alpha \int_{0}^{x} \W^*(u) du, \quad x \in \rr.
$$ 
To verify  $Q$ satisfies  \eqref{4.8}, we notice  $\W^*$ satisfies \eqref{4.32} in  {$[a^*,  b^*].$}  Moreover,  $\mathcal{G}Q(b^*) +C(b^*)=0$ and  
  for $x>b^*,$
\[[\mathcal{G}Q(x) +C(x)] -[ \mathcal{G}Q(b^*) +C(b^*)] =(\alpha+\delta_s)\left[ \frac{\theta_s}{\alpha + \delta_s} - p_s\right](x-b^*) >0.\]
Thus  $\mathcal{G}Q(x) +C(x) >0$ and    ${Q}'(x)=p_s$ on {$(b^*, \infty).$}  Consequently, $Q$ satisfies  \eqref{4.8} on $[b^*, \infty).$   Similarly, $Q$ satisfies  \eqref{4.8} on   $(-\infty,   a^*].$ Therefore, $Q\in C^2(\rr)$ satisfies the HJB equation \eqref{4.8}  and  by the verification Lemma \ref{lm4.1},  we obtain 
$Q(x) \leq V(x)\text{ for }  x\in\rr.$

Next, we show that   $X^*_x$  in \eqref{reflected}  yields the pay-off $Q(x).$  For this, first we consider $ x\in[a^*, b^*]$ and  apply It$\hat{o}$'s lemma to  $Q(X^*(t))e^{-\alpha t}$   
to obtain  
$$Q(x)= E\left[Q(X^*_x(T)) e^{-\alpha T}\right] +E\left(\int_{0}^{T} e^{-\alpha t} [C(X^*_x(t)) dt +p_b dL_{a^*}(t) + p_s dL_{b^*}(t)]\right).$$
Since $Q$ is bounded on $[a^*,   b^*]$, $\lim_{T\rightarrow \infty}  E[Q(X^*_x(T)) e^{-\alpha T}] =0.$ Hence, when $x\in[a^*, b^*],$  
$$Q(x)= E\left(\int_{0}^{\infty} e^{-\alpha t} [C(X_x^*(t)) dt +p_b dL_{a^*}(t) + p_s dL_{b^*}(t)]\right).$$ 
It is evident that  the process $(X^*_x,  L_{a^*},  L_{b^*})$  forms an admissible policy.  When $x>b^*$,  there is an initial jump of $X^*$ to $b^*$ so that $L_{b^*}(0+)= (x-b^*)$. Then it follows that $Q(x)=Q(b^*) +p_s (x-b^*).$   Similar analysis follows for the case $x<a^*.$ Consequently,  $Q(x) \geq V(x)$ holds and the process  $X^*_x$ together with the controls $L_{a^*}$ and $L_{b^*}$ describes an optimal policy. This completes the proof.
\Halmos

We now discuss the situations where one-sided reflected diffusion processes are optimal.


\proof[Proof of Theorem \ref{thm:solution-DCP} (iii) and (iv)]
To prove part (iii), we find a point $\tilde a^*<0$ and construct a  bounded strictly increasing function  $\W_L :[\tilde a^*, \infty) \rightarrow \mathbf{R}$ satisfying 
\begin{equation}
\begin{aligned}
& \mathcal{G}\W_L(x) - \gamma(x) \W_L(x)+C'(x)=0 \text{  for } x >\tilde a^*, \\
 & \W_L(\tilde a^*)=  -p_b, \quad \W'_L(\tilde a^*)=0,  \quad \W_L(\infty)=  T_s. 
\label{4.40}
\end{aligned}
\end{equation}
This function $\W_L(\cdot)$ will be obtained in the Appendix \ref{sec:proofprop-R} (see Remark \ref{rem:thm4.6}). Thereafter, we can essentially follow the proof of Theorem \ref{thm:solution-DCP} (ii), since this is essentially the case $b^*=\infty.$ Hence, we omit the details.
Proof of part (iv) is quite similar and is omitted. This completes the proof.
\Halmos

\subsubsection{Computing the optimal buffer sizes}\label{sec:num}
We first assume that both blocking cost parameters are below the thresholds, i.e., $0<p_s < T_s$ and $0<p_b < T_b$. We summarize the solution construction process in Appendix \ref{sec:proofprop-R} and find $a^*$ and $b^*$ in Theorem \ref{thm:solution-DCP} (ii). 
\begin{itemize}
    \item For $a<0$, we consider a function $W_a$ that satisfies the linear equation \eqref{4.18} on $(-\infty, 0],$ and $W_a(a) = -p_b$ and $W_a'(a)=0$. For each $a<0,$  we can extend  $W_a$ to $(0,  \infty)$ so that it satisfies  \eqref{4.18}  on  $(0, \infty)$  with the available initial data $W_a(0)$ and $W_a'(0).$ More precisely, $W_a$ satisfies  
\begin{align*}
& \frac{\sigma^2}{2} W''_a(x) + (\beta+\delta_b x) W'_a(x) -(\alpha +\delta_b)W_a(x) - \theta_b =0, \ x\le 0,\\
& \frac{\sigma^2}{2} W''_a(x) + (\beta-\delta_s x) W'_a(x) -(\alpha +\delta_s)W_a(x) + \theta_s =0, \ x>0,
\end{align*} 
with the given $W_a(a) = -p_b$ and $W_a'(a) =0$. 

\item We summarize some properties of $W_a(\cdot)$ for $ a < 0.$
\begin{itemize}
    \item When $|a|$ is sufficiently large, from Lemma \ref{lemA.20}, $\lim_{x\to\infty} W_a(x) = \infty$ and $W_a(\cdot)$ is strictly increasing on $(a, \infty).$
\item Let 
\begin{align}\label{c-value}
    c=\sup\{a<0: \lim_{x\to\infty} W_a(x) = \infty\}
\end{align} From Lemma \ref{lemA.22}, $c<0$ and $\lim_{x\to\infty}W_c(x) = T_s.$
\item From Lemma \ref{lemA.26}, for $a\in (c, 0)$, $\lim_{x\to\infty} W_a(x) = -\infty$ and $W_a(\cdot)$ has a unique maximum on $(0, \infty).$
\end{itemize}
\item For $a\in (c, 0)$, introduce 
\begin{align*}
    M(a) = \max_{x\ge a} W_a(x), \ \ r_a = \arg\max_{x\ge a} W_a(x).
\end{align*}
Then $a^*$ is chosen from the interval $(c, 0)$ such that $M(a^*) = p_s$ (this can be achieved because $p_s < T_s$) and $b^* = r_{a^*}.$
\item Let $W^*(t) = W_{a^*}(t), t\in \rr.$ Then $W^*$ satisfies all conditions in Proposition \ref{prop-R}.
\end{itemize}

Consider now the case when one cost parameter is below the threshold and the other one is above its threshold. In the aforementioned construction process, when $p_s$ increases and approaches the threshold $T_s$, the corresponding $a^*$ will approach $c$, and the $b^*$ will grow to infinity. More precisely, when $0<p_b < T_b$ and $p_s \ge T_s$, we have $\tilde a^* = c$, where $c$ is as in \eqref{c-value}. When $0<p_s < T_s$ and $p_b \ge T_s$, one can switch the two sides, and calculating the corresponding $c$.

\begin{figure}[ht]
    \centering
    \includegraphics[width=0.75\textwidth]{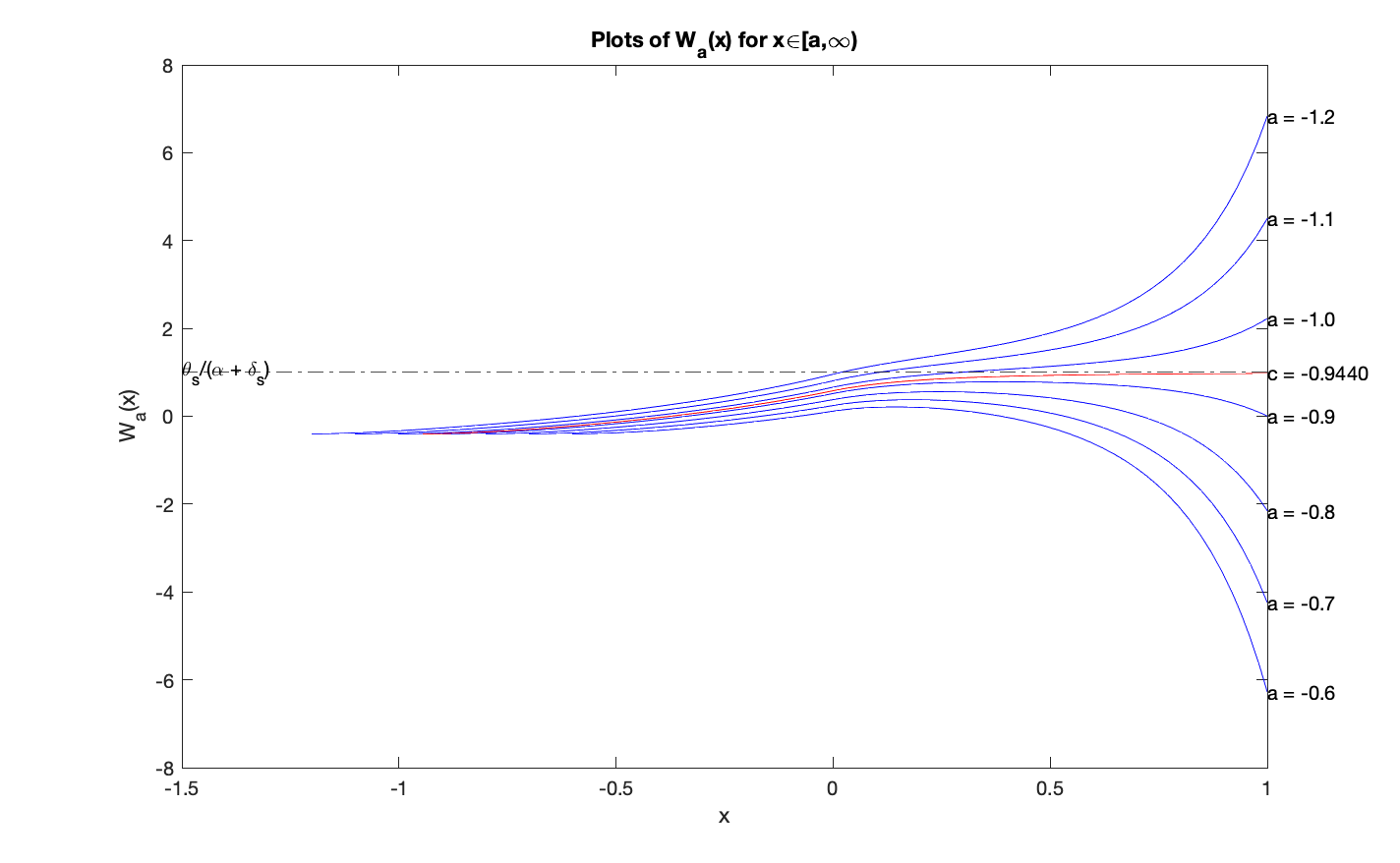}
    \caption{The functions $W_a(x), x\ge a,$ for different values of $a$.}
    \label{diff-a}
\end{figure}
\begin{example}\label{ex-1}
We set $\sigma^2 = 1, \beta =2, \alpha = 1, \delta_b = 2, \delta_s = 4, \theta_b = 4, \theta_s = 5$, and $p_b = 0.4, p_s = 0.1$. It is easily seen that $p_b < T_b = \theta_b/(\alpha+\delta_b) = 4/3$ and $p_s < T_s = \theta_s/(\alpha+\delta_s) = 1.$ We present two figures. In Figure \ref{diff-a}, the plots of $W_a(x), x\ge a$ for $a=-1.2, -1.1, \ldots, -0.7, -0.6$ are derived along with the value $c=-0.944$. One can observe that when $a<c$, $W_a(x)\to\infty$ as $x\to\infty$, while when $a\in (c, 0)$, $W_a(x)\to-\infty$ as $x\to\infty$. In Figure \ref{values}, we find the values $a^* = -0.5248$ and $b^* = 0.1104$. 
\begin{figure}[ht]
    \centering
    \includegraphics[width=0.75\textwidth]{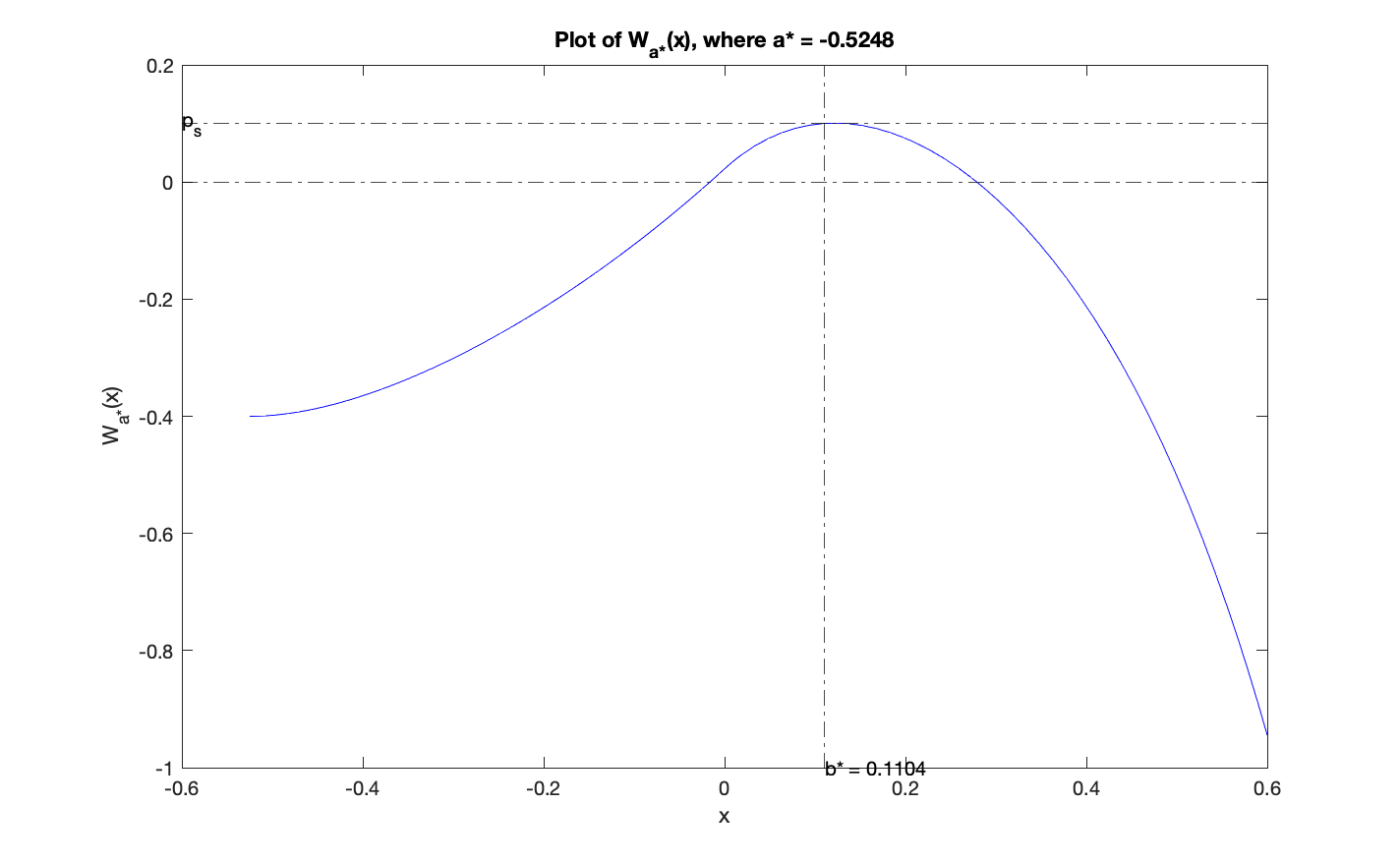}
    \caption{Finding the values $a^*$ and $b^*$.}
    \label{values}
\end{figure}
\end{example}

\begin{example}
We set $\sigma^2 = 1, \beta =2, \alpha = 1, \delta_b = 2, \delta_s = 4, \theta_b = 4$ and $\theta_s = 5$ as in Example \ref{ex-1}. We consider different values of $p_b$ and $p_s$ and observe the corresponding changes of the optimal buffer sizes. For convenience, we fix the value of $p_b$, and let the value of $p_s$ vary. In this setting, the value of $c$ does not change. We summarize the numerical results in the following table. We observe that as $p_s$ approaches $T_s$, both $-a^*$ and $b^*$ increase and in particular, $a^*$ goes to $c$.  
\begin{table}[h]\label{table-1}
\centering
\begin{tabular}{ |c|c|c|c|c|c|c| } 
 \hline
 $p_s$ & $p_b$ & $T_s$ & $T_b$ & $c$ & $a^*$ (buyer side) & $b^*$ (seller side)\\
 \hline
 $0.1$ & $0.4$ & $1$ & $4/3$ & $-0.9440$ & $-0.5248$ & $0.1104$ \\
 \hline
 $0.3$ & $0.4$ & $1$ & $4/3$ &$-0.9440$& $-0.6568$ & $0.1333$ \\
 \hline
 $0.5$ & $0.4$ & $1$ & $4/3$ &$-0.9440$& $-0.7707$ & $0.1935$ \\
 \hline
 $0.7$ & $0.4$ & $1$ & $4/3$ &$-0.9440$& $-0.8671$ & $0.2876$ \\
 \hline
 $0.9$ & $0.4$ & $1$ & $4/3$ &$-0.9440$& $-0.9345$ & $0.5501$ \\
 \hline
\end{tabular}
\caption{The value of $p_b$ is fixed to be $0.4$ and the value of $p_s$ varies from $0.1$ to $0.9$ approaching $T_s = 1$.}
\end{table}


\end{example}

\section{Asymptotic Optimality}\label{sec:as}
In this section, we establish the following two main results: In the first theorem (see Theorem \ref{thm5.4}), we prove  the value function  $V(x)$ of the DCP given in \eqref{value-fn}  is an asymptotic lower bound for the sequence of value functions $\{\hat{V}^n(x)\}$ of the QCPs given in \eqref{2.20}. In the second result, Theorem \ref{thm5.2} exhibits an asymptotically optimal sequence of controlled queueing processes where the associated sequence of cost functionals  converges to the value function  $V(x).$ Hence the lower bound  described in Theorem \ref{thm5.4} is achievable. 

{We simplify the form of the cost functional  $\hat{J}^n( \hat{X}^n(0), \hat{U}^n_s, \hat{U}^n_b)$ in \eqref{cost-1} using  Proposition \ref{prop-3.9} to match the form of the cost functional $J(X_x, U_s, U_b)$ in \eqref{cost-fn}. More precisely, from the proof of the following Theorem \ref{thm5.4}, the cost functional $\hat{J}^n( \hat{X}^n(0), \hat{U}^n_s, \hat{U}^n_b)$ can be written in the following form
\begin{equation}
\hat{J}^n( \hat{X}^n(0), \hat{U}^n_s, \hat{U}^n_b)= E \left(\int_{0}^{\infty} e^{- \alpha t} [ C(\hat{X}^n(t))  dt +p_s d\hat{U}^n_s(t) +p_b d\hat{U}^n_b(t)]\right) + \varepsilon_n,
\label{5.4-1}
\end{equation}
where $\lim_{n \rightarrow \infty} \varepsilon_n =0.$ We note that \eqref{5.4-1} is in agreement with the form of \eqref{cost-fn} for the DCP with a negligible term $\varepsilon_n$, and this plays an important role in establishing the asymptotic lower bound of the DCP for the QCP (Theorem \ref{thm5.4}). }

\begin{theorem}[Asymptotic lower bound]\label{thm5.4}
Let $(\hat{X}^n,   \hat{U}^n_s,  \hat{U}^n_b)$ be an admissible process for the $n^{\rm th}$ system. Then 
\begin{equation}
   \liminf\limits_{n \rightarrow \infty}\hat{J}^n( \hat{X}^n(0), \hat{U}^n_s, \hat{U}^n_b) \geq V(x), \quad  \text{ for all } x \in \rr, 
  \label{5.4}
\end{equation}
where $V(\cdot)$ is the value function of the DCP given in  \eqref{value-fn}.

\label{thm5.1}
\end{theorem}

\proof
{We first show \eqref{5.4-1}. In \eqref{cost-1}, using the Fubini's theorem,
\begin{equation}\label{qcp-sim-1}
\begin{aligned}
E\left(\int_{0}^{\infty} e^{- \alpha t}  d\hat{G}^n_s(t)\right) &= E\left( \int_{0}^{\infty} \int_{t}^{\infty} \alpha e^{- \alpha u} du  d\hat{G}^n_s(t) \right)\\
&= \alpha E \left(\int_{0}^{\infty} \hat{G}^n_s(u) e^{- \alpha u} du\right).
\end{aligned} 
\end{equation}
Using  \eqref{3.6} and  \eqref{3.22},  we have   $ E[ \sup_{t\in[0,T]}| \hat{G}^n_s(t)- \delta_s \int_{0}^{t} (\hat{X}^n(s))^+ ds |^2] \leq C(1+T^m) $ where the constant $C>0$ and the integer $m>1$ are independent of $n$ and $T.$  Hence using  \eqref{3.80} and  Fubini's theorem, we conclude that   
\begin{align}\label{qcp-sim-2}
\lim\limits_{n \rightarrow \infty} E\left( \int_{0}^{\infty} e^{- \alpha t} \left| \hat{G}^n_s(t)- \delta_s \int_{0}^{t} (\hat{X}^n(s))^+ ds \right| dt\right) =0.
\end{align}
By reversing the above procedure, we have   
\begin{align}\label{qcp-sim-3}
\alpha E \left(\int_{0}^{\infty} e^{- \alpha t}  \int_{0}^{t} (\hat{X}^n(s))^+ ds dt\right) =  E \left(\int_{0}^{\infty} e^{- \alpha t} (\hat{X}^n(t))^+ dt\right).
\end{align}
From \eqref{qcp-sim-1} -- \eqref{qcp-sim-3},   $\lim_{n \rightarrow \infty} |E \int_{0}^{\infty} e^{- \alpha t}  d\hat{G}^n_s(t)- \delta_s E \int_{0}^{\infty} e^{- \alpha t} (\hat{X}^n(t))^+ dt| = 0$, and similarly, $\lim_{n \rightarrow \infty} |E \int_{0}^{\infty} e^{- \alpha t}  d\hat{G}^n_b(t)- \delta_b E \int_{0}^{\infty} e^{- \alpha t} (\hat{X}^n(t))^- dt | = 0.$ 
Recall the following parameters from the DCP \eqref{cost-fn}, $\theta_s=c_s +r_s \delta_s > 0$ and $\theta_b=c_b +r_s \delta_b >0$,
and the effective running cost function $C(x)=  \theta_s x^+ + \theta_b x^-, x\in\rr,$ defined in \eqref{dcp-holding}. Hence, the cost functional can be written in the form of \eqref{5.4-1}.}
%
%
%

Now fix $x\in\rr.$ Let $(X_x, U_s, U_b)$  be an arbitrary admissible strategy in \eqref{ctrl-set}. We first establish $E[X_x(t)^2] \leq x^2 + K t $ for all $t \geq 0$, where $K>0$ is a constant independent of $t$  as well as the strategy  $(X_x, U_s, U_b).$ Using the It$\hat{o}$'s formula to $X_x(t)^2$  in \eqref{4.2} and the fact   $X_x(t) dU(t) \geq 0,$  we obtain 
\begin{equation}
 E[X_x(t)^2] \leq x^2 -2 E\left(\int_{0}^{t} [\delta_0 X_x(s)^2-\beta X_x(s)]ds \right),
 \label{5.6-1}
 \end{equation}
  where $ \delta_0=\min\{\delta_b, \delta_s\} >0. $  Notice that  $\delta_0 y^2-\beta y= \delta_0( y-c_0)^2-c_0^2 \geq -K$, where $ c_0={\beta}/{(2 \delta_0)} $ and $K= \delta_0 c_0^2= {\beta}/{(2 \delta_0)}.$  Using  this in the RHS of \eqref{5.6-1}, we obtain   $E[X_x(t)^2] \leq x^2 + K t $ for all $t \geq 0$.

Let $\epsilon >0$ be arbitrary and pick $T_0>0$ so that  $ \int_{T_0}^{\infty} e^{-\alpha t} (x^2 +1 +K t) dt  < \epsilon$. We now consider $T>T_0.$   Let   $(\hat{X}^n,   \hat{U}^n_s,  \hat{U}^n_b)$ be an admissible process which satisfies \eqref{3.92}.  By Theorem \ref{thm3.2}, we consider any limit point  $(X_x,   U_s, U_b)$  in $D^3[0,  T]$  which satisfies  \eqref{3.90}.  
Using \eqref{5.4-1}, integration by parts, and the Fatou's lemma, we obtain  
\begin{align*} 
& \liminf\limits_{n \rightarrow \infty} \hat{J}^n( \hat{X}^n(0), \hat{U}^n_s, \hat{U}^n_b)  = \liminf\limits_{n \rightarrow \infty} E \left(\int_{0}^{\infty} e^{- \alpha t} [ C(\hat{X}^n(t))  dt +p_s d\hat{U}^n_s(t) +p_b d\hat{U}^n_b(t)]\right) \\
& = \liminf\limits_{n \rightarrow \infty} \left[E \left(\int_{0}^{\infty} e^{- \alpha t} C(\hat{X}^n(t))  dt \right) +  \alpha  E \int_{0}^{\infty} e^{- \alpha t} [ p_s \hat{U}^n_s(t) +p_b \hat{U}^n_b(t)] dt\right]\\
& \ge E \left(\int_{0}^{T} e^{- \alpha t} C(X_x(t))  dt \right) +  \alpha  E \int_{0}^{T} e^{- \alpha t} [ p_s U_s(t) +p_b U_b(t)] dt \\
& =   E \left(\int_{0}^{T} e^{- \alpha t} [ C(X_x(t))  dt +p_s dU_s(t) +p_b dU_b(t)]\right).
\end{align*}  
Next, we extend $(X_x,   U_s, U_b)$ on $[0,T]$ to the  process  $(\tilde{X}_x,   \tilde{U}_s, \tilde{U}_b)$ on $[0, \infty)$ in $D^3[0,  \infty)$ by letting $ \tilde{U}_s(t) = U_s(T)$  and $ \tilde{U}_b(t) = U_b(T) $ for all  $t \geq T$   and   $(\tilde{X}_x,   \tilde{U}_s, \tilde{U}_b)$ satisfies  \eqref{3.90}.   Then $(\tilde X_x, \tilde U_s, \tilde U_b)$  is an admissible strategy for the DCP in \eqref{ctrl-set}  and   $J(\tilde X_x,  \tilde U_s, \tilde U_b)$  satisfies  
$J(\tilde X_x,  \tilde U_s, \tilde U_b) \leq E \left(\int_{0}^{T} e^{- \alpha t} [ C(X_x(t))  dt +p_s dU_s(t) +p_b dU_b(t)]\right)  + M \epsilon,$  
 where $M=\max \{ \theta_s, \theta_b \}.$  
 Hence $  \liminf_{n \rightarrow \infty} \hat{J}^n( \hat{X}^n(0), \hat{U}^n_s, \hat{U}^n_b)   +M  \epsilon >   J(\tilde X_x,  \tilde U_s, \tilde U_b) \geq V(x), $  for every $x$ and $\epsilon >0$. Letting $\epsilon\to0$, then \eqref{5.4} follows and  the proof  is complete. \Halmos

 {Our next aim is to  construct an  asymptotically optimal  sequence of processes $(\hat{X}^n,   \hat{U}^n_s,  \hat{U}^n_b)$ to achieve the lower bound in \eqref{5.4}.
  There will be  four different types of asymptotically optimal state process sequences  corresponding to  the  four types of optimal controls for the DCP constructed in Theorem \ref{thm:solution-DCP}.}

{  For each $n \geq 1,$  we consider a diffusion-scaled state process  $(\hat{X}^n,   \hat{U}^n_s,  \hat{U}^n_b)$  which satisfies \eqref{state-2} with time-independent  reflection barriers at $a^*_n <0< b^*_n$ and initial data $\hat{X}^n(0)=x_n.$  Here $ \hat{X}^n(t)$ takes values in a lattice of the form  $\{ {j}/{\sqrt{n}}: j \text{ is an integer} \}$ and hence, the initial position  $x_n,$ as well as the reflecting barriers $a^*_n$ and   $ b^*_n$ are assumed to take values in the same lattice as well.  We choose the sequences $\{a^*_n\}$ and $\{b^*_n\}$ in the lattice so that    $\lim_{n \rightarrow \infty} a^*_n =a^*$ and   $\lim_{n \rightarrow \infty} b^*_n =b^*,$ where  $a^* <0<b^*.$ In the following lemma, for the simplicity of the presentation, we assume $a^*_n\leq x_n \leq  b^*_n,$ but if $x_n$ is outside the interval $[a^*_n, b^*_n],$  then there  is  an initial jump at time $0-$ to the nearest point in $\{a^*_n,  b^*_n \}$  from $x_n.$ Since the jump size of this possible jump is $\min\{|x_n-a^*_n|, |x_n-b^*_n| \}$ and it is bounded, the conclusion of the lemma clearly holds for this general case.}

\begin{lemma} 
\label{local-times}
  Let $n \geq 1$ and $(\hat{X}^n,   \hat{U}^n_s,  \hat{U}^n_b)$  satisfy \eqref{state-2} with time-independent  reflection barriers at $a^*_n <0< b^*_n$ and initial data $\hat{X}^n(0)=x_n \in[a^*_n, b^*_n].$ 
We assume there is a $\delta>0$ so that 
$a^*_n <-\delta <0< \delta < b^*_n$ for all $n$, and for each $n,$ the process $\hat{X}^n(t) $ and    $a^*_n,  b^*_n$ take values in the lattice  $\{{j}/{\sqrt{n}}: j \text{ is an integer} \}$. Then the following moment bound holds:
\begin{equation}
 E[  \hat{U}^n_s(T)^2 +  \hat{U}^n_b(T)^2 ] \leq C(1+T^m),
\label{local-time-moment} 
\end{equation}
where the constant $C>0$ and the integer $m \geq 2$ are independent of $n$ and $T.$
\end{lemma}
\proof
{We notice $(\hat{X}^n,   \hat{U}^n_s,  \hat{U}^n_b)$ which satisfies \eqref{state-2} with reflection barriers at $a^*_n <0< b^*_n$   is the unique solution to the SP (see \cite{kruk})   on  $[a^*_n, b^*_n],$  with the input process $ \hat{\zeta}^n(t)  -\hat{G}^n(t) $  as given in \eqref{3.2} and \eqref{3.3}. To use the comparison theorem for the constraint processes, we consider for $t \geq 0,$ 
\[ \hat{Y}^n(t)= x_n+ \hat{A}^n_s(t) - \hat{A}^n_b(t)  - \tilde{U}^n_s(t) + \tilde{U}^n_b(t),\]
where $\tilde{U}^n_b(t)=n^{-1/2}\int_{0}^{t}1_{\{\hat{Y}^n(s)=a^*_n\}} dA^n_b(s),$  $\tilde{U}^n_s(t)=n^{-1/2}\int_{0}^{t}1_{\{\hat{Y}^n(u)=b^*_n\}}dA^n_s(u).$
Let 
\[\hat Z^n(t) = x_n + \hat A^n_s(t) - \hat A^n_b(t) - \hat G^n_s(t) - \check U^n_s(t) + \check U^n_b(t),\]
where
 $\check{U}^n_b(t)=n^{-1/2} \int_{0}^{t}1_{\{\hat{Z}^n(s)=a^*_n\}}dA^n_b(s),$ and   $\check{U}^n_s(t)=n^{-1/2}\int_{0}^{t}1_{\{\hat{Z}^n(u)=b^*_n\}}dA^n_s(u). $
Since both $\hat G^n_s(\cdot)$ and $\hat G^n_b(\cdot)$ are nondecreasing, using the comparison theorem (Theorem 1.7 of  \cite{kruk}) for $\hat Y^n$ and $\hat Z^n$, we have 
\begin{align}\label{eq1}
\tilde U^n_b \le \check U^n_b \le \tilde U^n_b + \hat G^n_s, \quad \check U^n_s \le \tilde U^n_s \le \check U^n_s + \hat G^n_s.
\end{align}
Next using the comparison theorem for $\hat X^n$ and $\hat Z^n$, we have 
\begin{align}\label{eq2}
\hat U^n_b \le \check U^n_b \le \hat U^n_b + \hat G^n_b, \quad \check U^n_s \le \hat U^n_s \le \check U^n_s + \hat G^n_b.
\end{align}
From \eqref{eq1} and \eqref{eq2}, we have 
\begin{align} \label{compare}
\hat U^n_b \le \tilde U^n_b + \hat G^n_s, \quad \hat U^n_s \le \tilde U^n_s + \hat G^n_b. 
\end{align}

We note that the processes $A^n_b$ and $A^n_s$ are independent. Hence  $\hat{A}^n_s(t) - \hat{A}^n_b(t)$ is a pure jump process with jump size $ {\pm1}/{\sqrt{n}}.$  
Next, we rely on the oscillation inequalities (see Proposition 4.1 part (ii) (c) in \cite{ward-kumar}  or \cite{dai-dai})  to obtain an upper bound for  $ \hat{U}^n_s(T).$ 
For any $f$ in $D[0, T],$ let $Osc(f, T)$ be defined by \eqref{Osc1}, then using  Proposition 4.1 part (c) of \cite{ward-kumar}
we obtain for some $\kappa >0$ independent of $n$ and $T$ so that 
\begin{equation}
\tilde{U}^n_s(T)\leq \kappa\left( Osc( x_n+\hat{A}^n_s - \hat{A}^n_b , T) + \frac{1}{\sqrt{n}}\right) \le  \kappa(2 \| \hat{A}^n_b\|_T +  2\| \hat{A}^n_s\|_T + 1).
 \label{Osc}
\end{equation}
A similar  estimate can be obtained for $ \tilde{U}^n_b(T).$  Next, using this estimate,  \eqref{compare}, \eqref{arrivals-M2} and Proposition \ref{prop-3.6} (i), we can easily obtain \eqref{local-time-moment}. This completes the proof.}
\Halmos

\begin{theorem}[Asymptotic optimality]\label{thm:AO}
Under Assumptions \ref{assump:initial} -- \ref{assump:admiss}, the following results hold. {Let $a^*, b^*, \tilde a^*$ and $\tilde b^*$ be as in Theorem \ref{thm:solution-DCP}.}
\begin{itemize}
\item[\rm (i)]When $p_s \geq T_s$ and  $p_b \geq T_b$, the sequence of processes
 $(\hat{X}^n,   \hat{U}^n_s,  \hat{U}^n_b)$  with the zero control policy $ \hat{U}^n_s(t) \equiv 0,$ and $ \hat{U}^n_b(t) \equiv 0 $  for all $t \geq 0$ is asymptotically optimal. 
  
\item[\rm (ii)] When $0 <  p_s< T_s$ and  $0 < p_b < T_b$, for $a^*_n < a^* < 0  < b^*<b^*_n$ satisfying $\lim_{n \rightarrow \infty} a^*_n =a^*$ and $\lim_{n \rightarrow \infty} b^*_n =b^*$, the associated sequence of the two-sided reflected processes   $(\hat{X}^n,   \hat{U}^n_s,  \hat{U}^n_b)$ with reflecting boundaries at $a^*_n$ and $b^*_n$  provides an asymptotically optimal sequence.

\item[\rm (iii)] When  $ p_s \geq  T_s$ and  $0 < p_b < T_b$, for $\tilde a^*_n  < \tilde a^*<0 $ satisfying $\lim_{n\to\infty} \tilde a^*_n = \tilde a^*$, the corresponding sequence of the one-sided reflected  processes $(\hat{X}^n,   \hat{U}^n_s,  \hat{U}^n_b)$ with reflecting boundary at $\tilde a^*_n$ yields an asymptotically optimal sequence.

\item[\rm (iv)] When  $ 0 < p_s < T_s$ and  $ p_b \geq  T_b, $ for $\tilde b^*_n > \tilde b^* > 0$ satisfying $\lim_{n\to\infty} \tilde b^*_n = \tilde b^*$, the corresponding sequence of the one-sided reflected  processes   $(\hat{X}^n,   \hat{U}^n_s,  \hat{U}^n_b)$ with reflecting boundary at $\tilde b^*_n$ yields an asymptotically optimal sequence.

\end{itemize}
\label{thm5.2}
\end{theorem}

\proof
Let $x$ be fixed. By Theorem \ref{thm:solution-DCP} (i),  zero control  strategy is optimal for the DCP when  $p_s \geq T_s$ and  $p_b \geq T_b. $ The corresponding state process $X_x$ is the unique strong solution to \eqref{4.26}.  Moreover, the value function    of the DCP  is given by  
 $V(x)= E \int_{0}^{\infty} e^{-\alpha s} C(X_x(s)) ds.$

We now consider the process $(\hat{X}^n,   \hat{U}^n_s,  \hat{U}^n_b)$ with the zero control policy  $ \hat{U}^n_s(t) \equiv 0,$ and $ \hat{U}^n_b(t) \equiv 0 $  for all $t \geq 0$. Then  by \eqref{3.92},
$\hat{X}^n(t)=  \hat{\zeta}^n(t)  +\epsilon_n(t) -  \int_{0}^{t}h(\hat{X}^n(s))ds,$
where $ \lim_{n \rightarrow \infty} \| \epsilon_n\|_T = 0$ in probability for any $T>0.$  Using  Theorem 11.4.5 of \cite{whitt}, 
 \eqref{3.4}, the Lipschitz continuity of  $h(\cdot)$ and  following  Theorem 2.11 of  \cite{talreja}, we conclude that $\hat{X}^n(\cdot)$ converges weakly to $X_x(\cdot)$ in $D[0, T].$ Let $\epsilon >0$  be arbitrary. Using 
\eqref{3.6} and the  Lipschitz continuity of $C(\cdot)$, we can find a large $T>0$ so that 
\[\sup_{n \geq 1} E \int_{T}^{\infty} e^{- \alpha t}  C(\hat{X}^n(t))  dt  < \epsilon.\] 
Since $\hat{X}^n$ converges weakly to $X_x$ in $D[0, T],$  using  \eqref{3.6}, we obtain 
\[\lim\limits_{n \rightarrow \infty}  E \int_{0}^{T} e^{- \alpha t}  C(\hat{X}^n(t))  dt = E \int_{0}^{T} e^{-\alpha t} C(X_x(t)) dt,\]
and therefore,   
 \[\liminf_{n \rightarrow \infty}  E \int_{0}^{\infty} e^{- \alpha t}  C(\hat{X}^n(t))  dt \le  E \int_{0}^{T} e^{-\alpha t} C(X_x(t)) dt + \epsilon \leq  V(x)+ \epsilon.\]
 Since $\epsilon >0$ is arbitrary,   $\liminf_{n \rightarrow \infty}  E \int_{0}^{\infty} e^{- \alpha t}  C(\hat{X}^n(t))  dt \leq V(x).$ Using \eqref{5.4}, we obtain  $\liminf_{n \rightarrow \infty}  E \int_{0}^{\infty} e^{- \alpha t}  C(\hat{X}^n(t))  dt = V(x).$  Hence part (i) follows.

To prove part (ii), consider $(X^*_x, L_{a^*},  L_{b^*}) $ satisfying  \eqref{reflected} with reflection barriers at $a^*$ and $b^*.$   It is an optimal control for the DCP in this parameter regime  as proved in Theorem \ref{thm:solution-DCP} (ii).  Our aim is to construct an asymptotically optimal sequence  $(\hat{X}^n,   \hat{U}^n_s,  \hat{U}^n_b)$  which converges weakly to  $(X^*_x, L_{a^*},  L_{b^*}) .$  { For each $n \geq 1,$  we consider a  state process  $(\hat{X}^n,   \hat{U}^n_s,  \hat{U}^n_b)$  which satisfies \eqref{state-2} with time-independent  reflection barriers at $a^*_n <0< b^*_n$ and initial data $\hat{X}^n(0)=x_n$  as described in Lemma \ref{local-times}.}  { Here we  choose  $\{a^*_n\}$ and $\{b^*_n\}$ so   that   $a^*_n <a^*<0<b^*< b^*_n,$   $\lim_{n \rightarrow \infty} a^*_n =a^*$ and   $\lim_{n \rightarrow \infty} b^*_n =b^*.$}  

 We assume  \eqref{initial} and consider the case  $x\in[a^*, b^*].$ Then the state equation can be written in the form  \eqref{3.92} and  $\hat{\zeta}^n$ converges weakly as in \eqref{3.4}.  Using   the Skorokhod representation theorem, we  assume that  all these processes are defined in a same probability space and  $\lim_{n \rightarrow \infty} \sup_{t\in[0,T]}| \hat{\zeta}^n(t) - (x+ \sigma B(t) +\beta t)|=0 \text{ a.s.}$ The reflected diffusion  process  $(X^*_x, L_{a^*},  L_{b^*}) $ is a strong solution to \eqref{reflected} in this probability space with respect to the same Brownian motion $B.$  Next, we use Proposition \ref{propB.4} in Appendix \ref{sec:sm} to conclude  $\lim_{n \rightarrow \infty} \| \hat{X}^n - X^*_x\|_T =0,$   $\lim_{n \rightarrow \infty} \| \hat{U}^n_s - L_b^* \|_T =0 $ and  $\lim_{n \rightarrow \infty} \| \hat{U}^n_b - L_a^* \|_T =0 $  for any $T>0.$ Next let $\epsilon >0$ be arbitrary. Using the moment bounds in \eqref{3.6} and \eqref{local-time-moment},  we can find  a $T>0$ so that  
\[\sup_{n \geq 1}E \int_{T}^{\infty} e^{- \alpha t} [ C(\hat{X}^n(t))  dt +p_s d\hat{U}^n_s(t) +p_b d\hat{U}^n_b(t)]  <\epsilon\] and similarly, 
\[ E \int_{T}^{\infty} e^{- \alpha t} [ C(X^*_x(t))  dt +p_s dL_{a^*}(t) +p_b dL_{b^*}(t)]  <\epsilon.\]
 Using the convergence of $(\hat X^n, \hat U^n_s, \hat U^n_b)$ to $(X^*_x, L_b^*, L_a^*)$ almost surely on $[0, T],$  and by an argument similar to part (i), we can conclude that $ | \hat{J}^n( \hat{X}^n(0), \hat{U}^n_s, \hat{U}^n_b) -   J(X^*_x,  L_{a^*}, L_{b^*}))| < 2 \epsilon $ and consequently,  $ \lim_{n \rightarrow \infty} \hat{J}^n( \hat{X}^n(0), \hat{U}^n_s, \hat{U}^n_b) =  J(X^*_x,  L_{a^*}, L_{b^*}).$  This complets the proof of part (ii).

Proofs of parts (iii) and (iv) uses only the one-sided     Skorokhod map  and hence the proofs are  much simpler than part (ii).  In each case, the proof is very similar to that of Theorem 4.2 in \cite{weera4}. Therefore, we omit it here.
\Halmos

\begin{appendix}

\section{Two-sided Skorokhod maps} \label{sec:sm}

In this sub-section, we establish several  results to supplement the work of   \cite{kruk, reed} and  Chapter 14 of \cite{whitt}.  These results enable us to obtain asymptotically optimal strategies in Theorem  \ref{thm5.2}. We first revisit the definition of the two-sided Skorokhod map (see Definition 1.2 in \cite{kruk}). We follow the notation in \cite{kruk}, and  let $ \Gamma_{a, b}: D[0, \infty) \rightarrow D[0, \infty)$ represent the two-sided Skorokhod map on the (time-independent) interval  $[a, b]$.

\begin{definition}\label{2-sp_constant}
Let the constants $a < b$. Given $\psi \in D[0,\infty)$, there exists a unique pair of functions $\phi \in D[0,\infty)$ and $\eta$ of bounded variation such that 
\begin{itemize}
\item[\rm (i)] for each $t\ge 0$, $\phi(t) = \psi(t) + \eta(t) \in [a, b]$;
\item[\rm (ii)] $\eta(0-)=0$, $\eta(0)\ge 0$, and $\eta$ has the decomposition $\eta = \eta^l - \eta^r$ satisfying that $\eta^l$ and $\eta^r$ are non-decreasing, and 
\[
\int_0^\infty 1_{\{\phi(s) > a\}} d\eta^l(s) = 0 \quad \mbox{and} \quad \int_0^\infty 1_{\{\phi(s) < b\}} d\eta^r(s) = 0. 
\]
\end{itemize}
\end{definition}

The map $\Gamma_{a,b}: D[0,\infty) \to D[0,\infty)$ that takes $\psi$ to the corresponding $\phi$ is referred to as the two-sided Skorokhod map on $[a, b]$, and the triple $(\phi, \eta^l, \eta^r)$ is referred to as the Skorokhod decomposition of $\psi$ on $[a,b].$ From the comparison properties of the Skorokhod map on $[a, b]$ (Theorem 1.7 of \cite{kruk}), the Skorokhod decomposition is unique. 

 
 \begin{lemma}
 Let $c<a<b<d.$ Then  for any $f$ in  $D[0, \infty),$
\begin{equation}
\| \Gamma_{c, d}(f) - \Gamma_{a, b}(f) \|_T \leq  3[ |a-c| + |b-d| ].
\label{7.10}
\end{equation} 
 \label{lmB.1}
 \end{lemma}
 
 \proof
Using equation (1.14) in \cite{kruk} and the notation  therein, for any $f$ in  $D[0, \infty)$,  
\begin{align}\label{eq:twosidesm}
\Gamma_{a, b}(f)= \Lambda_{a, b} \circ \Gamma_a(f),
\end{align} 
where for $f\in D[0, \infty)$ and $t\ge 0$,
$\Gamma_a(f)(t) = f(t) + \sup_{s\in [0, t]} [a - f(s)]^+,$
and for $g\in D[a, \infty)$,
$\Lambda_{a,b}(g)(t) = g(t) - \sup_{s\in[0,t]}\left((g(s) - b)^+ \wedge \inf_{u\in[s,t]} (g(u) -a) \right).$
 We notice that 
\begin{align*} 
\| \Gamma_a(f) - \Gamma_c(f) \|_T  = \sup_{t\in [0, T]} \left|\sup_{s\in [0, t]} [a - f(s)]^+ - \sup_{\tilde s\in [0, t]} [c - f(\tilde s)]^+\right| \leq  |a-c|,
\end{align*}
and 
\begin{align*}
& \| \Lambda_{a,b}(f)- \Lambda_{c, d}(f)\|_T  = \sup_{t\in [0, T]} |\Lambda_{a,b}(f)(t) - \Lambda_{c,d}(f)(t)|  \\
& = \sup_{t\in [0, T]} \left|\sup_{s\in[0,t]}\left((g(s) - b)^+ \wedge \inf_{u\in[s,t]} (g(u) -a) \right) - \sup_{\tilde s\in[0,t]}\left((g(\tilde s) - d)^+ \wedge \inf_{u\in[\tilde s,t]} (g(u) -c) \right)\right| \\
& \le |a-c| +|b-d|.
\end{align*}
From (1.16) of \cite{kruk}, we have $ \| \Lambda_{c, d}(f) - \| \Lambda_{c, d}(g)\|_T \leq 2 \|f-g\|_T $ for any $f$ and $g$ in   $D[0, \infty).$ Hence
\begin{equation}
\begin{aligned}
  \| \Gamma_{c, d}(f) - \Gamma_{a, b}(f) \|_T&  \leq   \| \Lambda_{c, d}(\Gamma_c(f)) -  \Lambda_{c, d}(\Gamma_a(f))\|_T  \\
  &\quad +  \| \Lambda_{c, d}(\Gamma_a(f)) -  \Lambda_{a,b}(\Gamma_a(f))\|_T \\
 &\leq 2 \| \Gamma_c(f) -\Gamma_a(f)\|_T + |a-c| +|b-d| \leq  3( |a-c| +|b-d| ).
 \end{aligned}
 \label{7.12}
 \end{equation} 
  \Halmos
 
In the next result, we consider  two convergent sequences  $\{a_n\}_{n\in\NN}$ and $\{b_n\}_{n\in\NN}$ so that 
$a_n<a<b < b_n,$  $a_n$ is increasing to $a$, and $b_n$ is decreasing to $b$ as $n\to\infty.$ 
  Let $\{Y_n\}_{n\in\NN}$ be a convergent sequence in $D[0, \infty)$ so that  $  \lim_{n \rightarrow \infty} \| Y_n-Y_\infty \|_T=0$  for each $T>0$, where  $Y_\infty$ is a continuous function. We introduce     
  \begin{equation}
  W_n(t)=Y_n(t) + \int_{0}^{t} h(\Gamma_{a_n, b_n}(W_n)(s)) ds,
   \label{7.22}
     \end{equation}
and
  \begin{equation}
  W_\infty(t)=Y_\infty(t) + \int_{0}^{t} h(\Gamma_{a, b}(W_\infty)(s)) ds,
   \label{7.24}
     \end{equation}  
     where $h$ is a Lipschitz continuous function. 
     

Consider the Skorokhod decomposition $(Z_n, \eta^l_n, \eta^u_n)$ of the function $W_n$ on the interval $[a_n, b_n]$, and the  Skorokhod decomposition $(Z_\infty, \eta^l_\infty, \eta^u_\infty)$ of the function $W_\infty$ on $[a, b].$  
 Next, we obtain the convergence results of the Skorokhod decompositions of $W_n$ and $W_{\infty}.$ This proof is closely related to that of Proposition 4.2 of \cite{reed}.
 \begin{proposition}
 There exists a constant $C_T>0$ which depends only on $T$ and the Lipschitz constant of  the function  $h$ such that 
    \begin{equation}
   \| Z_n-Z_\infty \|_T \leq C_T\left[  |a_n -a| + |b_n -b| + \| Y_n-Y_\infty \|_T \right].
     \label{7.26}
     \end{equation}
    Consequently,   $\lim_{n \rightarrow \infty} \| Z_n-Z_\infty \|_T=0$ and 
  \begin{equation}
   \lim\limits_{n \rightarrow \infty} \| \eta^l_n-\eta^l_\infty \|_T+ \| \eta^r_n-\eta^r_\infty \|_T=0.
   \label{7.28}
   \end{equation}
%
 \label{propB.4}
 \end{proposition}
 
 \proof
The existence and uniqueness of solutions to  \eqref{7.22} and \eqref{7.24} depend only on the Lipschitz continuity of the function $h$  as shown in Lemma 4.2  of \cite{reed}. Using the  Lipschitz continuity of $h$ in  \eqref{7.22} and \eqref{7.24}, we obtain 
$|W_n(t) - W_\infty(t)| \leq |Y_n(t) -Y_\infty(t)|  +C \int_{0}^{t} |\Gamma_{a_n, b_n}(W_n)(s)- \Gamma_{a, b}(W_\infty)(s)| ds,$  
where $C>0$ is a constant depending on the Lipschitz constant of  the function  $h$. Then, 
\begin{equation}
\begin{aligned}
|\Gamma_{a_n, b_n}(W_n)(s)- \Gamma_{a, b}(W_\infty)(s)| \leq& |\Gamma_{a_n, b_n}(W_n)(s)- \Gamma_{a, b}(W_n)(s)|\\ 
&+  |\Gamma_{a, b}(W_n)(s)- \Gamma_{a, b}(W_\infty)(s)|.
\label{7.30}
\end{aligned}
\end{equation}
Now using \eqref{7.12} and the Lipschitz property of  $\Gamma_{a,b}$,  for each $0 \leq t \leq T$ we obtain, 

 \begin{equation}
\|W_n - W_\infty\|_t \leq \|Y_n -Y_\infty\|_T  + 3T( |a_n-a| + |b_n-b|)+2C \int_{0}^{t} \|W_n- W_\infty\|_s ds. 
\label{7.32}
\end{equation} 
 Employing the Gronwall's inequality leads to 
 \begin{equation}
  \|W_n - W_\infty\|_T \leq (\|Y_n -Y_\infty\|_T  + 3T( |a_n-a| + |b_n-b|)e^{2CT},  
  \label{7.34}
  \end{equation}
 and hence  $    \lim_{n \rightarrow \infty}  \|W_n - W_\infty\|_T=0.$
Next we note that  
$\|Z_n - Z_\infty\|_T \leq  \|\Gamma_{a_n, b_n}(W_n)- \Gamma_{a, b}(W_n)\|_T 
+  \|\Gamma_{a, b}(W_n)- \Gamma_{a, b}(W_\infty)\|_T.$ 
  Using  \eqref{7.10} and the Lipschitz property of  $\Gamma_{a,b}$  we obtain 
$  \|Z_n - Z_\infty\|_T \leq  C[ |a_n -a| + |b_n -b| + \|W_n - W_\infty\|_T ]. $
Consequently, using \eqref{7.34}, \eqref{7.26} follows and thus the proof of part (i) is complete.

 To prove part (ii), we let  $K_n(t)= \eta^l_n(t)-\eta^r_n(t)$ and  $K_\infty(t)= \eta^l_\infty(t)-\eta^r_\infty(t) $ for all $t \geq 0$. By part (i), $  \lim_{n \rightarrow \infty} \| K_n-K_\infty \|_T=0.$  Next we follow the proof of Theorem 14.8.1 in \cite{whitt}. Let $\epsilon >0$ be  arbitrarily small satisfying  $0 < \epsilon < {(b-a)}/{100}.$ Let $t_1=\min\{T, \tau_1 \} $ where  $\tau_1 = \inf\{ t \geq 0 : Z_\infty(t)<a+\epsilon \}.$ Let  $t_2=\min\{T, \tau_2 \} $ where  $\tau_2 = \inf\{ t \geq t_1 : Z_\infty(t)>b-\epsilon\}.$ Inductively, let $t_{2j-1}=\min\{T, \tau_{2j-1} \} $ where  $\tau_{2j-1} = \inf\{ t \geq t_{2j-2}: Z_\infty(t)<a+\epsilon \}$ and  $t_{2j}=\min\{T, \tau_{2j} \} $ where  $\tau_{2j} = \inf\{ t \geq t_{2j-1} : Z_\infty(t)>b-\epsilon\}.$  Since $Y_\infty$ is a continuous function, $Z_\infty $ also continuous on  $[0, \infty).$  Consequently, there are only finitely many points $ 0 \leq t_1 <t_2 <\cdot \cdot  < t_m \leq T.$ Let $n_0 >1$ so that  $ \| Z_n-Z_\infty \|_T <{\epsilon}/{10}$ for any $n\ge n_0$. For $n \geq n_0,$ $l_\infty$ as well as $l_n$ increases only on a finite number of intervals of the form  $[ t_{2j-1},  t_{2j}]$ say $j=1, 2, ...r.$  On those intervals, $u_n$ and $u_\infty$ remain constant.  Moreover,  $\eta^l_n$ and $\eta^l_\infty$ remain constant on the intervals $( t_{2j},  t_{2j+1}).$  On each interval  $[ t_{2j-1},  t_{2j}]$, $j=1, 2, ...r,$  only the one sided Skorokhod map (see equation (1.4) in \cite{kruk}) will be applied. Hence, it is evident that 
 $  \lim_{n \rightarrow \infty} \| \eta^l_n-\eta^l_\infty \|_{[t_{2j-1},  t_{2j}]}=0$ on each interval.  Consequently,   $  \lim_{n \rightarrow \infty} \| \eta^l_n-\eta^l_\infty \|_T=0$ follows. This completes the proof.  
\Halmos

The following results are immediate from the above proposition.
 \begin{corollary}
 
 The Skorokhod decomposition  $(Z_n, \eta^l_n, \eta^r_n)$   converges to   $(Z_\infty, \eta^l_\infty, \eta^r_\infty)$ in $D^3[0, \infty)$ in Skorokhod $J_1$-topology.
 
 \label{cor-B.5}
 \end{corollary}

%
%
  
 \begin{corollary}
 Let $Y$ be in $D[0, \infty)$ be fixed.  Consider the Skorokhod decompositions   $(Z_n, \eta^l_n, \eta^r_n)$ and   $(Z_\infty, \eta^l_\infty, \eta^r_\infty)$ of $Y$ on $[a_n, b_n]$ and $[a, b]$, respectively. Assume that $ \lim_{n \rightarrow \infty} a_n=a $ and    $ \lim_{n \rightarrow \infty} b_n=b. $   Then 
    \begin{equation}
     \lim\limits_{n \rightarrow \infty} \| Z_n-Z_\infty \|_T=0,  
     \label{7.18}
     \end{equation}
     and
       \begin{equation}
   \lim\limits_{n \rightarrow \infty} \| \eta^l_n-\eta^l_\infty \|_T+ \| \eta^r_n-\eta^r_\infty \|_T=0.
   \label{7.20}
   \end{equation}
 \label{lmB.2}
 \end{corollary}

 At last we consider the two-sided SP on time dependent intervals. The following definition is from \cite{burdzy}. 
 
 \begin{definition}
Let $l, r\in D[0, \infty)$, where $l$ could take $-\infty$ and $r$ could take $\infty$. Given $\psi \in D[0,\infty)$ so that $l(t) < r(t)$ for $t\ge 0$, a pair of functions $(\phi, \eta) \in D^2[0, \infty)$ is said to solve the SP for $\psi$ on the time-dependent interval  $[l(\cdot), r(\cdot)],$  if and only if it satisfies the following properties. 
\begin{itemize}
\item[\rm (i)] for each $t\ge 0$, $\phi(t) = \psi(t) + \eta(t) \in [l(t), r(t)]$;
\item[\rm (ii)] $\eta$ has the decomposition $\eta = \eta^l - \eta^r$, where $\eta^l$ and $\eta^r$ are non-decreasing functions, and 
\[
\int_0^\infty 1_{\{\phi(s) > l(s)\}} d\eta^l(s) = 0 \quad \mbox{and} \quad \int_0^\infty 1_{\{\phi(s) < r(s)\}} d\eta^r(s) = 0. 
\]
\end{itemize}
\end{definition}
If there is a unique solution $(\phi, \eta)$ to the SP for $\psi$ on $[l, r]$, then we write $\phi = \Gamma_{l, r}(\psi)$, where $\Gamma_{l,r}$ will be referred to as the two-sided Skorokhod map on $[l, r]$, and the triple $(\phi, \eta^l, \eta^r)$ is referred to as the Skorokhod decomposition of $\psi$ on $[l,r].$ 

From Theorem 2.5 and Corollary 2.4 of \cite{burdzy}, if       $ \inf_{t\ge 0} (r(t)-l(t))>0$, {\textcolor{red}{then}} there is a unique pair $(\phi, \eta)$ that solves the SP for $\psi$, and from the composition properties (Section 3 in \cite{burdzy}), the corresponding Skorokhod decomposition $(\phi, \eta^l, \eta^r)$ is unique. 

Next we develop the  oscillation inequalities for the constrained processes of Skorokhod decomposition.  For a function $f$ in $D[0, T],$ we recall  its oscillation 
$ Osc(f, [t_1, t_2])$ in an interval $[t_1, t_2] \subset [0, T]$ is defined by \eqref{Osc1}, and the modulus of continuity $\omega(f, \delta, T)$ is given in \eqref{omega}. 

\begin{proposition} \label{prop-3.7}

Let the functions $l$ and $r$ be in  $D[0, \infty)$ and $\inf_{t\ge 0}(r(t)-l(t))>0 $ for all $t \geq 0.$  Given a function $\psi$ in  $D[0, \infty),$  let $(\phi, \eta)$ be the unique solution to the SP  in   $D^2[0, \infty)$ satisfying $l(t) \leq \phi(t) \leq r(t)$ for all $t \geq 0$ and  $\eta$ be a function of bounded variation on $[0, \infty)$ as described in Theorem 2.6 of \cite{burdzy}. Then  
for all $\delta >0$ and $T>0$,
\begin{equation}
\omega(\phi, \delta, T) \leq 4[ \omega(\psi, \delta, T) + \omega(l, \delta, T) + \omega(r, \delta, T)].
\label{3.60}
\end{equation}
\end{proposition}

\proof
Using Theorem 2.6 and Corollary 2.4 of \cite{burdzy}, we can write  $\phi(t)=\psi(t)-\Theta(\psi)(t)$ for  $t\geq 0,$ where 
\begin{equation}
\Theta(\psi)(t)=\max\{ b(t), h(t) \},
\label{3.62}
\end{equation}
where 
\begin{equation}
b(t)=\min \{ (\psi(0)-r(0))^+, \inf_{u\in [0, t]}(\psi(u)-l(u)) \},
\label{3.64}
\end{equation}
and 
\begin{equation}
h(t)=\sup_{s\in[0, t]} [ \min \{(\psi(s)-r(s)), \inf_{u\in [s, t]}(\psi(u)-r(u))\} ],
\label{3.66}
\end{equation}
for all $t \geq 0.$

 Since $\phi=\psi-\Theta(\psi),$   it is evident that
\begin{equation}
 \omega(\phi, \delta, T) \leq 4[ \omega(\psi, \delta, T) + \omega(\Theta(\psi), \delta, T)].
 \label{3.68}
 \end{equation}
 From  \eqref{3.62}, it follows that 
\begin{equation}
\omega(\Theta(\psi), \delta, T) \leq  \omega(b, \delta, T) + \omega(h, \delta, T).
\label{3.70}
\end{equation}
  Next, we estimate  $\omega(b, \delta, T)$  and   $ \omega(h, \delta, T)$  carefully.
Using \eqref{3.64}, we obtain  $|b(t)-b(s) \leq \sup\limits_{u\in [s, t]} | \psi(u) - \psi(s)| +  \sup\limits_{u\in [s, t]} |l(u) - l(s)|$ whenever $ 0 \leq s \leq t.$ Therefore,  $ \omega(b, \delta, T) \leq  \omega(\psi, \delta, T) +  \omega(l, \delta, T)$. Similarly, for  $ \omega(h, \delta, T),$  we can use \eqref{3.66} together with the simple inequality $| \min \{a, b\} -\min \{c, d\}| \leq |a-c| + |c-d| $ to obtain 
\begin{equation}
 \omega(h, \delta, T) \leq [ 2 \omega(\psi, \delta, T) + \omega(r, \delta, T)+ 
 \omega(l, \delta, T) ]. 
 \label{3.72}
 \end{equation}
Now combining  \eqref{3.68} through \eqref{3.72}, we obtain  \eqref{3.60}. This completes the proof.
\Halmos
\begin{remark}
Following the above proof, one can obtain the inequality  
$$  \sup_{u\in [s, t]} | \phi(u) - \phi(v)| \leq 4\left[ \sup_{u\in [s, t]} | \psi(u) - \psi(v)| + \sup_{u\in [s, t]} | l(u) - l(v)| +  \sup_{u\in [s, t]} | r(u) - r(v)| \right] 
$$ 
on any given interval $[s, t].$
\end{remark}



\begin{proposition}
Let  $(\phi, \eta^l, \eta^r)$ in $D^3[0, T]$ be the solution to the SP with  the input function $\psi$ in $D[0, T]$ and the time-dependent barriers $l$ and $r$ which are also  in $D[0, T].$  Assume  $\inf_{t\in [0, T]}[r(t)-l(t)] >0,$  and $l(0)\leq \psi(0) \leq r(0).$
Then  the following oscillation inequalities hold for any $0 \leq t_1 \leq t_2 \leq T$:
\begin{equation}
 Osc(\eta^l, [t_1, t_2]) \leq  C[ Osc(\psi, [t_1, t_2]) + Osc(l, [t_1, t_2])  ], 
 \label{7.281}
 \end{equation}
and
\begin{equation}
 Osc(\eta^r, [t_1, t_2]) \leq  C[ Osc(\psi, [t_1, t_2]) + Osc(r, [t_1, t_2])  ],
 \label{7.282}
 \end{equation}
where  $C>0$ is a generic constant independent of  $\psi, l$ and $r$ and  $[t_1, t_2].$

\label{lmB.3}
\end{proposition}

\proof

We focus on proving \eqref{7.281} the oscillation property for $\eta^l$ below and the proof of \eqref{7.282} is similar. The proof is divided into the following three steps. 

{\em Step 1.}  Oscillation inequality in $[0, b],$ where $b>0$ is a constant:\\
Using  Theorem 4.2 of \cite{dai-dai}, 
$ Osc(\eta^l, [t_1, t_2]) \leq \kappa (   Osc(\psi, [t_1, t_2]) + \| \Delta\eta^l \|_{[t_1, t_2]}).$  The constant $\kappa$ is independent  of the functions $\psi, l, r$ and $\phi,$ $t_1, t_2$ as well as the domain $[0, b].$
Since set of  piecewise constant functions are dense  in $D[0, T]$ with respect to uniform norm (see \cite{whitt}, Theorem 12.2.2.), we assume $\psi$ is piecewise constant. Then $(\phi, \eta^l, \eta^r)$ also piecewise constant  with the same points of discontinuity. We can use the construction of   $(\phi, \eta^l, \eta^r)$  in the equation (8.2) of Chapter 14  of  \cite{whitt}.  Let $ 0 \leq s_1 <s_2<...<s_m \leq T $ be the  points of discontinuity of $\psi.$ Then, $ \eta^l(s_i)-\eta^l(s_i-) \leq |\psi(s_i) -\psi(s_i-)|$ holds at each jump point $s_i$. Hence we obtain $ Osc(\eta^l, [t_1, t_2]) \leq 2 \kappa (   Osc(\psi, [t_1, t_2]))$  when $\psi$ is a piecewise constant function.  This inequality can be generalized for any $\psi$ in $D[0, T]$ using Proposition \ref{propB.4} and  \eqref{7.26}  by taking the function $h$ to be identically zero. Hence $ Osc(\eta^l, [t_1, t_2]) \leq 2 \kappa (   Osc(\psi, [t_1, t_2]))$ follows.

{\em Step 2.}  Oscillation inequality in $[0, \beta(t)],$ where $\beta(\cdot)>0$ is  in $D[0, T]:$\\
 We assume  $\inf_{t\in[0, T]} \beta(t) >0.$ Consider the triple $(\phi, \eta^l, \eta^r)$ corresponding to  the input function  $\psi.$ We pick $t_1 < t_2$ in $[0, T].$  We pick  $b>0$ so that  $0 < b \leq \inf_{t\in[0, T]} \beta(t).$  For convenience assume $t_1$ is not a point of discontinuity of  $\beta.$ Now consider the interval $[t_1, T].$  If  $ \phi(t_1)>b$ choose $c_0=\phi(t_1)-b$  and $\tilde\psi= \psi -c_0.$ If  $ \phi(t_1)\leq b$ choose $\tilde\psi=\psi.$ 
Now $(\phi, \eta^l(t)-\eta^l(t_1), \eta^r(t)-\eta^r(t_1))$ is the solution to the SP for $\psi$ on $[t_1, T]$ and  let  $ (\tilde\phi, \tilde\eta^l, \tilde\eta^r)$ be the solution to the SP for $\tilde\psi$ on $[t_1, T].$  Then by the comparison theorem in Proposition 3.5 of [6], (with $\nu$ identically zero), we have 
$\eta^l(t_2)-\eta^l(t_1) \leq \tilde\eta^l(t_2)-\tilde\eta^l(t_1).$ Next we compare the solution  $(\tilde\phi, \tilde\eta^l, \tilde\eta^r)$ of $\tilde\psi$  on  the time dependent interval $[0, \beta]$ with the solution  $(\check{\phi}, \check\eta^l, \check\eta^r)$ of $\tilde\psi$  on   $[0, b].$  We use the comparison theorem in Proposition 3.3 of  \cite{burdzy} to conclude $ \tilde\eta^l(t_2)-\tilde\eta^l(t_1)\leq \check{\eta}^l(t_2)-\check{\eta}^l(t_1). $  Using Step 1, we obtain  $ \check{\eta}^l(t_2)-\check{\eta}^l(t_1)\leq 2 \kappa (   Osc(\psi, [t_1, t_2]).$ Consequently, 
$  \eta^l(t_2)-\eta^l(t_1) \leq 2 \kappa (   Osc(\psi, [t_1, t_2]).$ 
 Hence $ Osc(\eta^l, [t_1, t_2]) \leq 2 \kappa (   Osc(\psi, [t_1, t_2])$  holds in $[0, \beta(\cdot)],$ when   $\beta(\cdot)$ is  in $D[0, T]$ and $\inf_{t\in[0, T]} \beta(t) >0.$ 

{\em Step 3.} Oscillation inequality in  $[l(t),  r(t)]$ with $\inf_{t\in[0, T]}[r(t)-l(t)] >0:$ \\
We can use the proof of Lemma 2.2 and equation (2.26) in \cite{slaby} to observe that the corresponding constraining processes are unchanged  under translation by the function $\alpha$ there. Let $\beta(t)= r(t)-l(t)$ for all $t$ in $[0, T].$  Hence  using Step 2, we obtain  $ Osc(\eta^l, [t_1, t_2]) \leq 2 \kappa (   Osc(\psi-l,  [t_1, t_2])).$ But  $ Osc(\psi-l,  [t_1, t_2])\leq  Osc(\psi, [t_1, t_2])+  Osc(l, [t_1, t_2]).$  Hence we have $ Osc(\eta^l, [t_1, t_2]) \leq 2\kappa (   Osc(\psi, [t_1, t_2])+  Osc(l, [t_1, t_2]))$ and this yields \eqref{7.281} and the proof is complete.
\Halmos

For a given $\delta >0,$  by choosing intervals $[t_1,  t_2]$ of length less than  $\delta,$ and using \eqref{7.281} and \eqref{7.282},   we can obtain the following corollary.

\begin{corollary}\label{cor:3}
Let $(\phi, \eta^l, \eta^r)$ in $D^3[0, T]$ be the solution to SP for $\psi$ in $D[0, T]$ and the time-dependent barriers $l$ and $r$   in $D[0, T]$ so that   $\inf_{t\in[0, T]}[r(t)-l(t)] >0,$  and $l(0)\leq \psi(0) \leq r(0).$  Let $\omega(f, \delta, T)$ be as in \eqref{omega}. 
Then
\begin{equation}
 \omega(\eta^l, \delta, T) \leq  C[ \omega(\psi, \delta, T) + \omega(l, \delta, T)  ],  
 \label{7.291}
 \end{equation}
 and
 \begin{equation}
 \omega(\eta^r, \delta, T) \leq  C[ \omega(\psi, \delta, T) + \omega(r, \delta, T)  ],
 \label{7.292}
 \end{equation}
where  $C>0$ is a generic constant independent of  $\psi, l,  r$ and $[0,  T].$

\label{CorB.5}
\end{corollary}

\section{Constructing the solutions of the HJB} \label{sec:proofprop-R}

\subsection{Proof of Proposition \ref{prop-W}}\label{sec:zero-soln}
Our approach here is to find two  bounded solutions in the domains  $(-\infty, 0) $  and  $ (0, \infty) $  so that we can  paste them smoothly at the origin.
On $(0,  \infty),$   \eqref{4.18} can be written as 
$\frac{\sigma^2}{2} \W''(x) + (\beta-\delta_s x)\W'(x) -(\alpha + \delta_s) \W(x) + \theta_s =0.$ 
Since {\textcolor{red}{$W(x)\equiv T_s$}} is a particular solution, $\W$ has  the representation  $\W(x)=\U(x) + {\theta_s}/{(\alpha +\delta_s)},$  where $\U$ satisfies the homogeneous equation 
\begin{align}\label{homog-eq}
\frac{\sigma^2}{2} \U''(x) + (\beta-\delta_s x)\U'(x) -(\alpha + \delta_s) \U(x)=0, \quad x > 0.
\end{align}
Since we expect $\W$ to be bounded, we seek for a bounded non-trivial solution for $\U.$   A fundamental set  for this homogeneous equation consists of one bounded function and another unbounded function  as $x$ tends to infinity. We need only  to seek for this bounded solution.  Such a solution exists and has a stochastic representation (see Section 50 of Chapter V in \cite{R-W}): Consider a solution to the Ornstein-Uhlenbeck type equation $Z(t)=x +\sigma B(t) +  \int_{0}^{t}[\beta -  \delta_s Z(u)]du,$  where $x \geq 0,$   $B$  is a standard Brownian motion and  consider the stopping time  $\tau_0= \inf\{ t \geq 0: Z(t)=0 \}.$ Then  $\Psi_0(x)= E_x[e^{-(\alpha+ \delta_s)\tau_0}]$ for $x \geq 0$ is such a bounded solution to the homogeneous equation \eqref{homog-eq} and any other bounded solution in $(0, \infty)$  is a constant multiple of $\Psi_0.$  Moreover, $\Psi_0(0)=1$ and it  is a strictly decreasing function on $(0,  \infty).$ 
We expect $\W$ to be of the form $\W(x) = k_s \Psi_0(x) + {\theta_s}/{(\alpha +\delta_s)}$ in the interval $(0, \infty),$  where $k_s$ is a constant which needs to be determined. Since $\lim_{x \rightarrow \infty} \Psi_0(x)=0, $ it follows that  $\lim_{x \rightarrow \infty} \W(x)=   {\theta_s}/{(\alpha +\delta_s)}.$   

 We can perform a similar analysis on the interval $(-\infty,  0). $ Notice that  {\textcolor{red}{$W(x)\equiv -T_b$}} is a particular solution to \eqref{4.18} on $(-\infty, 0).$  Let $Z(t)=x +\sigma B(t) +  \int_{0}^{t}[\beta -  \delta_b Z(s)]ds,$ where $x \leq 0,$  and  $B$ is a standard Brownian motion. Introduce the stopping time  $\tau_0= \inf\{ t \geq 0: Z(t)=0 \}$ and let  $\Phi_0(x)=  E_x[e^{-(\alpha+ \delta_b)\tau_0}]$ for $x \leq 0.$  Then $\Phi_0$ is a bounded solution to the homogeneous  equation 
 \[ 
 \frac{\sigma^2}{2} \U''(x) + (\beta-\delta_b x)\U'(x) -(\alpha + \delta_b) \U(x)=0, \quad x <0.
 \]  
 Moreover,  $\lim_{x \rightarrow -\infty} \Phi_0(x)=0,$  $\Phi_0(0)=1,$ and it is strictly increasing in $(-\infty,  0)$ (for details, see Section 50 of Chapter V in \cite{R-W}).  Then similar to  above analysis, we expect $\W$ to be of the form  $\W(x) = k_b \Phi_0(x) - {\theta_b}/{(\alpha +\delta_b)}$ when  $x<0.$ The constant $k_b$  needs to be determined. Then $\lim_{x \rightarrow -\infty} \W(x)=  - {\theta_b}/{(\alpha +\delta_b)}$ also follows. 

Finally,  to determine the constants  $k_s$ and $k_b$, we impose the  ``smooth fit conditions''  for $\W$ across the origin as they are described by  $\W(0-)=\W(0+)$   and    $\W'(0-)=\W'(0+).$ We use  $\Psi_0(0)=\Phi_0(0)=1$ and   $\Psi_0'(0+) <0 < \Phi_0'(0-).$   Then, $\W(0-)=\W(0+)$ implies $k_b=k_s+T_s + T_b.$ The condition  $\W'(0-)=\W'(0+)$ implies $k_b=k_s {\Psi_0(0+)}/{\Phi_0(0-)}.$  By solving these  equations and  using  the condition $\Psi'_0(0+) <0 < \Phi'(0-),$   we obtain  
\[
k_s=-\frac{T_s + T_b }{1- \frac{\Psi'_0(0+)}{\Phi'_0(0-)}} <0, \quad \mbox{and} \quad k_b=k_s \frac{\Psi'_0(0+)}{\Phi'_0(0-)} >0.
\]
Consequently, for this pair of  $k_s$ and $k_b$,  the above described  $\W\in C^1(\rr),$  strictly increasing function on $\rr$  which satisfies  \eqref{4.18}. Moreover, $\W''$ is continuous  everywhere except at $x=0.$  It also satisfies \eqref{W-right} and \eqref{W-left} at the origin. Hence the proofs of parts(i), (ii) and (iii) are complete and it remains to show the uniqueness of $\W.$

To prove the uniqueness, we consider   $Z(t)=x +\sigma B(t) +  \int_{0}^{t}[\beta -  h(Z(s))]ds,$    where   $B$  is a standard Brownian motion. The condition $xh(x) \geq 0$ implies the explosion time of $Z$ is infinite. Since $\W\in C^1(\rr)\cap C^2(\rr-\{0\}),$   $\W''(0-)$ and $\W''(0+)$ are finite, we can apply  It$\hat{o}$'s lemma (see \cite{karatzas2}) to $\W(Z(T))e^{-\int_{0}^{T}(\alpha+\gamma(Z(t))dt }$ to obtain the representation
$\W(x)= E_x \left(\int_{0}^{\infty} e^{-\int_{0}^{s}(\alpha+\gamma(Z(t))dt } C'(Z(s)) ds\right).$  Hence the uniqueness of $\W$  follows.
\Halmos

\subsection{Proof of Proposition \ref{prop-R}}\label{sec:reflect-soln}
Let the cost parameters $p_b$ and  $p_s$  satisfy  
$ 0< p_b < T_b$      and    $   0< p_s < T_s. $
Throughout, we use  continuity properties of the solutions with respect to initial data and other parameters and we refer to Chapter V of \cite{hartman}.   First we gather a few useful facts about the differential equation  \eqref{4.18}. 
Recall that $W(x)=-T_b$ is a constant solution on the interval  $(-\infty,  0)$ and on $(0,  \infty),$ $W(x)=T_s$ is a constant solution for   \eqref{4.18}.  
\begin{lemma}
Let $\W$ satisfy \eqref{4.18} in the neighborhood of a point $x=c$  and assume $\W'(c)=0.$ Then the following hold.
\begin{itemize} 
\item[\rm (i)] If $c<0$ and $\W(c) >  {-\theta_b}/{(\alpha +\delta_b)},$ then $\W''(c)>0$ and  $x=c$ is a strict local minimum.
\item[\rm (ii)] If $c>0$ and $\W(c) <  {\theta_s}/{(\alpha +\delta_s)},$ then $\W''(c)<0$ and  $x=c$ is a  strict local maximum.
\end{itemize}
\label{lemA.1}
\end{lemma}
 
\proof 
If $c<0$  and $\W$ satisfies \eqref{4.18} together with $\W'(c)=0,$ then we observe that  
$\frac{\sigma^2}{2} \W''(c)=(\alpha+\delta_b)\left(\W(c)+  \frac{\theta_b}{\alpha +\delta_b}\right).$ 
Then the conclusion of part (i) is straightforward. The proof of part (ii) is similar and is omitted.
\Halmos

\begin{remark}\label{remA.1}
Notice that  if $c<0,$  $\W'(c) =0$  and $\W(c) =  {-\theta_b}/{(\alpha +\delta_b)},$ then by the uniqueness of solutions to the differential equation \eqref{4.18}, it follows that $\W(x) =  -T_b$ for all $x<0.$ Similar conclusion holds when    $c>0,$ $\W'(c) =0$  and $\W(c) = T_s.$ 

\end{remark}

For each  $a<0,$  we let  $\W_a$  be the solution to   \eqref{4.18} on $(-\infty,  0]$ with the initial data  $\W_a(a)=-p_b$ and $\W'_a(a)=0.$
Since $-T_b$ is a particular solution to   \eqref{4.18} on the interval $ (-\infty,  0],$  we can express $\W_a$  as follows:
\begin{equation}
\W_a(x) = U_a(x) -  T_b, \quad x\in (-\infty, 0)
\label{A.2} 
\end{equation}
where  $U_a$ is a solution to the homogeneous equation given below:
\begin{equation}
\frac{\sigma^2}{2} U''_a(x) + (\beta-\delta_b x) U'_a(x) -(\alpha +\delta_b)U_a(x) =0, \quad \text{  for  } x<0,   
\label{A.4}
\end{equation}
with the boundary conditions

\begin{equation}
U_a(a)\equiv d_b= T_b  - p_b >0,   \quad \text{  and }       U'_a(a)= 0.
\label{A.6}
\end{equation}

For the homogeneous differential equation  \eqref{A.4},  there exists  a fundamental set of solutions $\{ \Phi_0,  \Phi_1 \}$ on $( -\infty,  0]$ satisfying  $\Phi_0(-\infty)=0,$   $\Phi_1(-\infty)= \infty,$  $\Phi_0(0)=1,$   $\Phi_1(0)=1.$ These solutions  $\Phi_0$ and $\Phi_1$ have  stochastic representations  
as  described  in  Chapter V of \cite{R-W} (see Section 50 and in particular Proposition 50.3.).  Moreover, $\Phi_0$ is strictly increasing and $\Phi_1$ is strictly decreasing in the domain $(-\infty,  0).$

\begin{lemma}  Let $a<0$ then the following hold for each  solution $W_a.$
\begin{itemize}
\item[\rm (i)] Each $W_a$ is strictly decreasing on $(-\infty,  a),$ and  strictly increasing on  $(a, 0).$ Moreover, $W_a(0)$  and $W'_a(0)$  are finite.
\item[\rm (ii)] $W_a(-\infty)= \infty,$  and    $ \lim_{a \rightarrow -\infty} W_a(0)=\infty$.
\end{itemize}
\label{lemA.8}
\end{lemma}
\proof
By \eqref{A.2}, it suffices to prove the lemma for  $U_a.$ To show part (i), we observe that for any $a<0,$  if $U'_a (c)=0$ and $U_a(c)>0 $ for some $c<0,$ then  by \eqref{A.4}, $U''_a(c)>0$ and  $x=c$ is necessarily a strict local  minimum. In particular, for$a<0,$ $U_a(a)=d_b>0,$ and by \eqref{A.6},  $x=a$ is a strict local minimum. With the absence of any local maxima  when  $U_a(x)>0,$ we observe that $U_a$ is strictly  decreasing  in $ (-\infty, a)$ and strictly increasing in  $(a,  0).$  Moreover, the coefficients of the differential equation  \eqref{A.4} can be linearly and continuously  extended to $[0, \infty).$ Hence  $U_a(0)$  and $U'_a(0)$ are  finite.   This proves part (i).

For part (ii) let $a<0.$  Using the fundamental set we can write $U_a(x)= p_a \Phi_0(x) + q_a \Phi_1(x)$ for all $x \leq 0$, where $p_a$ and $q_a$ are constants. Since  $U_a$ has a positive, strict local minimum at $x=a,$   t $U_a$ is strictly decreasing in  $(a-\delta, a)$ and  strictly increasing  in $(a, a+\delta)$ for some $\delta>0. $  Since $\Phi_0$ is strictly increasing and $\Phi_1$ is strictly decreasing, it follows that  $p_a>0$ and $ q_a>0.$ But $\Phi_0$ is a bounded function on $(-\infty, 0)$ and $\Phi_1(-\infty )= \infty.$  Hence, $U_a(-\infty)=\infty$ and consequently, $W_a(-\infty)=\infty$ holds. 

To show  $ \lim_{a \rightarrow -\infty} U_a(0)=\infty,$ let $a_1 <a_2 <0.$  By \eqref{A.6},  $U_{a_1}$ and $U_{a_2}$ meet at a point in the interval $(a_1,  a_2)$ and $U_{a_1}(a_2) > U_{a_2}(a_2)=d_b.$ But by the uniqueness of solutions, two solutions to a homogeneous equation cannot meet more than once. Therefore,  $U_{a_1}(0) > U_{a_2}(0)>0.$ Consequently, $ \lim_{a \rightarrow -\infty} U_a(0)$ is finite or else it is $+\infty.$ Suppose $ \lim_{a \rightarrow -\infty} U_a(0)=L>0$ is finite. Then $U_a(0)=p_a +q_a \leq L$ for all $a<0.$
Moreover,   $p_a>0, q_a>0,$  and  $\Phi'_1(0) <0,$ $U'_a(0)\leq p_a \Phi'_0(0) \leq L \Phi'_0(0)$ for each $a<0.$
By integrating \eqref{A.4} and using  $U_a(a)=d_b>0,$ and $U_a$ is a non negative  increasing function in $(a, 0)$, we obtain the  bound $\frac{\sigma^2}{2} U'_a(0) \geq (\beta -\delta_b a)d_b -\beta U_a(0).$ Since $U_a(0)$ is bounded, RHS of this inequality tends to $+\infty$,  when $a$ tends to  $-\infty.$  But  $U'_a(0)$  is bounded above by     the constant $ L \Phi'_0(0)$ which is finite. This yields  a contradiction. Hence, $L= \infty$ and this completes the proof.
\Halmos

From the above discussion, $W_a$ satisfies the linear equation \eqref{4.18} on $(-\infty, 0),$ and $W_a(0)$ and $W_a'(0)$  are finite. For each $a<0,$  we can extend  $W_a$ to $(0,  \infty)$ so that it satisfies  \eqref{4.18}  on  $(0, \infty)$  with the available initial data {\textcolor{red}{for} }$W_a(0)$ and $W_a'(0).$ More precisely, on the interval  $(0, \infty),$ $W_a$ satisfies   
\begin{equation}
\frac{\sigma^2}{2} W''_a(x) + (\beta-\delta_s x) W'_a(x) -(\alpha +\delta_s)W_a(x) + \theta_s =0.
\label{A.11}
\end{equation} We let 
\begin{equation}
U_a(x) = W_a(x) - T_s, \quad x\in (0, \infty).
\label{A.12} 
\end{equation}
Then  $U_a$  satisfies the homogeneous equation 
\begin{align}\label{A.121}
\frac{\sigma^2}{2} U''_a(x) + (\beta-\delta_s x) U'_a(x) -(\alpha +\delta_s)U_a(x) =0, \quad x\in (0, \infty),
\end{align}
together with the boundary data  
$U_a(0)=W_a(0) - T_s,$ and $U_a'(0)=W_a'(0).$ 
Similar to $U_a$ on $(-\infty, 0)$, there is a fundamental set  of solutions  $\{ \tilde\Psi_0, \tilde\Psi_1 \}$ for the above homogeneous equation on $(0, \infty)$ so that  $\tilde\Psi_0(0)=\tilde\Psi_1(0)=1$, $ \tilde\Psi_0(\infty)=0$, $\tilde\Psi_1(\infty)=\infty,$   $\tilde\Psi_0$ is strictly decreasing and  $\tilde\Psi_1$ is strictly increasing on $(0, \infty).$ (See Chapter V of \cite{R-W}  for a stochastic representation of  $\tilde\Psi_0$ and  $\tilde\Psi_1.$)  Hence each $U_a$ can be written as 
\begin{equation}
U_a(x) = r_a \tilde\Psi_0(x) +t_a \tilde\Psi_1(x), \quad x\in (0,\infty),
\label{A.14} 
\end{equation}
where  $r_a$ and $t_a$ are constants. 

In the next two results,  we show that when $a<0,$ there are two types of solution profiles of $W_a$  which exhibits very different behavior on $[0, \infty).$   When  $|a| $  is very large,  $W_a(x)$ tends to  $+\infty$ as $x \rightarrow \infty.$   When $|a|$ is small,  $W_a(x)$ tends to  $-\infty$ as $x \rightarrow \infty.$   These two  profiles are separated by the solution  curve of $W_c$, where  $a=c<0$ is a special point.  The solution $W_c$ is bounded on $[0, \infty)$  and it approaches   $T_b$  as   $x \rightarrow \infty.$

\begin{lemma}There exists a $\delta >0$ so that if  $-\delta < a <0,$ then  
\begin{itemize}
\item[\rm (i)] $W_a$ has a local maximum  in $(0, \infty)$;
\item[\rm (ii)]     $ \lim_{x \rightarrow \infty} W_a(x)=-\infty.$  
\end{itemize}
\label{lemA.16}
\end{lemma}
\proof
We first consider $W_0$ which corresponds to $a=0$. Then $W_0$ satisfies \eqref{4.18}  everywhere except at $x=0$  and $W_0(0)=-p_b,  \text{  and  } W'_0(0)= 0.$ Moreover, $W_0$ and $W'_0$ are continuous everywhere,   $W''_0$ has a jump discontinuity at $x=0,$
$\frac{\sigma^2}{2} W''_0(0-)=(\alpha +\delta_b)\left(  T_b  -p_b\right) >0,$
and   
$\frac{\sigma^2}{2} W''_0(0+)=-(\alpha +\delta_s)\left(  T_s  +p_b\right) <0.$
Following the proof of  Lemma \ref{lemA.1} (i), $W_0$ is strictly decreasing on $(-\infty,  0).$ On $(0, \infty)$, since $W'_0(0)=0$ and  $W_0''(0)<0,$ by a similar argument,  $W_0$ is strictly decreasing. 

We can then choose $x_0>0 $ and $\epsilon_0>0$ so that $W_0(x_0)<-p_b-2\epsilon_0.$ Note that both $W_0$ and $W_a$ satisfy \eqref{4.18}  on $(-\infty, 0).$    Hence, using the continuity of the solutions in initial data, for any $\epsilon >0,$ there exists a $\delta_1 >0$ so that  $|W_a(0)-W_0(0)| + |W'_a(0)-W'_0(0)| < \epsilon$ whenever $-\delta_1< a<0.$  The functions  $W_0$ and $W_a$ satisfy \eqref{A.11} in $(0, \infty).$  Again using the same continuity property, we can find a $\delta_2 >0$ so that $ |W_a(x_0)-W_0(x_0)|<\epsilon_0$  whenever $|W_a(0)-W_0(0)| + |W'_a(0)-W'_0(0)| < \delta_2.$ Now, by letting  $\epsilon=\delta_2,$ we can find a $\delta>0$ so that  $ W_a(x_0)<-p_b-\epsilon_0$ whenever  $-\delta <a <0.$ Hence, each $W_a$ is strictly increasing in $(a, 0)$, $W_a(a)=-p_b$ and $W_a(x_0) < -p_b-\epsilon_0$. Thus,  $W_a$ has a local maximum  in $(0, x_0).$  This completes part (i).\\
 From part (i), $W_a$ has a local maximum at a point $x=c$ in $(0, x_0).$ Thus,  $W'_a(c)=0,$ 
$W''_a(c)= \frac{2(\alpha +\delta_s)}{\sigma^2} \left(W_a(c) - T_s\right) \leq 0,$
which yields that $W_a(c)  \leq T_s.$ But $W_a(c)$ cannot be equal to $T_s$; otherwise $W_a(x) = T_s$ becomes the unique solution as explained in Remark  \ref{remA.1}, which contradicts with the fact that $W_a$ is strictly decreasing on $(0,\infty)$. Hence $ W_a(c)  < T_s$ and $W_a'(c)<0 $   for $x>c.$  Using \eqref{A.11}, then it follows that $W_a$ is strictly concave and strictly decreasing for large $x>0.$ Hence, part (ii) follows. \Halmos

Next we consider the properties of the solution $W_a$ when $a<0$ and $|a|$ is very large. From Lemma \ref{lemA.8} (ii), we can find a constant $K>0$ so that 
\begin{align}
\label{large_a}
\mbox{when  $a<-K,$} \quad W_a(0)>  T_s.
\end{align}

\begin{lemma} Let $a<-K$ where  $K>0$ is given in \eqref{large_a}. Then 
\begin{itemize}
\item[\rm (i)] $W_a$ is strictly increasing on the interval  $(a, \infty)$;
\item[\rm (ii)] $ \lim_{x \rightarrow \infty} W_a(x)=\infty.$  
\end{itemize}
\label{lemA.20}
\end{lemma}
\proof
By Lemma \ref{lemA.8} (i), $U'_a(0) \geq 0$ (equivalently, $W'_a(0) \geq 0$). The differential equation \eqref{A.4} can be continuously and linearly extended to $(0, \infty).$  Hence we consider the corresponding extended solution $U_a$ as given in \eqref{A.2} and \eqref{A.12}. Next, we show that $U'_a(0) >0.$ Suppose $U'_a(0)=0.$ Then  by \eqref{large_a}  and \eqref{A.121}, we have 
$U''_a(0) = \frac{2}{\sigma^2} (\alpha + \delta_s) U_a(0) > 0.$  Thus $U_a$ has a strict local minimum at the origin. This is a contradiction  since $U_a$ is strictly increasing in $(a, 0).$ Hence $U'_a(0)>0$ and consequently,  $W'_a(0)>0.$ Thus, we can find a $\delta>0$ so that  $W'_a(x)>0$ when  $0\leq x< \delta.$  Let $c=\sup\{x>0: W'_a(u)>0 \text{ for } 0<u<x \}.$
Suppose $c$ is finite. Then $W'_a(c)=0,$ and  $ W_a(c) >  {\theta_s}/{(\alpha +\delta_s)}.$   Hence, $ W''_a(c)>0$   and $x=c$ is a strict local minimum. This is a contradiction since $W'_a(x)>0$ in $(0, c).$  Hence $c=\infty $ and part (i) follows.

By part (i) and  \eqref{A.12}, $U_a$ is strictly increasing on $(a, \infty).$ Then $ \lim_{x \rightarrow \infty} U_a(x)$  exists and suppose it  is finite.  Then,  $t_a=0$  and $r_a<0$ in \eqref{A.14}.  Thus, $ \lim_{x \rightarrow \infty} U_a(x)=0$  and by \eqref{A.12},  $ \lim_{x \rightarrow \infty} W_a(x)= T_s.$ It contradicts with \eqref{large_a} since $W_a$ is increasing on $(a, \infty)$. Hence, $ \lim_{x \rightarrow \infty} U_a(x)= \infty$  and by  \eqref{A.12},  part (ii) follows.
\Halmos

\begin{lemma}
There exists a point $c<0$ so that
\begin{itemize}
\item[\rm (i)] the solution $W_c$  is strictly increasing   in $(c,  \infty)$;   
\item[\rm (ii)]    $W_c$ is bounded   in $(c, \infty)$  and $ \lim_{x \rightarrow \infty} W_c(x)= T_s.$  
\end{itemize}
\label{lemA.22}
\end{lemma}
\proof
Let  $c=\sup\{a<0:    \lim_{x \rightarrow \infty} W_a(x)=\infty \}$. 
By Lemmas \ref{lemA.16} and \ref{lemA.20}, such a point $c$ exists and $-\infty<c<0.$   Moreover,  $W_c$ is strictly increasing on $(c, 0)$  and thus $W_c'(0) \geq 0 .$ Suppose  $W_c'(0)=0,$  then $W''(0-)>0$ from \eqref{A.11} and by the fact that $W_c(0) >-p_b> -T_b. $  This is a contradiction since $W_c$ is strictly increasing in $(c, 0).$  Hence $W'_c(0) >0$ and there is a $\delta>0$ so that $W'_c >0 $ over $(c, \delta).$ We let $x_0=\sup\{x>0 : W'_c (u)>0 \text{ for all }  0<u< x \}.$  Clearly, $x_0>\delta.$ Suppose  that $x_0$ is finite.  Then  $W'_c(x_0)=0.$  If $W_c(x_0)>  {\theta_s}/{(\alpha +\delta_s)},$ using  \eqref{A.11},  $W''_c(x_0)>0.$ This is a contradiction since $W_c$ is strictly increasing in $(0, x_0).$ Thus $W_c(x_0) \le  T_s.$ Now if $W_c(x_0)=  T_s,$ by the uniqueness of solutions, $W_c$ is a constant on $(0, \infty) $ and this contradicts with $W'_c(0) >0.$ Hence the only possibility is $W_c(x_0) <  T_s,$ in which case $W''_c(x_0) <0 $ and $x=x_0$ is a strict local maximum and  $W_c$ is strictly decreasing on $(x_0, \infty).$ Similar to the proof of  Lemma \ref{lemA.16} (ii), we see that  $ \lim_{x \rightarrow \infty} W_c(x)=-\infty.$ Now using the continuity of the solutions $W_a$ with respect to the parameter $a$, we can find a $\delta_2>0$ so that for each $c-\delta_2 <a<c,$  $ \lim_{x \rightarrow \infty} W_a(x)=-\infty.$ This contradicts with the definition of  $c.$ Hence $x_0$ is infinite and part (i) follows.

From part (i), $ \lim_{x \rightarrow \infty} W_c(x)=L \leq \infty$ exists. Suppose  $L=\infty.$ We can pick  $M>0$ so that $W_c(M) >1+T_s.$  Using the continuity of $W_a$ and $W'_a$ with respect to $a$, we can find a $\delta_3>0$ so that for each $c<a<c+ \delta_3,$  $W_a(M) > T_s,$ and $W_a'(x)>0$ on $(0, M).$  By the proof in Lemma \ref{lemA.20}, each such $W_a$  increases to $\infty$ and  this contradicts with the definition of $c$.  Hence $L <\infty$  and $W_c$ is bounded in $[c, \infty)$  and equivalently, $U_c$ defined in \eqref{A.12} is also bounded. By \eqref{A.14}, we can write  $U_c(x) = r_c \tilde\Psi_0(x) +t_c \tilde\Psi_1(x)$ for $x>0.$ Since $U_c$ is bounded, we have $t_c=0$ and thus,  $U_c(x) = r_c \tilde\Psi_0(x)$ for $x>0.$ But  $ \lim_{x \rightarrow \infty} \tilde\Psi(x)=0$ and consequently $ \lim_{x \rightarrow \infty} U_c(x)=0.$ Using this together with  \eqref{A.12}, we obtain $ \lim_{x \rightarrow \infty} W_c(x)= T_s.$   This completes the proof.
\Halmos

In the next lemma, we establish the uniqueness of the point $c.$
\begin{lemma}
For $a<0$, let $W_a$ be the corresponding solution of  \eqref{4.18} so that   $ \lim_{x \rightarrow \infty} W_a(x) $ is finite.  Then $a=c$ and  $W_a(x)\equiv W_c(x)$ for all $x.$
\label{lemA.24}
\end{lemma}

\proof
Suppose $a\neq c$  and $W_a$ satisfies the assertion.  Without loss of generality,  let $c<a.$ We note that $W_a(a) = W_c(c) = -p_b$ and  from Lemma \ref{lemA.8}, $W_a$ and $W_c$ are strictly decreasing on $(-\infty, a)$ and $(-\infty, c)$ respectively, and are strictly increasing on $(a, 0)$ and $(c, 0)$ respectively. Hence $W_c$ and $W_a$ meet at a point  $x_0,$  where $c<x_0 < a$, and moreover, $W_c(a)>W_a(a)=-p_b.$ Next, following the argument of Lemma \ref{lemA.22} (ii), we conclude that $ \lim_{x \rightarrow \infty} W_a(x)=  \lim_{x \rightarrow \infty} W_c(x)=T_s.$ 

We now consider  $H(x)=W_c(x)-W_a(x)$ for  $ x_0 \leq x <\infty.$  Then $H(x_0)=0,$  $H(\infty)=0$ and $H(a)>0.$  Moreover, $H$ is a bounded solution to the homogeneous equation   $\mathcal{G}H(x) - \gamma(x) H(x)=0.$ We let $m > x_0$ so that  $H(m)=\max_{x\in[x_0, \infty)} H(x).$  Then $H(m)>0$ and $H'(m)=0.$
When $m \neq 0,$  using the differential equation for $H$ at $x=m$, we obtain $H''(m)>0$ and thus, $x=m$ is a local minimum, which contradicts with the definition of $H(m)$. Hence $m=0.$  Again, by the differential equation, we obtain $H''(0-) >0$ and thus $H$ is strictly decreasing  in  $(-\delta, 0)$ for some $\delta>0,$  which is again a contradiction because $H(m)$ is the maximum. Hence $c=a$ and the conclusion of the lemma follows. \Halmos

Let  $c=\sup\{a<0:    \lim_{x \rightarrow \infty} W_a(x)=\infty \}$ be as in Lemma \ref{lemA.22}. 

\begin{lemma}
For $a \in (c, 0),$  $W_a(x) <W_c(x)$  for all $ x>a,$  $W_a$ has a unique local maximum on $(0,  \infty)$ and thereafter it is decreasing to $-\infty$, i.e., $ \lim_{x \rightarrow \infty} W_a(x)= -\infty.$  Moreover, if $c<a_1<a_2<0,$ then  $W_{a_2}(x)<W_{a_1}(x)<W_c(x)$  for all $x>a_2.$
\label{lemA.26}
\end{lemma}

\proof
Similar to the argument in the proof of Lemma \ref{lemA.24}, $W_c$ and $W_a$ meet at a point  $x_0 \in (c, a)$, and $W_c(a) > W_a(a)$. Suppose that $W_c$ and $W_a$ meet at a point $z>a.$ We consider $H(x)=W_c(x)-W_a(x)$ for  $ x_0 \leq x \leq z$, and define $m = \arg\max_{x\in [x_0, z]} H(x)$. Following an argument similar to the proof of Lemma \ref{lemA.24}, we can obtain a contradiction on $H(m)$, and thus $W_a(x) <W_c(x)$ for all $x>a.$  Consequently,  $W_a(x) < T_s $ for all $x>a$. By Lemma \ref{lemA.1},  $W_a$ can have at most one local maximum and no local minima. If $W_a$  has no local maxima, then it is increasing and bounded above by ${\theta_s}/{(\alpha +\delta_s)}$, and $ \lim_{x \rightarrow \infty} W_a(x) $ is finite. Using   Lemma \ref{lemA.24}, we have $a=c$ and this is a contradiction. Therefore, each $W_a$ has exactly one local maximum in $(0, \infty)$ and thereafter  $W_a$ is decreasing. Again using Lemma \ref{lemA.24}, we must have $ \lim_{x \rightarrow \infty} W_a(x)= -\infty.$  The proof of the last assertion is very similar to the above proof of $W_a(x) <W_c(x)$  for all $ x>a.$ And it is omitted. This completes the proof. 
\Halmos

For $a \in (c, 0)$, assume that each $W_a$ achieves the unique local maximum at $x=r_a.$  We introduce the function $M:(c, 0) \rightarrow  (-p_b,  T_s)$ by   

\begin{equation}
M(a)\equiv \max\limits_{x\in [a, \infty)} W_a(x) = W_a(r_a).
\label{max}
\end{equation}

\begin{proposition}  The following results hold for the function $M$ defined in \eqref{max}. 
\begin{itemize}
\item[\rm (i)] The function $M$ is a continuous strictly decreasing function. 
\item[\rm (ii)]  $ \lim_{a \rightarrow c+} M(a)= T_s$  and   $ \lim_{a \rightarrow 0-} M(a)=-p_b$ and $M$ is an on-to function.
\item[\rm (iii)]  Given  $0<p_s < T_s,$  there is a unique $a^*$ in $(c, 0)$ so that  $M(a^*)=p_s.$ Let  $b^*\equiv r_{a^*}>0.$  Then the  corresponding $W_{a^*} $  is strictly increasing  in $(a^*, b^*)$   and $W_{a^*}(b^*)=M(a^*)=p_s >0$  and the triple  $(W_{a^*}, a^*, b^*)$ satisfies all the conditions in Proposition \ref{prop-R}.
\end{itemize}
\label{prop-max}
\end{proposition}
\proof
The strictly decreasing property of the function $M$ follows from Lemma \ref{lemA.26}. To show the continuity of $M,$ we fix  $a_0$ so that $c<a_0<0$. 
Since the solution to \eqref{4.18} is continuous  with respect to the initial data, for 
a given $\epsilon_1>0,$ there exists a $\delta>0$ so that  $a \in (a_0 -\delta, a_0 + \delta )\subset (c, 0)$ implies  
$|W_a(0)-W_{a_0}(0)|+|W'_a(0)-W'_{a_0}(0)|< \epsilon_1.$ 
Let $b\in (c, a_0-\delta)$, and choose $z_0>0$ so that $W_b(z_0)<0.$ Next, $W_a$ satisfies \eqref{A.11} on $(0, \infty)$ and the coefficients of this differential equation can be smoothly extended to $(c, \infty).$ Using the continuity property of initial data for this extended linear equation, for any $\epsilon > 0$, we can find  $\delta_2>0$ so that  
$|W_a(0)-W_{a_0}(0)|+|W'_a(0)-W'_{a_0}(0)|< \delta_2$ guarantees 
$\sup_{x\in [0, z_0]} |W_a(x)-W_{a_0}(x)|<\epsilon.$   
Now choose $\epsilon_1=\delta_2$.  Then $|a_0- a|<\delta$ implies  
$\sup_{x\in[0, z_0]} |W_a(x)-W_{a_0}(x)|<\epsilon. $ 
By Lemma \ref{lemA.26}, $W_b(z_0)<0$ implies $W_a(z)<0$ for all $a \in (b, 0)$ and $z \ge z_0$. Consequently, when  $a-\delta<a_0 < a$, from Lemma \ref{lemA.26}, 
$0 < W_a(r_{a_0}) \le W_a(r_a) < W_{a_0}(r_{a}) \le W_{a_0}(r_{a_0}).$
Thus, 
$|M(a)-M(a_0)| = W_{a_0}(r_{a_0}) - W_a(r_{a_0}) \le \sup_{x\in [0, z_0]} |W_a(x) - W_{a_0}(x)| < \epsilon.$
Similarly, when $a_0-\delta< a <a_0,$  it  follows that $|M(a)-M(a_0)|  \le \sup_{x\in [0, z_0]} |W_a(x) - W_{a_0}(x)| < \epsilon.$ This proves the continuity of $M$ at $a_0.$  The proof of part (i) is complete.

Since $M$ is strictly decreasing and $T_s$ is an upper bound,  $ \lim_{a \rightarrow c+} M(a)$ exists and it is less than or equal to  $T_s. $  Again by the continuity property of the initial data, $W_c$ is close to $W_a$  when $a$ is very close to $c.$  Using this together with Lemma \ref{lemA.22} (ii),  it follows that   $ \lim_{a \rightarrow c+} M(a)= T_s. $ The proof of  $ \lim_{a \rightarrow 0-} M(a)=-p_b$ is similar.  These limits together with part (i) implies that $M$ is an onto function.

Finally, part (iii) is a consequence of parts (i) and (ii). This completes the proof. 
\Halmos

\begin{remark}\label{rem:thm4.6}
For the proof of Theorem \ref{thm:solution-DCP} (iii), let $(c, W_c)$ be as in the Lemma \ref{lemA.22}. Take $\tilde a^*=c,$ and  let $W_L$ identically equal to $W_c.$  Then the pair  $(\tilde a^*,  W_L)$ satisfies  \eqref{4.40}. 
\end{remark}

\end{appendix}

{\bf{Acknowledgements:}}   Ananda Weerasinghe  is partially supported by the Simons Foundation grant 317381 and  grateful for their support.

\vskip 1cm
\begin{tabular}{lcl}
Xin Liu& \hspace{0.5cm} & Ananda P. Weerasinghe\\
O-324 Martin Hall&& 396 Carver Hall  \\
 School of Mathematical and Statistical Sciences&&Department of Mathematics\\
 Clemson University&&Iowa State University\\
Clemson, SC  29634,  USA&&  Ames, IA 50011, USA.\\
\texttt{xliu9@clemson.edu}&& \texttt{ananda@iastate.edu} 
\end{tabular}

\begin{thebibliography}{10}
 
 
\bibitem{aaa}
{\sc  Akan, M., Alagoz, O., Ata, B., Erenay, F.S. and Said, A.}, {\em A broader view of designing the liver allocation system,} Operations research, 60 (2012), pp.757-770.


\bibitem{amr}
{\sc  Atar, R., Mandelbaum, A. and  Reiman, M. I.}, {\em Scheduling
a multi-class queue with many exponential servers:Asymptotic
optimality in heavy traffic ,} Ann. Appl. Probab., 14
 (2004), pp. 1084-1134.
 
 

\bibitem{billingsley}
{\sc Billingsley, P.}, {\em  Convergence of Probability Measures}, second edition, John Wiley and sons, New York (1999).



\bibitem{boxma}
{\sc Boxma, O. J., Israel, D., Perry, D.  and Stadje, W.},  {\em A new look at organ transplantation models and double matching queues}, Probability in the engineering and information sciences 25 (2011), pp 135--155.



\bibitem{burdzy}
{\sc Burdzy, K., Kang, W. and Ramanan, K.}, {\em The Skorokhod
problem in a time-dependent interval}, Stoch. Proc. and their  Appl., 119 (2009), pp. 428-452.


\bibitem{conolly}
 {\sc Conolly, B. W., Parthasarathy, P. R., and  Selvaraju, N.,}, {\em Double-ended queues with impatience,}
 Computers \& Oper. Res., 29  (2002),  pp. 2053--2072.

\bibitem{dai-dai}
 {\sc Dai, J. G. and  Dai, W.}, {\em  A heavy traffic limit theorem for a class of open queueing networks with finite buffers,}
Queueing Syst., 32 (1999),  pp. 5--40.

\bibitem{dai}
 {\sc Dai, J. G. and  He, S.}, {\em Customer abandonment in many-server queues,}
 Math. Oper. Res., 35 (2010),  pp. 347-362.

\bibitem{kurtz}
{\sc Ethier, S. N. and Kurtz, T. G.}, {\em Markov  Processes:
Characterization and Convergence}, Wiley, New York (1986).

\bibitem{fleming}
{\sc Fleming, W. H. and Soner, H. M.}, {\em Controlled Markov
Processes and Viscosity Solutions,} Second edition, Springer-Verlag,
New York (2006).


\bibitem{GW}
{\sc Gurvich, I., and  Ward, A.,} {\em On the dynamic control of matching queues}, Stoch. Syst.,
 4 (2014), pp. 479--523.


\bibitem{halfin}
{\sc Halfin, S. and  Whitt, W.}, {\em  Heavy traffic limits for
queues with many exponential servers}, Oper. Res. 29 (1981), pp.
567-588.

\bibitem{harrison1}
{\sc Harrison, J. M.}, {\em Assembly-like queues, } Journal of Appl. Probab. 10 (1973), pp. 354-367.

\bibitem{harrison}
{\sc Harrison, J. M.}, {\em Brownian Motion and Stochastic Flow
Systems}, Wiley Publications, New York (1985).

\bibitem{hartman}
{\sc Hartman, P.}, {\em Ordinary Differential Equations}, Wiley
Publications, New York (1964).

\bibitem{jacod}
{\sc Jacod, J., and  Shiryaev, A.}  {\em Limit theorems for stochastic processes}. Springer Science \& Business Media (2003).



\bibitem{karatzas2}
{\sc Karatzas, I. and Shreve, S.E.}, {\em Brownian Motion and
Stochastic Calculus}, Springer-Verlag, New York (1988).


\bibitem{kashyap}
{\sc Kashyap, B. R. K.}, {\em The double-ended queue with bulk service and limited waiting space}, Oper. Res. 14  (1966), pp. 822-834.

\bibitem{perry0}
{\sc Kaspi, H. and Perry, D.},
{\em Inventory systems of perishable commodities,}
 Adv. Appl. Prob. 15  (1983),  pp. 674--685.
 
 \bibitem{kl}
{\sc Khademi, A. and Liu, X.},
{\em Asymptotically optimal allocation policies for transplant queueing systems,}
Working Paper (2020).
 
\bibitem{kocaga}
{\sc Kocaga, Y. L. and  Ward, A.,} {\em Admission control for a multi-server queue with abandonment}, Queueing Systems,
 65 (2010), pp. 275-323.



\bibitem{taksar}
{\sc Krichagina, E. V. and Taksar, M.}, {\em Diffusion approximation
for GI/G/1 controlled queues}, Queueing Systems Theory  Appl., 12
(1992), pp. 333-367.

\bibitem{kruk}
{\sc Kruk, L., Lehoczky, J., Ramanan,K. and Shreve, S.}, {\em An
explicit formula for the Skorokhod map on [0,a]}, Ann. Appl.
Probab.,  (2007), pp. 669-682.

\bibitem{lllz}
{\sc Lee, C., Liu, X., Liu, Y. and Zhang, L.}, {\em Optimal control of a time-varying double-ended production queueing model}, working paper (2020).

\bibitem{L}
{\sc Liu, X.,} {\em Diffusion approximations for double-ended queues with reneging in heavy traffic,} Queueing  Syst.,  (2018),   pp. 1--39



\bibitem{LGK}
{\sc Liu, X.,  Gong, Q.,  and  Kulkarni,  V. G.,} {\em Diffusion models for double-ended queues with renewal arrival processes,} Stoch. Syst.,  5 (2015), pp. 1--61.


\bibitem{meyer}
{\sc Meyer, P. A.}, {\em Un cours sur les integrales  stochastiques,
Seminaire de  Probabilities X, Lecture Notes in Math. 511},
Springer, New York (1974).

\bibitem{ptw}
{\sc Pang, G., Talreja, R. and Whitt, W.}, {\em Martingale proofs of
many-server heavy-traffic limits for Markovian queues,} Probability
Surveys, 4 (2007), pp. 193-267.



\bibitem{perry}
{\sc Perry, D. and  Stadje, W.,}
{\em Perishable inventory systems  with  impatient demands,}
 Math. Methods of Oper. Res., 50  (1999),  pp. 77--90.
 

\bibitem{plumbeck}
{\sc Plumbeck, E. L., and  Ward, A. R.},
{\em Optimal control of a high-volume assemble-to-order system,}
 Math. Oper. Res., 31  (2006),  pp. 453--477.

\bibitem{prabhakar}
 {\sc Prabhakar, B., Bambos, N.,   and  Mountford, T. S.,}, {\em The synchronization of Poisson processes and  queueing 
networks with service and synchronization nodes,} Adv.  Appl. Probab., 32 (2000), pp. 824--843.

\bibitem{protter}
 {\sc Protter, P.}, {\em Stochastic Differential Equations,} Second edition, Springer-Verlag (2004).
 


\bibitem{reed2}
{\sc Reed, J. E.},
{\em  The  $G/GI/N$  queue  in the Halfin-Whitt  regime,}
Ann.  Appl. Probab.  19  (2009),  pp. 2211--2269.

\bibitem{reed}
{\sc Reed, J. E. and  Ward, A. R.},
{\em Approximating the {$GI/GI/1+GI$} queue with a nonlinear drift
  diffusion: hazard rate scaling in heavy traffic,}
 Math. Oper. Res., 33  (2008),  pp. 606--644.

\bibitem{R-W}
{\sc Rogers, L. C. G. and  Williams,  D.},
{\em Diffusions,  Markov Processes and  Martingales,}
Volume 2,  Cambridge University Press, (2000).

\bibitem{slaby}
{\sc Slaby, M},
{\em Explicit representation of the Skorokhod map with time dependent boundaries},
Probability and Mathematical Statistics 30 (2010): pp. 29-60.

\bibitem{song}
{\sc  Song,  J. S., and Zipkin, P.},
{\em Supply chain operations: Assemble-to-order  and configure-to-order systems,}
Handbook in Oper. Res. and Mgmt. Sci., Vol. XXX  (2003),  pp. 561--593.


\bibitem{talreja}
{\sc  Talreja, R.  and Whitt,  W.},
{\em Heavy-traffic limits for waiting times in many-server queues with abandonment,}
Ann.  Appl. Probab.  19  (2009),  pp. 2137--2175.


\bibitem{ward-kumar}
{\sc   Ward, A. R. and Kumar, S.},
{\em Asymptotically Optimal Admission Control of a Queue with Impatient Customers,}
 Math. Oper. Res., 33  (2008),  pp. 167--202.





 \bibitem{weera2}
 {\sc Weerasinghe, A.}, {\em A bounded variation control problem for diffusion processes},
 SIAM J. Control and Opt., 44 (2005), pp. 389-417.
 
 \bibitem{weera4}
 {\sc Weerasinghe, A.  and  Mandelbaum,  A.}, {\em Abandonment versus blocking in many server queues:asymptotic optimality  in the QED regime},
Queueing Syst., 75 (2013), pp. 279--337.
 
  
 
 \bibitem{weera3}
  {\sc Weerasinghe, A.}, {\em Diffusion approximations for $G/M/n+GI$ queues with state dependent service rates},
 Math. Oper. Res., 39 (2014), pp. 207--228.
 

\bibitem{whitt}
{\sc   Whitt,  W.},
{\em Stochastic-Process limits,}
Springer Series in Operations Research, Springer-Verlag, New York. (2002).

 
 
\bibitem{whitt2}
{\sc   Whitt,  W.},
{\em Proof of the martingale FCLT,}
Probability Surveys, 4 (2007), pp 268--302.

 
 
 
 
 
 
 
\end{thebibliography}
\end{document}